\ifdefined\pdfinfoomitdate\pdfinfoomitdate=1\fi
\ifdefined\pdftrailerid\pdftrailerid{}\fi
\ifdefined\pdfsuppressptexinfo\pdfsuppressptexinfo=15\fi
\documentclass[a4paper,11pt,onecolumn]{amsart}
\usepackage[a4paper,margin=1.02in,includeheadfoot]{geometry}
\usepackage[T1]{fontenc}
\usepackage{lmodern,mathtools,amssymb,graphicx,booktabs,array,tikz}
\usepackage[expansion=false]{microtype}
\usepackage[colorlinks=true,allcolors=blue]{hyperref}
\hypersetup{pdftitle={Whitehead doubles and the failure of the AJ conjecture},
  pdfsubject={Both clasp signs, both torus-knot chiralities, and all integer twists}}
\newtheorem{theorem}{Theorem}[section]
\newtheorem{corollary}[theorem]{Corollary}
\newtheorem{lemma}[theorem]{Lemma}
\newtheorem{proposition}[theorem]{Proposition}
\newtheorem{conjecture}[theorem]{Conjecture}
\theoremstyle{definition}
\newtheorem{definition}[theorem]{Definition}
\theoremstyle{remark}

\newcommand{\Q}{\mathbb Q}
\newcommand{\C}{\mathbb C}
\newcommand{\Z}{\mathbb Z}
\newcommand{\D}{\mathcal D}
\newcommand{\cat}{G_{\mathrm{Cat}}}
\DeclareMathOperator{\ord}{ord}
\DeclareMathOperator{\Res}{Res}
\DeclareMathOperator{\Vol}{Vol}
\newcommand{\tr}{\operatorname{tr}}
\newcommand{\sfL}{\operatorname{sf}_L}
\newcommand{\SL}{\mathrm{SL}_2}
\newcommand{\diag}{\operatorname{diag}}
\newcommand{\qbinom}[2]{\genfrac{[}{]}{0pt}{}{#1}{#2}_q}
\newcommand{\imunit}{\mathrm i}
\newcommand{\torcolor}{\mathcal T}

\newcommand{\cheb}{\mathsf U}
\makeatletter
\renewcommand{\section}{\@startsection{section}{1}%
  \z@{1.05\linespacing\@plus\linespacing}{.5\linespacing}%
  {\normalfont\scshape\centering}}
\makeatother

\allowdisplaybreaks
\title[Whitehead doubles and the AJ conjecture]{Whitehead doubles and the failure of the AJ conjecture}
\author{Jinzheng Li}
\address{Particle Theory and Cosmology Group, Center for Theoretical Physics of the Universe, Institute for Basic Science, Daejeon 34126, Republic of Korea}
\email{jinzhengli@ibs.re.kr}

\author{Yuqi Li}
\address{C. N. Yang Institute for Theoretical Physics, Stony Brook University, Stony Brook, New York 11794, USA}
\email{yuqi.li@stonybrook.edu}

\author{Hao-Yu Sun}
\address{Weinberg Institute, The University of Texas at Austin, 2515 Speedway, C1600, Austin, Texas 78712, USA}
\address{SpaceX, USA}
\email{hkdavidsun@utexas.edu}
\date{November 18, 2026}
\begin{document}
\begin{abstract}
For marked Whitehead doubles of torus knots with coprime parameters \(2\le r<s\), the AJ conjecture fails for both clasp signs, both companion chiralities, and every integer twist. The zero-framed signed unreduced colored Jones sequence has a nonzero rationally saturated annihilator. Every nonzero primitive representative specializes nontrivially at \(t=-1\) with the meridian variable \(M\) independent. With \(L\) the longitude variable, the specialization contains the irreducible quadratic \(E=L^2-(M^4-4M^2+1)L+M^4\) for the negative clasp and \(M^4E(M^{-1},L)\) for the positive clasp. The corresponding quadratic is coprime to the full classical \(A\)-polynomial, including \(L-1\), over \(\Q(M)[L^{\pm1}]\), so the discrepancy persists when each irreducible factor is counted once. The curve \(E=0\) comes from the marked Whitehead link exterior with its axis cusp complete. For the displayed monic polynomial, \(\pi m(E)\) is the unsigned volume of that two-cusped exterior. We use the actual Jones sequence to select a right divisor of the satellite transfer operator. Rational factor tests show that the quadratic survives in this divisor, and the two marked companion fillings exclude it as a component of the satellite boundary curve.
\end{abstract}
\maketitle
\section{Introduction}\label{sec:introduction}

Whitehead doubles of nontrivial torus knots give a uniform family of
counterexamples to the AJ conjecture. For either clasp, either companion chirality,
and every integer twist, the specialized minimal recurrence has an
irreducible quadratic missing from the classical boundary curve.
For the negative clasp this is the outer boundary curve in the marked
Whitehead link shape chart with its axis cusp complete. The axis
boundary torus is where the companion exterior attaches to form the
satellite. The problem is to show that this curve
survives on the quantum side after passing to the actual minimal
recurrence, and to see why the gluing excludes it as a classical
component.

The AJ conjecture of Garoufalidis~\cite{S01} compares the colored
Jones recurrence with the classical \(A\)-polynomial.
There are substantial families of AJ-positive satellites, including
the cable families of Ruppe--Zhang, Ruppe, and Tran under the
restrictions discussed in Section~\ref{sec:fingerprints}
\cite{RuppeZhang2015,Ruppe2015,TranCable2017}.
Counting each irreducible factor once addresses multiplicity discrepancies
such as those in Zhang's connected-sum examples~\cite{S06}.
The Whitehead quadratic survives this reduction.

\subsection{The recurrence and boundary comparison}\label{sec:comparison}

Use the zero-framed unreduced Jones polynomial \(\Phi_K(N;t)\) in
dimension color \(N\) for \(\mathfrak{sl}_2\). The quantum dimension
\([N]\) and the positive-color reduced polynomial are given by
\[
 [N]=\frac{t^{2N}-t^{-2N}}{t^2-t^{-2}},\qquad
 j_{K,N}(q)=\frac{\Phi_K(N;t)}{[N]}\quad(N>0),\qquad q=t^4.
\]
For the zero-framed unknot \(U\), \(\Phi_U(N;t)=[N]\), so
\(j_{U,N}(q)=1\) for \(N>0\). We use \(\Phi_K(1;t)=1\),
\(\Phi_K(0;t)=0\), and \(\Phi_K(-N;t)=-\Phi_K(N;t)\) for \(N>0\).
The operators \(M,L\) act by multiplication by \(t^{2N}\) and forward
shift of \(N\), so \(LM=t^2ML\). We use left coefficients and frequently
write \(u=L^{-1}\), \(\mathbb K=\Q(t,M)\) and
\(\D=\mathbb K[u^{\pm1};\sigma]\), where \(\sigma(M)=t^{-2}M\).
The forward shift acts on coefficients by \(\tau=\sigma^{-1}\).

Before localization the coefficient ring is
\(\Z[t^{\pm1},M^{\pm1}]\). Let \(I_K\) be the left ideal of operators
over this ring, with Laurent powers of \(L\), that annihilate the
signed sequence at every integer. Its rational saturation is
\(\D I_K\). A rational operator belongs to it exactly when left
multiplication by a nonzero Laurent-polynomial coefficient makes it
a member of \(I_K\).

For a generator of \(\D I_K\), shift the lowest \(L\)-exponent to zero
and clear denominators on the left. Divide the coefficients by their
content, the gcd of the nonzero coefficients in
\(\Z[t^{\pm1},M^{\pm1}]\). The resulting primitive minimum \(\alpha_K\)
is determined up to the coefficient units \(\pm t^iM^j\).
Its order is the difference between the largest and smallest shift
exponents, hence the \(L\)-degree after the shift to zero. The evaluation
\(M=t^{2N}\) is used to act on sequences. In
\(f_K=\alpha_K(-1,M,L)\) we instead specialize the coefficients at
\(t=-1\) with \(M\) independent.

The full classical polynomial defines the boundary eigenvalue curves in
the preferred meridian and longitude coordinates. Each irreducible factor
occurs once, including \(L-1\) for the abelian line.
We write \(A_K\) for its image in
\(R=\Q(M)[L^{\pm1}]\). This horizontal localization retains the
components that project dominantly to the meridian, since
meridian-only factors become units. We write \(\doteq\) for
association in \(R\).
The usual AJ equality, with the same abelian convention and up to
meridian factors, implies the radical equality
\begin{equation}\label{eq:AJ}
 \sqrt{(f_K)}=(A_K).
\end{equation}
For nonzero \(f_K\), this requires its square-free part to be associated
to \(A_K\) in \(R\).
An irreducible nonunit of \(R\) dividing \(f_K\) and coprime to \(A_K\)
contradicts this necessary condition. It suffices to bound the classical
curve by a containing polynomial coprime to that factor.

\begin{samepage}
\subsection{The four-family theorem}

The marking specifies the component orientations, the preferred
meridian--longitude pairs and the zero-framed companion embedding.
The negative clasp is the Whitehead pattern of Section~\ref{sec:mark}. The twist \(\ell\) is counted at the site in
Figure~\ref{fig:filling}, with sign fixed by \eqref{eq:twistgluing}.
Its action on the odd channel is \(q^{\ell n(n+1)}\). The positive clasp uses the reflected
clasp and the same signed disk-twist convention. Both use the faithful
zero-framed embedding of the companion and the preferred outer
meridian--longitude pair. Put \(T_+=T(r,s)\), \(T_-=T(r,-s)\), with
the orientation of \(T_-\) transported by reflection. Define
\begin{equation}\label{eq:E}
\begin{aligned}
 E(M,L)&=L^2-(M^4-4M^2+1)L+M^4,\\
 E^\vee(M,L)&=M^4 E(M^{-1},L)
              =M^4L^2-(M^4-4M^2+1)L+1 .
\end{aligned}
\end{equation}

\begin{theorem}[Both clasp signs and all integer twists]\label{thm:four}
Let \(2\le r<s\) be coprime integers, let \(\ell\in\Z\), and use the
marked satellite \(S=W_{\epsilon,\ell}(T_\eta)\) for each clasp sign
\(\epsilon\) and each companion chirality \(\eta\).
The actual signed rationally saturated annihilator is nonzero. With
\(M\) independent at specialization,
\[
 0\ne f_S=\alpha_S(-1,M,L),\qquad
 G\mid f_S,\qquad \gcd_R(G,A_S)=1,
\]
where \(G\) is specified by
\[
 \begin{array}{c|cc}
   &T_+&T_-\\ \hline
   \text{negative clasp}&E&E\\
   \text{positive clasp}&E^\vee&E^\vee
 \end{array}
\]
and \(G\) is irreducible in \(R\). In each case the radical
equality \eqref{eq:AJ} fails, so \(S\) is a counterexample to the AJ conjecture.
\end{theorem}
\end{samepage}

In the nondegenerate shape chart of Appendix~\ref{sec:geometry},
completing the axis cusp imposes both peripheral conditions on
that cusp. Their solutions project to \(E=0\) in the preferred outer
coordinates of~\eqref{eq:E} (Section~\ref{sec:shape}).
The same curve occurs in the transfer symbol of
Section~\ref{sec:symbol}, where \(Z+2\) is proportional to \(E\)
off \(x^2=1\) and \(L=-x\), with \(x=M^2\ne0\). Here \(Z\) is the trace of the axis longitude.
The geometric identification uses the full cusp conditions in
Section~\ref{sec:shape}. For the fixed monic integral representative
of~\eqref{eq:E}, Appendix~\ref{sec:geometry} also gives
\[
 \pi m(E)=\Vol(X_{\mathrm{Wh}})=4\cat ,
\]
where \(X_{\mathrm{Wh}}\) is the complete two-cusped Whitehead link
exterior, \(\Vol\) is its unsigned hyperbolic volume, and \(\cat\)
is Catalan's constant. The Mahler measure is the known square-weave
period in the normalization of~\eqref{eq:haar}
\cite[Sec.~4.1, Thm.~12]{CKL}.

\subsection{Selection of the quadratic}\label{sec:selection-intro}

Left multiplication of an annihilator can introduce extra factors
at specialization. Finding \(E\) in a constructed operator leaves
the central problem of proving that it survives in the actual minimal
recurrence.

For the untwisted negative-clasp double \(S=Z(T(2,3))\) of the
right-handed trefoil, let \(B_6\) be the longitude polynomial of the
axis filling of signed index \(6\), defined in~\eqref{eq:B}.
Theorem~\ref{thm:one} gives
\[
 \ord_L(\alpha_S)=16,\qquad f_S\doteq(L-1)^2B_6E,\qquad
 A_S\doteq(L-1)B_6 .
\]
Counting \(L-1\) once still leaves \(E\) as an extra irreducible factor.
The uniform theorem gives survival of the quadratic. The full
specialization and order above use the exact calculation in
Appendix~\ref{sec:whitehead}. For all untwisted cases,
Corollary~\ref{cor:untwisted-reverse} gives the divisibility
\(A_S\mid f_S\).

For the negative clasp put \(a=\eta rs-\ell\), where \(\eta=1\) for
the right-handed companion and \(\eta=-1\) for its mirror. The finite
satellite transform gives the inhomogeneous equation \(P_aF=b\) on the
reduced sequence \(F\) defined by the Jones invariant at positive colors.
The proof places \(b\) in a forcing module, a shift-stable module
of tail sequences over the rational coefficient field.
Sequences are identified on a tail if they agree for all sufficiently
large positive colors. Growth estimates show that the class of \(F\)
in the quotient is nonzero. Lemma~\ref{lem:relative} then selects
a nonunit right divisor of \(P_a\) that right-divides every tail
annihilator of \(F\).

The rational factor analysis determines what this right divisor must
retain. At isolated transverse intersections of specialized algebraic
factors, Theorem~\ref{thm:graph} requires an integral logarithmic
residue for the determinant of a quotient that separates them.
The computed half-integer residues obstruct those separations.
At nonzero effective index the rational operator is irreducible,
using separate rank-one tests at \(a=1\). At zero effective index the
proof excludes linear right divisors, which suffices to retain the
quadratic in the selected divisor. The factor survives the return to
the signed unreduced sequence.

On the classical side, a representation restricts to the companion
with one of two boundary relations. These give signed axis fillings
whose containing polynomial is coprime to \(E\)
(Section~\ref{sec:torusclass}). Thus \(E=0\) is absent as a
component of the satellite boundary curve, although individual
boundary pairs may lie on it. Whole-knot reflection gives the other
clasp.

The proof of the AJ discrepancy uses Zheng's exact finite-color
formula at generic \(q\)~\cite{Zheng}.
Appendix~\ref{sec:geometry} gives the elliptic normalization
\(Y^2=X^3-X\) and relates the regulator primitive to the stationary
action. For the negative clasp at zero twist and each fixed
right-handed companion \(T(2,d)\) with \(d\ge3\) odd, the additional
analytic estimates of Section~\ref{sec:physical} show that the
dominant window of the actual root-of-unity sum gives a nonzero
leading term and selects the lift \(w=\imunit\).

The index \(a=0\) occurs at \(\ell=\eta rs\), whereas the untwisted case
has \(\ell=0\) and \(a=\eta rs\). At \(a=0\) the selection argument
degenerates and requires finer coefficient information. In particular,
for the positive companion the extremal summands tie and the actual
leading coefficient \(N\) separates the sequence from the forcing module.
Table~\ref{tab:selection} gives the selection mechanisms for both
companion signs.

\begin{table}[htbp]
\centering\small
\begin{tabular}{@{}ll>{\raggedright\arraybackslash}p{.59\linewidth}l@{}}
\toprule
Companion & Index & Separation of the actual sequence from the forcing module
 & Selection \\
\midrule
\(T_+\) & \(a\ge1\) & Quadratic color valuation at \(t=0\).
At \(\ell=0\) the forcing is rational. & \S\ref{sec:select-pos-plus} \\
\(T_+\) & \(a=0\) & Leading coefficient \(N\) versus controlled coefficient
expansions of the finite Gaussian module. & \S\ref{sec:select-pos-nonpos} \\
\(T_+\) & \(a<0\) & Quadratic gap in the coefficient orders of the
uninverted expansion. & \S\ref{sec:select-pos-nonpos} \\
\(T_-\) & \(a<0\) & Quadratic color valuation at \(t=\infty\).
 & \S\ref{sec:select-neg-minus} \\
\(T_-\) & \(a\ge0\) & Quadratic gap in the coefficient orders after
inversion and the stated rescaling. & \S\ref{sec:select-neg-nonneg} \\
\bottomrule
\end{tabular}
\caption{Sequence selection for the negative clasp, with \(T_+\)
right-handed and \(T_-\) its mirror. The coefficient comparisons use the common forcing identity and coefficient estimates of
Sections~\ref{sec:select-gaussian} and~\ref{sec:select-coefficients}.}
\label{tab:selection}
\end{table}

Section~\ref{sec:operators} constrains the operator's rational
factors, and Section~\ref{sec:survival} uses the actual sequence to
select a right divisor. Sections~\ref{sec:irr-pos} and~\ref{sec:irr-neg}
use the two simple
intersections for \(a\ge2\) and \(a<0\), respectively. At \(a=1\)
the \(E/B_1\) intersection remains simple. The other split requires
the right and dual-left rank-one tests of Section~\ref{sec:irr-one}.
At \(a=0\), Section~\ref{sec:zero-factor} excludes first-order right
divisors. This weaker conclusion suffices. The selected right divisor
has degree two or three, and its specialization contains the quadratic
(Section~\ref{sec:selected-divisor}). Section~\ref{sec:torusclass}
proves classical exclusion. Sections~\ref{sec:signedreturn}
and~\ref{sec:wholemirror} give the signed return and whole-knot reflection.

\par\medskip
\noindent\begin{minipage}{\textwidth}
  \tableofcontents
\end{minipage}\par
\section{Identifying the framed pattern}\label{sec:mark}
The quantum sum uses Zheng's picture and the group calculation uses
Schilder's. We compare their preferred peripheral pairs so that the two
calculations apply to the same marked satellite. Write \(Z(K)\) for the untwisted satellite in Zheng's marking, and $\mu,\lambda$ for the positive meridian and oriented zero-linking longitude. The solid torus containing the pattern is the exterior of the axis; its meridian is the axis longitude. In gluing-adapted bases the zero-framed companion map sends $(\mu_C,\lambda_C)$ to $(\lambda_A,\mu_A)$.

Schilder records a classical \(A\)-polynomial formula for twisted
Whitehead doubles of torus knots in his own marking, with an attribution
to Ruppe's thesis~\cite[arXiv v1, Theorem~7.1]{Schilder}.
Proposition~\ref{prop:mark} compares the untwisted marked diagrams.
The matrices of the proposition carry the source pairs from Schilder's
link to Zheng's, with the target axis orientation adapted to gluing.
For the finite sum one can use the marked identification and the
signed gluing \eqref{eq:twistgluing} directly.

\begin{figure}[htbp]
\centering
\input{figures/marked-pattern}
\caption{Crossing correspondence on the outer pattern $P$ and unknotted
axis $A$. The left panel reproduces the left diagram of Schilder's
Fig.~1~\cite{Schilder}.
The right panel uses Zheng's diagram on p.~4~\cite{Zheng} at twist
$\ell=0$, with two parallel strands at the twist site. The added
crossing numbers refer to auxiliary walks used to compare the
diagrams. The peripheral pairs use the transported component
orientations. On each torus $\mu$ links its oriented component once positively, and
$\lambda$ is the oriented longitude of linking number zero. The
zero-framed gluing is $\mu_C\mapsto\lambda_A$, $\lambda_C\mapsto\mu_A$.}\label{fig:marks}
\end{figure}

\begin{proposition}\label{prop:mark}
Let $D_0$ be Schilder's untwisted pattern. For every oriented zero-framed companion $K$, $Z(K)$ is the mirror of $D_0(K^*)$, where $K^*$ is the mirror with its pushed-forward orientation. The map from Schilder's marked link to Zheng's reverses the ambient orientation. Its matrices in column bases $(\mu,\lambda)$ are $D=\diag(-1,1)$ on the outer torus and $-D$ on the gluing-adapted axis.
\end{proposition}
\begin{proof}
In Schilder's left drawing \cite[Fig.~1]{Schilder} number the upper and lower clasp crossings $1,2$. Number the inner bottom horizontal strand's two crossings with the axis from left to right $3,4$, and the outer strand's $5,6$. In Zheng's $WL(0)$ \cite[p.~4]{Zheng}, straighten the indicated twist region. Number the left and right clasp crossings $1,2$, the inner right strand's upper and lower belt crossings $3,4$, and the outer strand's $5,6$. Walk into crossing 1 from upper left in Schilder and down its left upper arm in Zheng. Walk the Schilder axis down its left side and the Zheng belt along its upper half from left to right. These are auxiliary walks for reading Figure~\ref{fig:marks}. Writing $o$ for an overpass and $u$ for an underpass, both walks give the cyclic words
\[
 P:\ 1o\,2u\,3o\,4u\,2o\,1u\,6u\,5o,\qquad
 A:\ 3u\,5u\,6o\,4o .
\]
In right/up page coordinates the signs of the determinant of the over and under tangents are respectively $(+,+,-,-,+,+)$ and $(-,-,+,+,-,-)$. Thus the correspondence preserves overpasses and reverses every vertex rotation. Map vertex disks, then edge strips, then the eight complementary disk faces of this connected planar graph. Extending with heights fixed over a thickening and radially over the two remaining balls gives an orientation-reversing homeomorphism of the labeled links. Figure~\ref{fig:filling} shows the belt and twist region.

With transported component orientations the positive meridian reverses and the oriented longitude is preserved, since linking numbers reverse and self-linking zero is preserved. Both longitudes are null homologous in the exterior and both meridian classes are nonzero, ruling out any added meridian shear. The self-crossing writhes of the pattern, $+2$ and $-2$, give precisely the corresponding blackboard-to-preferred corrections. Reverse the target axis orientation to adopt the gluing convention; this composes its raw matrix with $-I$. The compatibility needed to glue the two mirror maps is
\begin{equation}\label{eq:gluesquare}
 (-D)J_{\rm gl}=J_{\rm gl}D,\qquad
 J_{\rm gl}=\begin{pmatrix}0&1\\1&0\end{pmatrix}.
\end{equation} Collar isotopies realize the boundary square, so the pattern map glues to the mirror map of the companion. The outer meridian is still a meridian up to sign, giving the asserted extension over filling.
\end{proof}

Schilder writes the generators $x,y$ for outer and axis meridians and capitals for inverses. His relation and preferred words are
\[
 \Omega=yxYxyXyx=\operatorname{reverse}(\Omega),\qquad
 \lambda_x=YxyXyxYX,\quad\lambda_y=YXyxYxyX .
\]
In Tillmann's right-handed link let $a,b$ be the axis and outer
generators, with $A=a^{-1}$ and $B=b^{-1}$. Write
$R_T=(babABAba)(abABAbab)^{-1}$ for the relator. His preferred words are $\lambda_0=AbaBABab,\ \lambda_1=BabABAba$ \cite[Ch.~5, (5.10), (5.20)--(5.22)]{Tillmann}. The substitution $\phi(x)=b$, $\phi(y)=a^{-1}$ is invertible on free generators. By cancellation it sends the relator (left word times inverse right word) to $a^{-1}R_Ta$, and the two longitudes to $b\lambda_1 b^{-1}$ and $a\lambda_0^{-1}a^{-1}$. It therefore preserves the entire peripheral subgroups with matrices $I,-I$ respectively.

To realize this comparison topologically, Tillmann's complete hyperbolic triangulation gives an irreducible, boundary-incompressible core with two cusps. It is sufficiently large by Waldhausen's Lemma~1.1.6(a)~\cite{Waldhausen}. Schilder's link of two unknots is nonsplit: a split pair would have free group, whereas the isomorphic target group contains a cusp $\Z^2$. Its exterior is therefore irreducible, and a compressible torus boundary in an irreducible manifold would give a solid torus, excluded by the second cusp. Waldhausen's Corollary~6.5 realizes the peripheral isomorphism; both boundary matrices have determinant $+1$, giving the orientation. This identifies outer cusp 1 and axis cusp 0. On a corresponding common eigenline the mirror takes the right-handed coordinates $(u,v)$ to $(u^{-1},v)$. The other common eigenline inverts both eigenvalues, giving the
equivalent pair $(M,L)=(u,v^{-1})$.

At generic $q$ reflection inverts the bracket parameter and exchanges the two crossing smoothings. The Jones--Wenzl projectors and quantum trace transform by that inversion, writhe changes sign, and the quantum dimension is invariant. Thus $j_{Z(K),N}(q)=j_{D_0(K^*),N}(q^{-1})$ in zero framing. With the reduced invariant of Section~\ref{sec:comparison}, Zheng's
channel input before the twist weight is
$\Phi_K(2n+1;t)=[2n+1]j_{K,2n+1}(q)$. Division by $[N]$ reduces the
external color. These choices precede any root specialization.

\subsection{The signed full twist}

\begin{figure}[ht]
\centering
\begin{minipage}[c]{.36\textwidth}
 \centering
 \begin{tikzpicture}
  \node[anchor=south west,inner sep=0pt] (diagram) at (0,0)
    {\includegraphics[width=.93\linewidth]{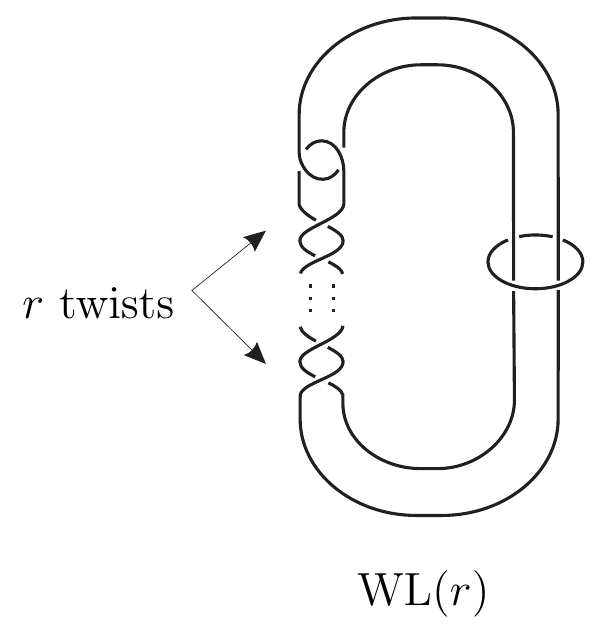}};
  \begin{scope}[x={(diagram.south east)},y={(diagram.north west)}]
    \fill[white] (.01,.477) rectangle (.29,.56);
    \node[inner sep=0pt,font=\scriptsize] at (.147,.516) {$\ell$ twists};
    \fill[white] (.565,.019) rectangle (.815,.103);
    \node[inner sep=0pt,font=\scriptsize] at (.69,.06) {$WL(\ell)$};
  \end{scope}
 \end{tikzpicture}
\end{minipage}\hfill
\begin{minipage}[c]{.57\textwidth}
 \small
 \textbf{Right belt.} This is the axis. Its preferred pair is
 $(\mu_{\rm ax},\lambda_{\rm ax})$.
 \smallskip

 \textbf{Left twist site.} The fixed clasp is above the twist site.
 The signed twist is $\ell$, with $\ell=0$ for the untwisted pattern.
 \smallskip

 \textbf{Axis filling.} The slope $\mu_{\rm ax}\lambda_{\rm ax}^{\,k}$
 inserts $-k$ full twists, counted with parallel test strands.
 The outer boundary keeps its preferred pair.
\end{minipage}
\caption{Zheng's Whitehead-link diagram (arXiv v4, p.~4)
\cite{Zheng}, with his twist label $r$ written as $\ell$.
The positive twist is the one in~\eqref{eq:twistgluing}. The axis belt surrounds the two strands of the
outer component. For the untwisted trefoil the nonabelian companion
restriction gives the filling slope
$\mu_{\rm ax}\lambda_{\rm ax}^{\,6}$.}\label{fig:filling}
\end{figure}

Write \(Z_\ell(K)=W_{-,\ell}(K)\), with \(Z(K)=Z_0(K)\).
The full twist at the indicated site acts on dimension \(2n+1\)
by \(q^{n(n+1)}\) in Zheng's variable convention~\cite[(2.3), (2.9)]{Zheng}.
The signed gluing is
\begin{equation}\label{eq:twistgluing}
 \mu_C\longmapsto\lambda_A,\qquad
 \lambda_C\longmapsto\mu_A\lambda_A^{-\ell}.
\end{equation}
The preferred outer pair is unchanged under the disk twist.
Write \(Z_\ell=Z_\ell(U)\) for the pure closure, where \(U\) is the
unknot and the patterned solid torus is embedded with zero framing.
The unknot longitude \(\lambda_C\) bounds a disk in its exterior,
so \eqref{eq:twistgluing} gives \(Z_{-k}\) by axis filling along
\(\mu_A\lambda_A^k\).

To check the sign, close two parallel test strands directed down the
twist site of Figure~\ref{fig:filling} to parallels of the core. The over-
and under-tangents in right/up page coordinates have directions
\((-a,-b)\) and \((a,-b)\) with \(a,b>0\), so both crossings of a positive
full twist contribute positively. This increases their linking by one.
A meridian disk in the \((1,k)\) filling differs by \(k\) oppositely
oriented longitudinal caps, correcting the linking of the test curves
by \(-k\). The disk twist on the axis side fixes \(\lambda_A\) and sends
\(\mu_A\) to \(\mu_A\lambda_A\). Moving the \(\ell\) diagrammatic twists
to the gluing composes the untwisted map with the inverse
\(\ell\)-th power of this action, which gives \eqref{eq:twistgluing}.
These disk relations hold with signed integers. The actual pattern has
linking zero with the axis, so the same linking correction leaves its
preferred longitude unchanged.

\section{The finite satellite transfer}\label{sec:family}

The index \(n\) labels the companion color of dimension \(2n+1\),
\(j\) indexes the intermediate sequence \(D_j\), and \(N\) is the
satellite dimension color. The inner finite transform carries the
recurrence for the twisted companion input to an equation for \(D\).
The outer transform gives \(P_aF=b\) on a tail of the reduced
satellite sequence \(F\).

For the negative clasp fix \(m=rs\), identify the companion sign
\(\eta\) with \(1\) or \(-1\), and set \(a=\eta m-\ell\).
Let \((z;q)_k=\prod_{i=0}^{k-1}(1-zq^i)\) and
\([j]=(q^{j/2}-q^{-j/2})/(q^{1/2}-q^{-1/2})\).
The centered torus sum \(V_n\) omits the torus framing prefactor.
The input to the finite transform includes both that prefactor, for
the zero-framed companion invariant, and the pattern twist.
\begin{equation}\label{eq:commoninput}
 V_n(q)=\sum_{k=-n}^n q^{rk(sk+1)}[2sk+1],\qquad
 v_n=q^{-a n(n+1)}V_n(q^\eta)
    =q^{\ell n(n+1)}\Phi_{T_\eta}(2n+1;t).
\end{equation}
For \(\eta=1\) this is Zheng's zero-framed torus formula
\cite[(4.1)]{Zheng} multiplied by \(q^{\ell n(n+1)}\).
For \(\eta=-1\) the companion input is obtained by inverting \(q\)
in the torus formula. Quantum dimensions are invariant under that
inversion. The negative-clasp kernel stays fixed in both cases.

\subsection{Two forms of the finite sum}

The outer output is \(F_N=j_{S,N}(q)=\Phi_S(N;t)/[N]\) for \(N>0\).
In Zheng's formula~\cite[(2.7)--(2.9)]{Zheng} this is
\(F_N=\sum_{n=0}^{N-1}\beta_{N,n}v_n\), where
\begin{equation}\label{eq:claspfinite}
 \beta_{N,n}=
 \frac{q^{N^2-1}(1-q)}{1-q^N}
 \sum_{i=0}^{N-n-1}q^{-N(i+n)}
 \frac{(q^{i+1};q)_n(q^{N-i-n};q)_n}{(q;q)_n}.
\end{equation}
The division by \([N]\) is part of \(\beta_{N,n}\). The sum has
\(F_1=1\). We return to this form for the coefficient estimates in
Section~\ref{sec:select-coefficients}. To transfer the companion
equation we use the equivalent triangular form, valid for arbitrary inputs \(v_0,\ldots,v_{N-1}\),
\begin{equation}\label{eq:universalbailey}
\begin{aligned}
 F_N &=\sum_{j=0}^{N-1} \mathcal C_j(N)D_j,\qquad
 D_j=\sum_{n=0}^j K_{jn}v_n,\\
 \mathcal C_j(N)&=\prod_{d=1}^j(1-q^{d+N})(1-q^{d-N}),\\
 K_{jn}&=\frac{q^j(1-q)(q;q)_j(-1)^nq^{n(n+1)/2}}
                  {(q;q)_{j-n}(q;q)_{j+n+1}} .
\end{aligned}
\end{equation}

In Zheng's channel trace~\cite[(2.4)]{Zheng}, the Hopf-belt
input of dimension \(p\) is \([p(2n+1)]\).

The two-color convention of Beliakova, Blanchet and L\^e
\cite[Secs.~1.1 and~1.4, p.~225]{BBL} divides the zero-framed invariant
by \([N][p]\). In the mirrored marking its reduced Whitehead value
is \(\sum_j b_j \mathcal C_j(N)\mathcal C_j(p)\), with
\(b_j=q^j(1-q)/(q^{j+1};q)_{j+1}\). To obtain this coefficient,
invert \(q\) in the printed coefficient
\((-1)^j q^{-j(j+1)/2}(1-q)/(q^{j+1};q)_{j+1}\)
and in both \(\mathcal C_j\) factors. Each inverted \(\mathcal C_j\) contributes
\(q^{-j(j+1)}\). Multiplying by \([p]\) restores the dimension of
the belt input and leaves the external reduction by \([N]\).

With \(\qbinom{d}{k}=(q;q)_d/[(q;q)_k(q;q)_{d-k}]\) for \(0\le k\le d\)
and zero otherwise, the finite binomial identity at \(Y=q^{p/2}\) is
\[
 (Y-Y^{-1})\mathcal C_j(p)=(-1)^{j+1}Y^{-2j-1}q^{j(j+1)/2}
                           (Y^2q^{-j};q)_{2j+1}.
\]
Its \(Y^{2n+1}\) coefficient is
\((-1)^n q^{n(n+1)/2}\qbinom{2j+1}{j-n}\).
The polynomial is odd under \(Y\mapsto Y^{-1}\), so these
coefficients express \([p]\mathcal C_j(p)\) in the inputs \([p(2n+1)]\).
Multiplying by \(b_j\) gives \(K_{jn}\).

The \(N\) inputs with \(p=1,\ldots,N\) span the \(N\)-dimensional
channel space: after division by \([2n+1]\) they are monic
Chebyshev polynomials of degrees \(0,\ldots,N-1\) evaluated at
the distinct points \(q^{n+1/2}+q^{-n-1/2}\). This proves
\eqref{eq:universalbailey} by linearity.

\subsection{Transport in the companion channel}\label{sec:inner}

Let \(Q=q^n,\ J=q^j,\ \widehat Q=(qQ)^{-1}\),
and put \(g(n)=q^{\ell n(n+1)}\).
Set \(\beta_i=-m+\lambda_i\) with
\((\lambda_1,\lambda_2,\lambda_3,\lambda_4)
=(s+r,-s+r,s-r,-s-r)\), and write
\[
\begin{aligned}
 H_+(q,Q)&=\sum_{i=1}^4 c_i(q)Q^{\beta_i},\qquad
 H_\eta(q,Q)=H_+(q^\eta,Q^\eta),\\
 (c_1,c_2,c_3,c_4)&=
 \frac{(q^{\beta_1+1},-q^{\beta_2},-q^{\beta_3},q^{\beta_4+1})}{q-1}.
\end{aligned}
\]
Removing the two terms \(k=\pm(n+1)\) in the centered torus sum
gives
\begin{equation}\label{eq:twistedinner}
 v_{n+1}=q^{-2a}Q^{-2a}v_n+g(n+1)H_\eta(q,Q).
\end{equation}
Indeed the unchanged terms have ratio \(q^{-2a}Q^{-2a}\).
For the positive companion the two new terms have respective
dimension numerators \(q^{s(n+1)+1}-q^{-s(n+1)}\) and
\(q^{1-s(n+1)}-q^{s(n+1)}\); their powers give the four terms
displayed in \(H_+\). The companion inversion gives \(H_-\).
The paired exponents satisfy \(\beta_i+\beta_{5-i}=-2m\), with matching
signs and extra powers. Consequently
\[
 H_\eta(q,q^{-2}Q^{-1})=q^{2\eta m}Q^{2\eta m}H_\eta(q,Q),
 \qquad H_\eta(q,q^{-1})=2 .
\]
These identities and \(g(-n-1)=g(n)\) extend
\eqref{eq:twistedinner} to \(v_{-n-1}=-v_n\), with
\(v_0-v_{-1}=2\) at the middle step.

The reflection \(n\mapsto -n-1\) sends \(Q\) to \(\widehat Q\).
Extend the kernel by \(K_{j,-n-1}=-K_{jn}\) on \([-j-1,j]\) and
by zero outside. The two odd extensions make \(K_{jn}v_n\) even.
Cancellation of consecutive factorials gives
\[
 K_{j,n+1}(1-q^2 JQ)=-q K_{jn}(Q-J).
\]
Thus for a weight \(f(Q)\) the sum of
\(K_{jn}(1-qJQ)f(Q)v_n\) on the reflected interval equals the
sum of \(-q K_{jn}(Q-J)f(qQ)v_{n+1}\).
The lower endpoint has \(1-qJQ=0\), the upper \(Q-J=0\).
Insert \(f(Q)=Q^{a-1}\) and \eqref{eq:twistedinner}. The same-input
multiplier on the left becomes
\(Q^{a-1}-qJQ^a+q^{-a}Q^{-a}-q^{-a}JQ^{-a-1}\).
Pairing the reflected terms gives the half-interval multiplier
as the sum of this expression at \(Q,\widehat Q\).

Use the trace polynomials normalized by \(S_0=2,\ S_1=Z,\ S_{k+1}=ZS_k-q^{-1}S_{k-1}\), extended
to integer \(k\). Then \(S_k(Q+\widehat Q)=Q^k+\widehat Q^k\)
and \(S_{-k}=q^k S_k\). Multiplication by \(Q+\widehat Q\) in
the kernel sum corresponds to
\[
 \mathcal Z=\frac{1}{qJ}+J-\frac{1-J}{J}\mathsf S,\qquad
 (\mathsf S D)_j=D_{j-1}.
\]
The shift \(\mathsf S\) acts on the intermediate index \(j\).
Indeed \(K_{j-1,n}/K_{jn}=(1-J/Q)(1-qJQ)/(q(1-J))\) for \(j>0\)
with a zero at \(n=j\). At \(j=0\) the \(\mathsf S\)-coefficient vanishes,
and \(\mathcal Z\) acts by \(1+q^{-1}\).
This proves the exact inner equation for every integer \(a\):
\begin{equation}\label{eq:innerforcing}
\begin{aligned}
 p_aD&=e,\qquad
 p_a=S_{a-1}(\mathcal Z)+(1-qJ)S_a(\mathcal Z)-qJ S_{a+1}(\mathcal Z),\\
 e_j&=\sum_{n=0}^j K_{jn}g(n)\mathcal R_\eta(Q,J),\\
 \mathcal R_\eta(Q,J)
 &=-q^{\eta m+\ell}\big(A(Q,J)-A(\widehat Q,J)\big),\\
 A(Q,J)&=(Q-J)Q^{\eta m+\ell-1}H_\eta(q,Q).
\end{aligned}
\end{equation}
The prefactor is \(-q^a g(n+1)/g(n)\) with its \(Q\)-power
put into \(A\). The reflected forcing gives a difference because
\(K_{jn}g(n)\) is odd. At \(j=0\), the left side is \(2(1-q)\),
using \(S_k(1+q^{-1})=1+q^{-k}\). The same value follows on the
right from \(Q=J=1\) and \(H_\eta(q,q^{-1})=2\).

\subsection{The operator in the satellite color}\label{sec:outer}

To carry the inner equation to the satellite color, we express
the weight from adjunction in a basis that gives shifts of \(N\).
Put \(x=M^2\), so \(x=q^N\) in evaluations. Take
\[
 d=\begin{cases}2a+3,&a\ge0,\\2b+2,&a=-b<0,\end{cases}
 \qquad \Delta(J)=\prod_{k=0}^{d-1}(1-xJ/q^k).
\]
Here \(d\) sets the pole window and the maximum shift in the
coefficient solve for \(P_a\). The adjoint of \(\mathcal Z\) in the sum against \(\mathcal C_j(N)\) is
\[
 Z^*h=\Big(\frac1{qJ}+J\Big)h(J)
     -\frac{(1-qJ)(1-qxJ)(1-qJ/x)}{qJ}\,h(qJ).
\]
This follows by reindexing \((\mathsf S D)_j\) and using
\(\mathcal C_{j+1}/\mathcal C_j=(1-qxJ)(1-qJ/x)\).
The endpoints vanish by \(1-J=0\) and \(\mathcal C_N(N)=0\).
Clearing the common denominator gives a polynomial map
\begin{equation}\label{eq:adjoint}
 \mathcal A_d h=\Delta Z^*(h/\Delta)
   =Jh+\frac{h(J)-(1-qJ)(1-qJ/x)(1-xJ/q^{d-1})h(qJ)}{qJ}.
\end{equation}
The numerator of the quotient vanishes at \(J=0\), and the map raises
degree by at most two. Adjunction places the left multipliers of
\(p_a\) on the seeds before applying the trace polynomials. The
cleared numerator is
\begin{equation}\label{eq:adjointnumerator}
 \mathcal N_a
   =S_{a-1}(\mathcal A_d)1+S_a(\mathcal A_d)(1-qJ)
                            -S_{a+1}(\mathcal A_d)(qJ).
\end{equation}
The bound on its degree is \(d\). For \(a\ge1\) this follows
term by term. At \(a=0\), combine the first and last terms
as \(qS_1(\mathcal A_d)(1-J)\). For \(a=-b\) rewrite the
numerator as
\[
 q^b\{qS_{b+1}(\mathcal A_d)1+S_b(\mathcal A_d)(1-qJ)
                                  -S_{b-1}(\mathcal A_d)J\}.
\]
The three degree bounds are \(2b+2,2b+1,2b-1\), respectively.

Define \(R_0=1\) and
\[
 R_i(J)=\prod_{k=0}^{i-1}
       \frac{1-x/q^k}{1-q^{k+1}/x}\,
       \frac{1-q^{k+1}J/x}{1-xJ/q^k}\qquad (1\le i\le d).
\]
The cleared functions \(\Delta R_i\) are a polynomial basis.
At \(J=q^{i-1}/x\) the columns \(\Delta R_k\) for \(k<i\) vanish,
and \(\Delta R_i\) has a nonzero pivot. Evaluate in descending order and then at
\(J=0\). The unique solution defines the actual coefficient row
\begin{equation}\label{eq:physicalrow}
 \mathcal N_a=\sum_{i=0}^d a_i(q,x)\Delta R_i,\qquad
 P_a=\sum_{i=0}^d a_i(q,M^2)u^i .
\end{equation}
In particular both chiralities give identical coefficients at a
fixed \(a\).

For \(N>d\), the factors in the evaluation of
\(\Delta(q^vJ)\) have exponents \(N+j+v-k>0\) for \(v\ge0\).
The basis constants have denominators \(1-q^{k+1-N}\), with
nonzero exponents. At a pivot the exponents of its factors include
\(k+i-2N\) for \(k<i\) and \(i-1-k\) for \(k\ge i\), also
nonzero. At \(J=0\) the remaining column is \(\Delta R_0=1\).
All factorial indices are in range. The cleared product relation
\(\mathcal C_j(N)R_i(q^j)=\mathcal C_j(N-i)\) implements the shift of the external
color, including zero terms for \(j\ge N-i\).
It follows that
\begin{equation}\label{eq:physicalforcing}
 (P_aF)(N)=\sum_{0\le n\le j<N}
      \mathcal C_j(N)K_{jn}g(n)\frac{\mathcal R_\eta(q^n,q^j)}{\Delta(q^j)}
      \qquad (N>d).
\end{equation}
Thus \(P_a\) is an inhomogeneous row on the actual reduced tail.

\subsection{The characteristic factor}\label{sec:symbol}

The relations \(J=q^j\) and \(x=q^N\) specify evaluation in the
finite sums. For the commuting symbol we first specialize at \(q=1\)
with \(x,J\) independent. The shift basis becomes
\[
 R_i\big|_{q=1}=\left(\frac{J-x}{1-xJ}\right)^i .
\]
Thus \(u=(J-x)/(1-xJ)\) is the commuting backward-shift coordinate.
This identifies the rational charts \((x,J)\) and \((x,u)\) at
generic \(x\). In these coordinates \(\mathcal A_d\) at \(q=1\)
acts by multiplication by \(Z\), where
\[
 J=\frac{x+u}{1+xu}=\frac{xL+1}{L+x},\qquad
 K=x+x^{-1},\quad Z=1+K-KJ+J^2,\quad e=1+J-K .
\]
Thus \(Z-2=(J-1)e\). Let \(T_0=1,\ T_1=Z-1+e\) and extend
\(T_{k+1}=ZT_k-T_{k-1}\) to all integers \(k\).
At \(q=1\) the numerator of the row is
\begin{equation}\label{eq:commonPi}
 \Pi_a=S_{a-1}(Z)+(1-J)S_a(Z)-JS_{a+1}(Z)
       =(1-J)(Z+2)T_a .
\end{equation}
Both sides obey the same invertible trace recurrence in \(a\).
At \(a=0,1\) equality follows from \(S_{-1}=Z\) and
\((J-1)e=Z-2\).
With this shift basis the symbol of \(P_a\) is
\((1-xJ)^{-d}\Pi_a\). In the preferred longitude chart,
\begin{equation}\label{eq:pullbacks}
 1-J=\frac{(1-x)(L-1)}{L+x},\qquad
 Z+2=\frac{(x+1)^2 E}{x(L+x)^2},\qquad
 1-xJ=\frac{(1-x^2)L}{L+x}.
\end{equation}
For positive \(k\) denote the cleared twist polynomials by
\begin{equation}\label{eq:B}
 B_k=x^k(L+x)^{2k}T_k,\qquad
 C_k=x^k(L+x)^{2k-1}T_{-k}.
\end{equation}
The classical longitude polynomials \(B_k,C_k\) are indexed
by \(k>0\). The finite transform weight is \(\mathcal C_j(N)\).
Up to units in \(R\), the symbol of the satellite operator has the factors
\[
 \begin{array}{cl}
 a\ge 1 & (L-1)EB_a,\\
 a=0 & (L-1)E,\\
 a=-b<0 & (L-1)EC_b .
 \end{array}
\]
The factor \(E\) is present in each symbol. To retain it in the
minimal recurrence, we first constrain the rational factors of
\(P_a\), then determine which right divisor is selected by the
actual sequence in \eqref{eq:physicalforcing}.

\section{Rational quotients at a crossing}\label{sec:residues}
The commuting symbol splits into algebraic factors. We ask which
blocks can lie in separate rational factors of the difference operator.
At a transverse crossing, a separating right factor would give a
quotient with an integral logarithmic residue of its determinant.
Theorem~\ref{thm:graph} makes this test under joint regularity and
isolation on the whole meridian fiber. A nonintegral residue then forces
the two blocks into the same part of every factorization.

In this section work over the coefficient field
$\mathcal F=\C((h))(z)$, with $\sigma z=e^{ch}z$ for $c\in\C^\times$,
and $\delta=cz\partial_z$. The constants are fixed.
The operators form $\mathcal D=\mathcal F[u^{\pm1};\sigma]$,
with left coefficients and $ua=\sigma(a)u$. The same criterion
applies to an embedded rational subfield with its induced shift and
residue field. Form the block partition over that residue field,
then embed a putative factorization for the local calculation.

For $P=SR$, the right factor gives the quotient of left modules
$\mathcal D/\mathcal DP\to\mathcal D/\mathcal DR$.
The shift on the ambient module acts in coefficient columns by
$v\mapsto A\sigma(v)$. The left factor gives the kernel of this quotient:
right multiplication by $R$ embeds $\mathcal D/\mathcal DS$ into
$\mathcal D/\mathcal DP$. Its dual is a quotient of the ambient dual,
whose shift matrix is $(A^{-1})^t$ in the dual basis.

At $z_0\ne0$ put $\mathcal O=\C[[h]][z]_{(h,z-z_0)}$.
Both $\sigma$ and $\sigma^{-1}$ preserve this ring.
The local hypothesis is
\begin{equation}\label{eq:joint}
 A,\ A^{-1}\in{\rm Mat}_d(\mathcal O).
\end{equation}
It controls the quantum matrix and inverse jointly in the meridian and deformation parameter, including meridians that move with $h$.

\begin{lemma}[Quotient determinant]\label{lem:integral}
For a rational quotient of such a module let $D=D_0+hD_1+\cdots$ be its shift determinant in a rational basis, with $v_h(D)=0$. Then $D_0$ is a unit at $z_0$ and
\begin{equation}\label{eq:detindex}
 {\rm ind}_{z_0}(D)=\Res_{z_0}\frac{D_1}{D_0}\frac{dz}{cz}\in\Z.
\end{equation}
\end{lemma}
\begin{proof}
\emph{1. Image lattices and finite support.}
To record the determinant divisor, use the images of the ambient
lattices. Write $Q$ for the full-row-rank matrix of the quotient map
on coefficient columns, and $B$ for the semilinear shift matrix
in the quotient basis. The intertwining identity is $B\sigma(Q)=QA$. Split the nonzero numerator and denominator polynomials of $D$, the entries of $Q$ and its maximal minors in a common finite Puiseux field $K'=\C((h^{1/e}))$. For each $\beta\in K'$ with limit $z_0$ denote the DVR $K'[z]_{(z-\beta)}$ by $O_\beta$ and put $\Lambda_\beta=QO_\beta^d$. This image is a free lattice of full quotient rank. Its determinant lattice is $(z-\beta)^{e_\beta}O_\beta$ in the fixed determinant line. Outside finitely many $\beta$ the entries are integral and a maximal minor is a unit, so $e_\beta=0$.

\emph{2. Transport and cancellation of the local divisor.}
Set $b=e^{ch}$. Substitution $z\mapsto bz$ identifies $O_{b\beta}$ with $O_\beta$. Every denominator that is a unit in $\mathcal O$ stays nonzero at $\beta$. Thus \eqref{eq:joint} sends the ambient lattice onto itself, and intertwining gives the equalities
\begin{equation}\label{eq:transport}
 B\sigma(\Lambda_{b\beta})=\Lambda_\beta,\qquad
 \ord_\beta D=e_\beta-e_{b\beta}.
\end{equation}
The difference of exponents in \eqref{eq:transport} is a
gauge divisor. To cancel it, take $g=\prod_{\beta\to z_0}(z-\beta)^{-e_\beta}$, a finite product. At every point of this neighborhood, including shifted neighbors of its support, $\sigma(g)/g$ has the divisor in \eqref{eq:transport}. Write $D=(\sigma(g)/g)V$. No zero or pole of $V$ tends to $z_0$.

\emph{3. Reduction and the logarithmic residue.}
We next check that the remaining factor $V$ contributes no logarithmic
residue at $z_0$.
A finite root of $V$ with another limit yields a linear factor with unit reduction at $z_0$. A root tending to infinity yields a scalar times $1-z/\beta$, also a unit after reduction. Factoring numerator and denominator this way shows that after removing its scalar Gauss order all Puiseux coefficients of $V$ are regular at $z_0$ and its leading coefficient is a unit there. In fact $v_h(V)=v_h(D)=0$ since the shift preserves Gauss order. The reduction of $g$ up to scalar is $g_0=(z-z_0)^j$, where $j=-\sum_\beta e_\beta$. The first shift term of $\sigma(g)/g$ is $h\delta g_0/g_0$: substituting $bz$ in a coefficient of strictly positive Puiseux order changes it only one order in $h$ later. Thus $D_0$ is a unit at $z_0$; in \eqref{eq:detindex} the sole polar contribution is $d\log g_0$, of residue $j$.
\end{proof}

\begin{theorem}[Crossing criterion]\label{thm:graph}
Let $P=P_0+hP_1+\cdots$ be monic, of positive order, with $\deg_u P_0=\deg_u P$, nonzero constant
term $P_0(u=0)$, and no repeated irreducible factor in $P_0$ over the
residual coefficient field. Choose two distinct monic irreducible blocks $G,H$ of $P_0$, with a common zero $(z_0,\lambda)$, $z_0\lambda\ne0$. Assume the companion of $P$ belongs to $\mathrm{GL}_d(\mathcal O)$, $G_uH_u\ne0$, and $J=G_uH_z-H_uG_z\ne0$ at this point. On the entire fiber over $z_0$ this must be the only meeting of distinct blocks, with no third block there. Set $T=P_0/(GH)$, $v=cz_0$. A rational right factor retaining $G$ and omitting $H$ requires
\begin{equation}\label{eq:crossindex}
 I_G=-\frac{P_1/T-\lambda v H_uG_z}{v\lambda J}\bigg|_{(z_0,\lambda)}\in\Z,
 \qquad I_H=-1-I_G .
\end{equation}
Take the residual blocks as vertices and join two when a crossing
satisfies these hypotheses and has nonintegral index. Each right
factor selects a union of graph components. In particular a connected graph forces irreducibility.
\end{theorem}
\begin{proof}
\emph{1. Endpoints and block partition.} The Gauss valuation (minimum coefficient order) is multiplicative: normalized reductions multiply in the ordinary polynomial domain because $\sigma$ acts trivially on the residue. In a monic factorization both monic factors have valuation at most zero and their sum is zero. Their constants have valuation zero because their product does. Each positive degree and both endpoints survive. Since $P_0$ has no repeated factors, each block occurs in exactly one of the two reductions. Applying the ordinary Gauss argument also at $z-z_0$ to these residual monic factors makes their coefficients regular and their constant terms units there. In particular there are no missing roots at infinity or zero on this fiber.

\emph{2. The first Ore coefficient.} Write $R_0=GA$, $S_0=HB$, so $T=AB$. Expanding the actual Ore product gives
\begin{equation}\label{eq:orelinear}
 P_1=S_1R_0+S_0R_1+(u\partial_uS_0)\delta R_0 .
\end{equation}
Over the local power series ring the simple root of $G$ is $u=g(z)$, $g'=-G_z/G_u$. At this root $R_1$ is the quotient of $P_1-(u\partial_uS_0)\delta R_0$ by $S_0$. The numerator at the node equals $P_1-\lambda vH_u G_zT$; the denominator has derivative $BJ/G_u$ along the branch. This identifies the possible pole at the chosen intersection. To recover the determinant residue, we must also control the other roots on the same meridian fiber.

\emph{3. The rest of the meridian fiber.}
In the completed ring $\C[[z-z_0]]$ the simple root splits off as $u-g$ from $R_0$. Every root of $R_0/(u-g)$ at $z_0$ lies outside $S_0=0$: a shared root would be another mixed collision, and $A(\lambda)B(\lambda)\ne0$. Thus the resultant of $R_0/(u-g)$ and $S_0$ is a unit. Modulo the complementary monic factor of $R_0$ the denominator is a unit since there is no further mixed root. This permits repeated roots within that complement.

The projector for the isolated root is $R_0(u)/[(u-g)R_{0,u}(g)]$. At $u=0$, after division by $R_0(0)$, its value is $-1/[gR_{0,u}(g)]$. Multiplying this projector value by the pole coefficient of $R_1(g)$ and using $R_{0,u}(g)=G_u A$ at the node gives $-(P_1/T-\lambda v H_uG_z)/(\lambda J)$ as the residue coefficient of $R_1(0)/R_0(0)$ in $z-z_0$. The weight $dz/(cz)$ divides it by $v$. This is \eqref{eq:crossindex} for the determinant $(-1)^{\deg R}R(0)$. Lemma~\ref{lem:integral} gives integrality. Interchanging $G,H$ gives the sum $-1$. Every proper partition of a connected graph cuts an edge. Jointly invertible scalar factors used in computing the local formula cancel at the node.
\end{proof}
\subsection{Coordinates, gauges and failure controls}
For $B'=C^{-1}B\sigma(C)$ the determinant changes by $\sigma(g)/g$, $g=\det C$. Normalizing its Gauss order as above adds $\ord_{z_0}g_0$ to the index. The rational-gauge invariant is the residue modulo $\Z$; jointly invertible gauges preserve the number itself. A singular presentation needs its own joint lattice check.

An unramified coordinate change pushes $cz\partial_z$ forward, preserving the one-form. An $h$-dependent regular change adds a regular term since $D_0$ is a unit. Rescaling $h$ rescales both numerator and shift vector. Extending constants by $h=s^e$ uses the first shift at order $s^e$. Pulling back $z-z_0=w^e$ multiplies the residue by $e$ and the lifted vector field can have a pole.

For the witness we take $h=q^{-1}-1$, $u=L^{-1}$, $\delta=x\partial_x$ with $x=M^2$; the rescaling from $\log q^{-1}$ has derivative one. These formulas hold on either unramified lift to $M$, where $\delta=(M/2)\partial_M$.

\emph{Dependence on the first deformation term.} Keep the specialized factors fixed while changing the
first quantum term. For $\sigma z=qz$, $q=e^h$, compare $u^2-2u+z(2-z)$ with $(u-z)(u-(2-z))$. Both have blocks $u-z,u+z-2$ crossing only at $(1,1)$, with $J=2$. The actual jets are respectively $0$ and $zu$, giving indices $(-1/2,-1/2)$ and $(-1,0)$. The first operator is irreducible, and its solution at $z=q^n$ with $f_0=f_1=1$ has precisely this minimum: the nonzero endpoint $q^n(2-q^n)$ permits backward recovery from any two zeros; a unit right gcd with a second annihilator would force a zero tail by a left B\'ezout identity.

\emph{Annihilator and minimum.} Consider the constant sequence. Its minimum is $u-1$, since a nonzero scalar cannot annihilate it. The larger row $(u^2-2u+z(2-z))(u-1)$ annihilates it but the two additional reduced factors disappear from its minimum.

\emph{Joint regularity.} A moving pole in rank one shows why joint regularity is required. The ratio $(z-1+(q-1)/4)/(z-1-(q-1)/4)$ has residual value one but first logarithmic coefficient $1/[2(z-1)]$, of index $1/2$ at $z=1$. Its moving pole prevents joint regularity.

We now return to sequences on the positive tail, identified when
they agree for all sufficiently large $N$. Embed $\mathbb K=\Q(t,M)$
by evaluation at $M=t^{2N}$. For a nonzero rational operator $P$, let
$v_{t+1}(P)$ be the minimum $(t+1)$-valuation of its coefficients.
Its initial polynomial is the specialization at $t=-1$ of
$(t+1)^{-v_{t+1}(P)}P$, with $M$ independent and $u=L^{-1}$ in $R$.
This normalization acts on the coefficient valuation.

The next lemma selects a divisor for a nonrational tail with
rational forcing, including zero forcing, and passes its initial
factors to the signed sequence. Lemma~\ref{lem:relative} generalizes
the selection step by using a submodule in place of the rational
scalars.

\begin{lemma}[The selected sequence]\label{lem:selected}
Suppose $F$ is a positive-tail sequence outside $\mathbb K=\Q(t,M)$, with $M=t^{2N}$, and $PF\in\mathbb K$ for an irreducible $P$ of positive shift span. Then $P$ right-divides every tail annihilator of $F$. On multiplying $F_N$ by $[N]$ and taking the odd extension, the
factors of the initial polynomial of $P$ survive in the specialization
of every primitive saturated all-integer annihilator. Such an annihilator exists.
\end{lemma}
\begin{proof}
For a nonzero polynomial $\sum a_j(t)M^j$ the distinct slopes in $\ord_t a_j+2Nj$ give a unique minimum eventually at either end. Rational evaluation therefore embeds the coefficient field in tails modulo eventual equality. If the greatest common right divisor of $P,H$ were a unit, left B\'ezout on $F$ with $HF=0$ would put $F$ in the field. Irreducibility now gives the division. If $PF=b\ne0$, $(L-\tau(b)/b)P$ annihilates $F$; for $b=0$ use $P$. Conjugation by the rational gauge $g=(M-M^{-1})/(t^2-t^{-2})$ has ratios $\tau^i g/g$ of valuation zero and residue one. The ratios therefore leave the initial factors unchanged under this
conjugation. Multiplicativity of the coefficient Gauss valuation
at $t+1$ preserves initial divisibility for rational left multiples.
Removing content from a polynomial representative leaves the initial
polynomial up to a unit in $R$.

For the odd extension, reflect a homogeneous gauged tail row by $(M,L)\mapsto(M^{-1},L^{-1})$. This is an Ore automorphism and its row kills the negative tail. A common left multiple exists by Euclidean division. After clearing coefficients multiply on the left by $\prod_{n\in S}(M-t^{2n})$ over its finite exceptional colors to obtain an all-integer annihilator. Conversely any saturated row restricts to the positive tail away from finitely many poles or zeros of its clearing multiplier. This proves survival without equating the two ideals. The stronger equality of minimum orders for the trefoil uses Appendix~\ref{app:domain}.
\end{proof}

\section{Rational factors of the satellite operators}\label{sec:operators}

We constrain the rational factors of the operator \(P_a\) from
\eqref{eq:physicalforcing}, which acts on the actual reduced satellite
sequence. We prove irreducibility over \(\Q(t,M)\) for \(a\ne0\)
and exclude first-order right divisors at \(a=0\).
Section~\ref{sec:survival} then uses the sequence to select a
right divisor in each case.

Garoufalidis and Koutschan detect possible right factors through
linear right factors of exterior powers
\cite[Thm.~4.8 and Sec.~7]{GaroufalidisKoutschan2013}.
The crossing criterion used here obstructs a partition of the
specialized factors by the first quantum correction to the operator.
When this calculation leaves a possible split at \(a=1\), the
right and dual-left rank-one tests address it directly.
The test at \(a=0\) only rules out linear right divisors.
For the knot \(7_4\), their multiplicity-retaining nonabelian
\(A\)-polynomial counts its nongeometric factor twice. Their
AJ comparison uses its reduced polynomial, which counts that
factor once, in the coordinates and mirror of their formula
\cite[Sec.~1.4 and Cor.~1.3]{GaroufalidisKoutschan2013}.
This is a multiplicity convention. Here \(E\) is an irreducible
factor absent from the full classical curve.

\subsection{Specialization and the full meridian field}\label{sec:twistirr}

We use \(q=e^\varepsilon,\ t=-e^{\varepsilon/4}\).
The shift is \(\sigma(M)=e^{-\varepsilon/2}M\) and the first
shift vector field on \(x\) is \(\delta=-x\partial_x\).
Since \(h=q^{-1}-1=-\varepsilon+O(\varepsilon^2)\), both the
first coefficient and the infinitesimal shift change sign
from the \(h\)-convention of Section~\ref{sec:residues}.
The normalized residue in \eqref{eq:crossindex} is unchanged.

\emph{Joint regularity of the coefficient solve.} At a fixed \(x_0\ne0,\pm1\), the map \eqref{eq:adjoint} acts on
polynomials over the joint local ring at \(q=1,x=x_0\).
The basis \(\Delta R_i\) reduces to
\[
 (J-x)^i(1-xJ)^{d-i}.
\]
Its matrix is the \(d\)-th symmetric power of an invertible
change of basis on linear forms since \(x_0^2\ne1\).
All denominators in the cleared basis coefficients are units there.
The actual solve \eqref{eq:physicalrow} is therefore jointly regular.
This argument works directly in the polynomial basis despite the
coalescence of auxiliary poles at \(q=1\).

\emph{Retention of the shift endpoints.} At \(J=x\) the specialized numerator is
\(x^{a-1}(1+x)(1-x^2)\) for \(a\ge0\), by \(S_k(K)=x^k+x^{-k}\).
At \(J=1/x\) use the same expression with \(x\) inverted.
These values isolate the two endpoint coefficients.
For \(a=-b<0\), put
\(Z_{\rm num}=(x+1)^2 E-2x(L+x)^2=x(L+x)^2Z\).
The recurrence gives
\[
 C_1=L+x^3,\qquad C_2=Z_{\rm num}C_1-x^2(L+x)^3,\qquad
 C_{k+1}=Z_{\rm num}C_k-x^2(L+x)^4 C_{k-1}\quad(k\ge2).
\]
It follows that \(C_b\) is monic of degree \(2b-1\), has
constant \(x^{4b-1}\), and
\(C_b(x,-x)=x^b(x^2-1)^{2b-1}\).
For the constant induction use \(Z_{\rm num}(x,0)=x^2(x^2+1)\).
The raw symbol in this case is exactly
\begin{equation}\label{eq:negative-residue}
 \overline P_{-b}=
 -\frac{(L-1)EC_b}{x^{b+1}(x-1)^{2b+1}(x+1)^{2b}L^{2b+2}} .
\end{equation}
Thus both shift endpoints are units at all such \(x_0\) in every
case. Monic normalization, and then pullback to \(x=M^2\),
give joint regularity of the companion and its inverse as in
\eqref{eq:joint}. At \(a=0\) the monic symbol in \(u\) is
\begin{equation}\label{eq:Eb}
 (u-1)E_{\rm back},\qquad
 E_{\rm back}=u^2+(-1+4/x-1/x^2)u+x^{-2}.
\end{equation}

\emph{Irreducibility of the residual blocks.} The discriminant of \(E\) is
\((x-1)^2(x^2-6x+1)\). The nonsquare part after \(x=M^2\)
has four simple roots, so \(E\), and hence \(E_{\rm back}\), is
irreducible over \(\C(M)\). Also \(E(M,1)=4x\).

For the other blocks let
\(\cheb_{-1}=0,\ \cheb_0=1,\ \cheb_{j+1}=z\cheb_j-\cheb_{j-1}\),
where \(z=Z\) and the meridian trace is \(K=x+x^{-1}\).
For \(a\ge2\) use the local polynomials
\(B=\cheb_{a-1},\,D=\cheb_a,\,A=B-D\).
The equation \(T_a=0\) is \(-A+eB=0\). The Chebyshev identity
\(D^2+B^2-zBD=1\) gives the formulas
\begin{equation}\label{eq:tracegraph}
 e=\frac{A}{B},\qquad J=1+(z-2)\frac{B}{A},\qquad
 K=z-\frac1{AB}.
\end{equation}
The identity gives coprimality of \(A,B\), so \(B=0\) yields no
point of the curve. If \(A=0\), then \(e=0\) and \(z=2\),
where \(A(2)=-1\). Thus these formulas describe the full graph
in the auxiliary variables. Moreover
\[
 (z+2)AB-1=-(\cheb_{2a}+\cheb_{2a-1})
\]
by the adjacent polynomial identity and the product formulas
\(\cheb_{2a}=D^2-B^2,\ \cheb_{2a-1}=2BD-zB^2\).
Writing \(z=2\cos\theta\) shows that the right side has
\(2a\) simple zeros, at \(\theta=2\pi j/(4a+1)\) for
\(1\le j\le2a\). None lies at \(AB=0\). They give simple zeros
of \(K+2\) with \(K-2\ne0\), so \(x^2-Kx+1\) is irreducible
over \(\C(z)\). At a simple root of \(B\), which exists for
\(a\ge2\), the function \(K\) has a simple pole and the small
branch \(x=K^{-1}+O(K^{-3})\) is unramified, of order one.
Thus \(x\) is nonsquare in the first quadratic extension and
adjoining \(M\) gives another quadratic extension.
Over \(\C[z,(AB)^{-1}]\) the graph is given by the monic
quartic \(M^4-K(z)M^2+1\), so its coordinate ring embeds
in the constructed field. The resulting equation for the
longitude is reduced and irreducible at generic \(M\).
Indeed \(T_a\) is monic of degree \(2a\) in \(J\),
\(B_a(x,-x)=x^a(1-x^2)^{2a}\), and at \(J=x\),
\(T_a=x^a\). These checks control the Möbius denominator
and give longitude degree \(2a\).

At \(a=1\),
\[
 T_1=J^2+(1-K)J+1,\qquad
 B_1=x^2L^2-(x^4-x^3-2x^2-x+1)L+x^2 .
\]
Its discriminant is
\((x-1)^2(x+1)^2(x^2-3x+1)(x^2+x+1)\).
The last two factors have disjoint simple roots, none zero.
Their product still has simple roots after \(x=M^2\).
Thus \(B_1\) is irreducible over \(\C(M)\).
Its monic constant is \(1\), whereas that of \(E\) is \(x^2\).

For \(a=-b<0\) take instead
\(B=\cheb_{b-1},\,D=\cheb_{b-2},\,A=B-D\).
Then \(T_{-b}=A-eB\), including \(b=1\), where \(A=B=1\).
The adjacent identity again gives \eqref{eq:tracegraph}.
The same exclusions hold, now with \(A(2)=1\).
The rational function \(K=(zAB-1)/(AB)\) has degree
\(2b-1\) since the numerator and denominator are coprime
of degrees \(2b-1,2b-2\). Hence \(z\), and also \(J\)
which recovers \(z\) together with \(K\), has degree \(2b-1\)
over \(\C(K)\). The extension \(\C(M)/\C(K)\) has degree four.
The compositum degree is a multiple of both degrees and is
at most their product, so the odd degree \(2b-1\) is retained
over \(\C(M)\). The equation \(T_{-b}\) has that degree in \(J\);
the longitude change is invertible at generic \(M\) and the
value \(C_b(x,-x)\) excludes any extra clearing factor.
This proves irreducibility of \(C_b\) over \(\C(M)\).

In each nonzero regime the three residual blocks are distinct.
For \(a\ge2\) they have degrees \(1,2,2a\).
At \(a=1\) use the two distinct quadratics just checked.
For \(b>1\) the degrees are \(1,2,2b-1\). For \(b=1\),
the two linear blocks \(L-1\) and \(L+x^3\) are distinct.
The monic Gauss argument in Theorem~\ref{thm:graph} then
forces any proper factorization to partition the residual
blocks, preserving shift degrees and both endpoints.

\subsection{The nodes and the physical first jet}\label{sec:jets}

The recurrence gives, for every integer \(a\),
\[
 T_a(2,e)=1+ae,\qquad T_a(-2,e)=(-1)^a(2a+1-ae).
\]
Call the varying block \(B_a\) when \(a>0\) and \(C_{-a}\)
when \(a<0\). Its intersection with \(L-1\) has
\[
 J=1,\quad K=2+1/a,\quad Z=2.
\]
Its intersection with \(E\) has
\begin{equation}\label{eq:enode}
 J=-\frac{2a-1}{2a+1},\quad K=-2-\frac{1}{a(2a+1)},\quad Z=-2 .
\end{equation}
At \(J=1\) the value \(Z=2\) separates \(L-1\) from \(E\).
At \(Z=-2\) the varying block fixes \(e=(2a+1)/a\), and then
\((J-1)e=-4\) fixes \(J\) uniquely. The two values of the meridian trace \(K\) are distinct and
neither is \(\pm2\). Each entire meridian fiber
has only the indicated meeting of distinct blocks. In the
\(J=\infty\) chart, \(\Pi_a\) has degree \(d\) and leading
coefficient \(-1\) for \(a>0\), \(1\) for \(a<0\).
Thus \(\Pi_a/(1-xJ)^d\) is nonzero there.
The endpoints \(u=0,\infty\) were checked above. All the
\(M\)-lifts on these fibers are unramified. This leaves the shift
derivatives, transversality, and first quantum coefficient to compute.

Write \(f_*=(1-xJ)^{-d}\), \(\mathcal H=\mathcal N_a/\Delta\),
\(A_*=1+J^2-KJ\) and \(D_*=x-x^{-1}\). A dot denotes
\(\partial_\varepsilon\) at zero. For \(\dot{\mathcal H}\) hold
\((x,J)\) fixed, and for \(\dot P_a\) hold \((x,u)\) fixed.
The fixed-\((x,J)\) derivative
of \(Z^*\) is \(s+\omega\partial_J\), with
\[
 s=2J^2-(1+K)J,\qquad \omega=-(1-J)A_*.
\]
In a power of an operator with residual multiplication by \(Z\),
inserting \(s+\omega\partial_J\) in each of the \(n\) positions
gives the derivative on a fixed seed \(v\):
\[
 \Big(nZ^{n-1}s+\tbinom n2 Z^{n-2}\omega Z_J\Big)v
       + nZ^{n-1}\omega v_J.
\]
At fixed \(Z\) the trace derivative follows from scaling the trace variable:
\[
 S_n(q,Z)=q^{-n/2}S_n(1,q^{1/2}Z),\qquad
 \dot S_n=(ZS_n'-nS_n)/2\qquad(n\ge0).
\]
These identities account for operator order as well as the explicit
\(q\)-dependence of each seed.

To recover the coefficient derivative of the satellite operator
from \(\dot{\mathcal H}\), subtract the contribution from the change
in the basis \(R_i\). Logarithmic differentiation in its \(k\)-th factor gives
the coefficient of \(k\) as
\(-1-J/(x-J)-xJ/(1-xJ)=(J^2-1)/A_*\) and the constant part as
\(\beta=1/(x-1)-J/(x-J)\). Consequently
\[
 \dot R_i/u^i=\rho i(i-1)+\beta i,\qquad \rho=(J^2-1)/(2A_*).
\]
Let \(p=\Pi_a=GH\) locally, where \(G=1-J\) or \(G=Z+2\)
at the corresponding node, and \(H\) includes the third block.
Use \(D_u=u\partial_u\) and \(D_x=x\partial_x|_u\). At a
mixed double zero \(p=D_up=0\), so differentiation of
\(\mathcal H=\sum_i a_iR_i\) gives
\begin{equation}\label{eq:physical-basis-jet}
 \frac{\dot P_a}{f_*}=\frac{\dot{\mathcal H}}{f_*}
                            -2\rho(D_uG)(D_uH).
\end{equation}
The coordinate derivatives are
\begin{equation}\label{eq:physical-derivatives}
 D_u=\frac{A_*}{D_*}\partial_J,\qquad
 D_x=D_*\partial_K+\frac{J^2-1}{D_*}\partial_J .
\end{equation}

At a simple mixed node use the slopes and normalized jet
\[
 s_G=-\frac{D_xG}{D_uG},\qquad
 s_H=-\frac{D_xH}{D_uH},\qquad
 \kappa=\frac{\dot P_a}{f_*(D_uG)(D_uH)}.
\]
The two indices in \eqref{eq:crossindex} are
\begin{equation}\label{eq:physical-indices}
   -\frac{\kappa-s_G}{s_H-s_G},\qquad
   -\frac{\kappa-s_H}{s_G-s_H}.
\end{equation}
Indeed, after monic normalization let the crossing factors be
\(G_0,H_0\) with cofactor \(T_0\). The local multipliers in \(G,H\)
give \((D_uG)(D_uH)=
(\bar a_d/f_*)T_0 u^2(G_0)_u(H_0)_u\) at the node.
The actual first monic coefficient is \(\dot P_a/\bar a_d\) there,
since the derivative of the normalization multiplies a zero.
This gives \eqref{eq:physical-indices} with \(\delta=-D_x\),
without making any operator gauge out of the local functions.
On a nonzero \(M\)-lift the vector field is
\(-\tfrac12M\partial_M\), so the determinant one-form and index
are unchanged.

\subsection{Positive effective index}\label{sec:irr-pos}

\emph{The \(E/B_a\) crossing.} For \(a\ge1\), put \(s_a=(-1)^a\) at this node.
Use
\[
 S_n(\pm2)=2(\pm1)^n,\quad
 S_n'(\pm2)=(\pm1)^{n-1}n^2,\quad
 S_n''(\pm2)=(\pm1)^n n^2(n^2-1)/6 ,
\]
which follow by differentiating the trace recurrence at \(2\)
and using parity. In differentiating \(p\) in \(Z,J\) take
these as independent variables before the symbol substitution.
At this node
\[
 p_Z=0,\quad p_{ZJ}=-s_a(2a+1),\quad
 p_{ZZ}=-2s_a a(2a-1)/3.
\]
The scalar adjoint term and common-weight derivative vanish.
The fixed-\(Z\) derivative of the traces contributes
\(-s_a(1+J)\), while the seed derivative \(-J(S_a+S_{a+1})\)
vanishes. Thus
\[
 \dot{\mathcal H}/f_*=
 \omega(p_{ZJ}+Z_Jp_{ZZ}/2)-s_a(1+J)
 =s_a\frac{56a^2+16a+2}{3(2a+1)^3}.
\]
With \(G=Z+2\), the free-coordinate derivatives of \(H=p/(Z+2)\)
are \(H_Z=p_{ZZ}/2,\ H_J=p_{ZJ}\).
Substitution in \eqref{eq:physical-derivatives} gives
\begin{align*}
 s_G&=-\frac{4a}{4a+1},&
 s_H&=\frac{24a^3+8a^2+2a}{10a^2+5a+1},\\
 (D_uG)(D_uH)&=
 -s_a\,\frac{2(4a+1)(10a^2+5a+1)}
                 {3a(2a+1)^2(8a^2+4a+1)},&
 \rho&=-\frac{4a^2}{2a+1}.
\end{align*}
The product and slope difference are nonzero. Equation
\eqref{eq:physical-basis-jet} makes \(\kappa=(s_G+s_H)/2\),
so both indices are \(-1/2\). In particular at \(a=1\)
the derivative product is \(160/351\), the raw jet divided
by \(f_*\) is \(106/351\), and the same simple crossing binds
\(E\) to \(B_1\).

\emph{The \((L-1)/B_a\) crossing.} At this node \(\omega=\rho=0\)
and \(p_Z=-4a\). The scalar adjoint, parameter and seed contributions
to \(\dot{\mathcal H}/f_*\) are respectively \(4a+4,-4a+2,-4\),
summing to \(2\). Here \(G=1-J\), \(H=(Z+2)T_a\), and
\((T_a)_Z(2,-1/a)=(a+1)(2a+1)/6,\ (T_a)_e=a\), from
\(T_a=\cheb_a+(e-1)\cheb_{a-1}\). For \(a\ge2\)
\eqref{eq:physical-derivatives} yields
\[
 (D_uG)(D_uH)=-\frac{2(a-1)}{3a},\qquad
 s_G=0,\qquad s_H=-\frac{6a}{a-1},\qquad
 \kappa=-\frac{3a}{a-1}.
\]
Thus the second edge also has two simple transverse branches and
indices \(-1/2,-1/2\). The two edges connect all three blocks.

\subsection{Negative effective index}\label{sec:irr-neg}

For \(a=-b\), \(b\ge1\), the numerator is the three-seed
expression after \eqref{eq:adjointnumerator}.
Use \(p=S_{b+1}+(1-J)S_b-JS_{b-1}\) at \(q=1\).
For clarity, differentiating the cleared map \(\mathcal A_d\) gives
\(\alpha_d+\omega\partial_J\), where
\[
\begin{aligned}
 T^{(1)}={}&-J(1-J/x)(1-xJ)
 -(J/x)(1-J)(1-xJ)\\
 &+(d-1)xJ(1-J)(1-J/x),\\
 \alpha_d={}&-(Z-J)-T^{(1)}/J .
\end{aligned}
\]
This is the derivative of the scalar product and dilation in
\eqref{eq:adjoint}. Inserting it in the power formula and
differentiating \(q^b\), the traces and the seeds gives
\begin{equation}\label{eq:negative-jet}
\begin{aligned}
 \dot{\mathcal N}_{-b}={}&\tfrac b2 p
  +\tfrac12 S_{b+1}-JS_b-\tfrac12 JS_{b-1}
  +(Z/2+\alpha_d)p_Z\\
  &+\omega (Z_Jp_{ZZ}/2+p_{ZJ}) .
\end{aligned}
\end{equation}
At either node \(p=0\), so \(\dot{\mathcal H}/f_*=\dot{\mathcal N}_{-b}\).

\emph{The \((L-1)/C_b\) crossing.}
At \(J=1,\ K=2-1/b\), one has \(p_Z=4b\),
\(Z/2+\alpha_d=1/b\), and \(\omega=\rho=0\).
Equation \eqref{eq:negative-jet} gives \(-2+4=2\).
Using \(G=1-J,\ H=(Z+2)T_{-b}\), the derivatives at the
node are \(H_J=-2(4b-1)(b+1)/(3b)\), \(H_K=4b\).
These follow by substituting \(p_Z=4b\),
\(p_{ZJ}=-(2b^2-2b+1)\), \(p_{ZZ}=(4b^3+2b)/3\) in
the derivative of \(p=(1-J)H\) after the symbol substitution.
Equation \eqref{eq:physical-derivatives} now gives
\[
 (D_uG)(D_uH)=-\frac{2(b+1)}{3b},\qquad s_G=0,\qquad
 s_H=-\frac{6b}{b+1},\qquad \kappa=-\frac{3b}{b+1}.
\]

\emph{The \(E/C_b\) crossing.}
At \(Z=-2\), put \(h_b=2b-1\) and \(s_b=(-1)^b\).
Here \(A_*=1/(bh_b)\), \(D_*^2=(8b^2-4b+1)/(b^2h_b^2)\),
\(Z_J=-(4b-1)/(bh_b)\) and \(\omega=-4/h_b^2\).
The same trace derivatives give \(p_Z=0,\ p_{ZJ}=s_b h_b\),
\(p_{ZZ}=-2s_b b(2b+1)/3\). Hence
\[
 \dot{\mathcal N}_{-b}=-s_b\frac{56b^2-16b+2}{3h_b^3}.
\]
Now \(G=Z+2,\ H=p/(Z+2)\) with \(H_Z=p_{ZZ}/2\),
\(H_J=p_{ZJ}\) in the free coordinates. Equations
\eqref{eq:physical-basis-jet}--\eqref{eq:physical-derivatives} give
\begin{align*}
 s_G&=-\frac{4b}{4b-1},&
 s_H&=-\frac{24b^3-8b^2+2b}{10b^2-5b+1},\\
 (D_uG)(D_uH)&=
 -s_b\,\frac{2(4b-1)(10b^2-5b+1)}
               {3b h_b^2(8b^2-4b+1)},&
 \rho&=\frac{4b^2}{h_b},\\
 s_H-s_G&=-\frac{6b h_b(8b^2-4b+1)}
                    {(10b^2-5b+1)(4b-1)}.
\end{align*}
All relevant derivatives are nonzero for \(b\ge1\).
The normalized jet is \(\kappa=(s_G+s_H)/2\).
For \(a\ge2\) or \(a<0\), the two simple crossings give the path
\[
 L-1\ \xleftrightarrow{\;(-1/2,-1/2)\;}\ 
 \begin{cases} B_a,&a\ge2,\\ C_{-a},&a<0 \end{cases}
 \ \xleftrightarrow{\;(-1/2,-1/2)\;}\ E .
\]
Each edge records the two indices for retaining either incident
block. For the two ranges in the display the resulting connected
graph obstructs factorization over \(\C((\varepsilon))(M)\), hence
over the embedded \(\Q(t,M)\).
At \(a=1\), separating \(L-1\) requires the rank-one argument below.
At \(a=0\) the varying block is absent and we only test first-order
right divisors.

\subsection{The rank-one split at \(a=1\)}\label{sec:irr-one}

At \(K=3\) the block \(T_1=(J-1)^2\) has a double shift root.
The simple crossing with \(E\) has already forced \(E\) and \(B_1\)
into the same factor of any partition. The only remaining split is
the block \(u-1\) against their product. This allows a rank-one
calculation on both sides of the factorization.

Let \(P\) be \(P_1\) made monic and \(H_0=\bar P/(u-1)\).
If a right factor is \(u-1+\varepsilon r_1+\cdots\), evaluating
\eqref{eq:orelinear} at \(u=1\) yields
\(r_1=\dot P(x,1)/H_0(x,1)\).
If a left factor is \(u-1+\varepsilon s_1+\cdots\), the
same product formula yields
\[
 s_1=\frac{\dot P(x,1)-\delta H_0(x,1)}{H_0(x,1)} .
\]
Use the raw row to compute these ratios. At \(u=1\),
\[
 H_{\rm raw}(x,1)=
 4(1-x)^{-5}\frac{x-1}{x+1}(3-K).
\]
It has a simple \(x\)-zero at each solution \(x_0\) of \(K=3\):
\((x-1)(x-x^{-1})/(x+1)=K-2\). Thus
\[
 \partial_x H_{\rm raw}(x_0,1)
   =-\frac{4(1-x_0)^{-5}}{x_0},\qquad
 \delta H_{\rm raw}(x_0,1)=4(1-x_0)^{-5}.
\]
The raw jet is \(2(1-x_0)^{-5}\) by the abelian crossing
calculation, where \(\rho=0\) also at \(a=1\).
Monic normalization scales both numerator and zero derivative
by the same unit at the point. The rank-one determinant starts
with \(1-\varepsilon r_1\) on the right and
\(1-\varepsilon s_1\) on the left. With shift vector field
\(-x\partial_x\), the normalized indices are
\[
 \Res_{x_0} r_1\,dx/x=-\tfrac12,\qquad
 \Res_{x_0} s_1\,dx/x=+\tfrac12 .
\]
The right index must be integral by Lemma~\ref{lem:integral}.
For the left factor use the dual quotient described in
Section~\ref{sec:residues}. Its ambient shift matrix \((A^{-1})^t\)
and inverse are jointly regular. The determinant of the rank-one
dual is reciprocal to the submodule determinant, so integrality on
the dual quotient implies integrality of the left index. Both residues contradict
integrality, excluding the remaining partition.
This proves irreducibility of the physical \(P_1\).

\subsection{First-order right divisors at \(a=0\)}\label{sec:zero-factor}

The three-pole solve \eqref{eq:adjointnumerator}--\eqref{eq:physicalrow}
at \(a=0\), after multiplication by
\(x(x-q)(x-q^2)(x^2-q^3)(x-1)\), gives the four coefficients
\begin{equation}\label{eq:zero-coeff}
\begin{aligned}
 b_0={}&q^5(x-q^2)(x-1),\\
 b_1={}&-q^3(x-q)^2(2q^3-q^2x-3qx+x^2),\\
 b_2={}&q^2(x-q^2)^2(q^3-q^2x-3qx+2x^2),\\
 b_3={}&-qx^2(x-q^3)(x-q).
\end{aligned}
\end{equation}
They are obtained by matching the cleared polynomial
\(qS_1(\mathcal A_3)(1-J)+2(1-qJ)\) in the basis
\(\Delta R_0,\ldots,\Delta R_3\).
We take \(P=b_3^{-1}\sum_i b_i u^i\). Its residue is
\eqref{eq:Eb}. A monic linear right divisor over \(\C(t,M)\)
must reduce to \(u-1\). Two distinct such divisors would have
a monic degree-two least common left multiple right-dividing
\(P\), by the right Euclidean algorithm. Its degree-two
residual divisor of \(\bar P\) would be \(E_{\rm back}\).
Degree-preserving reduction of its factorization by one of
the linear divisors would make \(u-1\) divide \(E_{\rm back}\),
contrary to \(E_{\rm back}(1)=4/x\). Thus the monic linear
right divisor is unique if it exists. The same argument works
over the even-meridian field, because \(E_{\rm back}\) is
already irreducible over \(\C(M)\).

Suppose \(u-r_0\) is a right divisor over \(\Q(t,M)\).
The involution \(M\mapsto-M\) commutes with the shift and fixes
\(P\); uniqueness implies \(r_0\in\Q(t,x)\).
Only after this descent does the order-four automorphism
\(t\mapsto\imunit t\) commute with the shift. On \(M\) the
multiplier \(t^{-2}\) changes sign under this map, whereas
the multiplier \(t^{-4}\) on \(x\) is fixed. Apply uniqueness over \(\C(t,x)\) and
intersect the fixed field with \(\Q(t,x)\) to get
\(r_0\in\Q(q,x)\). Here \(\sigma(x)=x/q\). The equation is
\[
 b_0+b_1r_0+b_2 r_0\sigma(r_0)+b_3 r_0\sigma(r_0)\sigma^2(r_0)=0.
\]
In particular \(r_0\ne0\).
In the basis \(1,u,u^2\) of the ambient cyclic module let \(C\)
be the column-action matrix, so the shift sends \(c\) to
\(C\sigma(c)\). The quotient to the module for \(u-r_0\)
is the row \(v^t\), where \(v=(1,r_0,r_0\sigma(r_0))^t\).
Its intertwining identity is \(v^t C=r_0\sigma(v^t)\).
Thus \(A=C^t\) has last row \((-b_0,-b_1,-b_2)/b_3\) and
first rows \((0,1,0),(0,0,1)\), and
\(Av=r_0\sigma(v)\).
The companion and inverse are integral and invertible
at every finite nonzero prime of \(\Q(q)[x]\) outside the orbit
of \(x-1\). Indeed the possible finite nonzero poles in
\eqref{eq:zero-coeff}, after monic division and inversion,
are at \(1,q,q^2,q^3\).

For such a good prime \(\mathfrak p\), let \(e_{\mathfrak p}\)
be the minimum valuation of the coordinates of \(v\). The
matrix equation gives
\[
 \ord_{\mathfrak p} r_0=e_{\mathfrak p}
                         -\ord_{\mathfrak p} \sigma(v).
\]
Only finitely many \(e_{\mathfrak p}\) are nonzero. On good
orbits their product of prime powers \(h\) has
\(\ord_{\mathfrak p}\sigma(h)=\ord_{\mathfrak p}\sigma(v)\).
Consequently \(r_0/(h/\sigma(h))\) is supported on the
remaining primes. With \(g=h^{-1}\), factorization in
\(\Q(q)[x]\) gives
\begin{equation}\label{eq:riccati-normal}
 r_0=c(q)x^\mu\prod_{k\in\Z}(x-q^k)^{n_k}\,\frac{\sigma(g)}{g},
\end{equation}
where \(\mu,n_k\) are integers and the product is finite.
This accounts for every finite divisor, including \(x=0\).
Infinity is controlled by the degrees of the resulting rational
functions.

Normalize \(g\) by a scalar parameter power so its reduction is a
nonzero rational function \(g_0(x)\). Since \(\overline{r_0}=1\)
and \(\overline{\sigma(g)/g}=1\),
\eqref{eq:riccati-normal} gives \(\mu=0\), \(\sum n_k=0\)
and \(c(1)=1\).
Put \(d_0=\ord_{x=0}g\), \(s=\sum k n_k\).
At \(x=0\), the right side is \(c(q)q^{s-d_0}\).
The coefficients of the Riccati equation there are
\((q^7,-2q^8,q^9,0)\), giving \(r_0(q,0)=q^{-1}\).
Hence \(c(q)=q^\kappa\) with \(\kappa=-1-s+d_0\in\Z\).

Now set \(q=e^\varepsilon\) and \(r_0=1+\varepsilon f(x)+\cdots\).
The monic residue gives coefficient \(4/x\) for \(f\) in the
first-order Riccati equation. The derivative of
\(\sum b_i/b_3\) at \(q=1\) is
\(-2(x+1)/(x(x-1))\) by \eqref{eq:zero-coeff}.
Thus \(f(x)=(x+1)/(2(x-1))\).
On the other hand \eqref{eq:riccati-normal} gives
\[
 f(x)=\kappa-\frac{s}{x-1}-x\frac{g_0'(x)}{g_0(x)} .
\]
At infinity this has integer limit given by \(\kappa\) minus the
degree difference of numerator and denominator of \(g_0\).
The actual first-order equation gives limit \(1/2\).
This contradiction excludes linear right factors over
\(\Q(t,M)\). The calculation specializes the rational function
before taking its degree at infinity, so moving poles of a
hypothetical factor are allowed. With the nonzero sequence class,
this exclusion is enough for the selection in
Section~\ref{sec:selected-divisor}.

\section{Selection by the actual sequence}\label{sec:survival}

The factor analysis limits the right divisors of \(P_a\).
To select one for the Jones sequence, we put the forcing in a
submodule and show that the actual sequence has a nonzero class
in the quotient.

Let \(\mathcal S\) be the space of \(\Q(t)\)-valued tail
sequences with the eventual equality of Section~\ref{sec:selection-intro}. The coefficient field
\(\mathbb K\) acts faithfully by \(M=t^{2N}\). Indeed, distinct
monomials \(t^iM^j\) have exponents \(i+2Nj\) after evaluation.
A finite list of distinct monomials has distinct exponents away
from finitely many \(N\). Thus every nonzero coefficient
is defined and invertible on a tail. With the backward shift, this
makes \(\mathcal S\) a left \(\D\)-module.

\begin{lemma}[Relative forcing]\label{lem:relative}
Let \(\mathcal U\subset\mathcal S\) be a left \(\D\)-submodule
with \(PF\in\mathcal U\), \(P\ne0\) and \(F\notin\mathcal U\).
The annihilator of \(F+\mathcal U\) has a nonunit generator \(R\)
that right-divides \(P\) and every actual tail annihilator of
\(F\). If \(P\) is irreducible, \(R\) is associated to \(P\).
\end{lemma}
\begin{proof}
The annihilator of the nonzero class is proper and contains \(P\).
It is a nonzero principal left ideal by skew Euclidean division.
Every row killing \(F\) kills its class. These inclusions give
the right divisions. An irreducible \(P\) and a nonunit \(R\)
in \(P=SR\) force \(S\) to be a unit.
\end{proof}

We first separate by color valuations in two ranges, then use
coefficients of controlled expansions in the other ranges.
Write \(\nu=\tfrac14\ord_{t=0}\), so \(\nu(q)=1\);
on expressions in \(q\) this is the lowest \(q\)-exponent.
These are growth valuations, independent of the \(t+1\) valuation
used for specialization. The torus input \(V_n(q)\) has
\(\nu(V_n)=0\) and leading coefficient \(1\) for every \(n\).
Indeed \(k=0\) contributes \(1\). For \(k=p>0\) the
least exponent is \(p(mp+r-s)>0\);
for \(k=-p<0\) it is \((rp-1)(sp-1)>0\).

\subsection{Color growth for the positive companion at \(a\ge1\)}\label{sec:select-pos-plus}

For \(\eta=1\), \(a\ge1\), the lowest \(n\)-term in \(D_j\)
is uniquely \(n=j\), since its \(n\)-dependent exponent is
\((\tfrac12-a)n(n+1)\). The \(\mathcal C_j(N)\) factor has valuation
\(j(j+1)/2-jN\) and leading coefficient \((-1)^j\).
The \(j\)-th outer summand leading coefficient is consequently
\(1\), with valuation \(j-(a-1)j(j+1)-jN\).
Its successive difference is \(1-N-2(a-1)j\). For \(N>1\) the
minimum is uniquely \(j=N-1\). Together with \(F_1=1\) this gives
\begin{equation}\label{eq:actual-eta-plus}
  \nu(F_N)=-(N-1)(aN-1),\qquad \operatorname{lc} F_N=1 .
\end{equation}

At \(\ell=0\) the forcing belongs to \(\mathbb K\).
Indeed the binomial calculation in Section~\ref{sec:family}
gives, on the full reflected interval,
\[
 \sum_n K_{jn} q^{bn}=b_j(1-q^{-b}) \mathcal C_j(b)
\]
for every integer \(b\). The formula follows by multiplying the
identity for \([b] \mathcal C_j(b)\) by \(q^{-b/2}(q^{1/2}-q^{-1/2})b_j\).
The right side vanishes for \(b=0\) or \(j\ge |b|\).
At \(\ell=0\), the exponents in \(A(Q,J)\) of
\eqref{eq:innerforcing} are \(\lambda_i,\lambda_i-1\).
Hence \(e_j=0\) for \(j\ge r+s+1\). The outer forcing is the
finite rational sum \(\sum_{j<r+s+1} \mathcal C_j(N)e_j/\Delta(q^j)\).
A nonzero rational function evaluated on the tail has eventually
affine valuation, so \eqref{eq:actual-eta-plus} excludes \(F\)
from that field.

For every integer \(\ell\), define a weighted forcing module
\(\mathcal W_\ell\). We use it for the other twists in this range
and for the negative companion in Section~\ref{sec:select-neg-minus}.
A fixed rational weight
\(R(t,M,J,Q)\) is admissible if it has a presentation with
denominator nonzero on every point of the full finite triangle
for all large \(N\). The module is the left \(\D\)-span of the
embedded field and the sequences
\[
 N\longmapsto\sum_{0\le n\le j<N}
 \mathcal C_j(N)K_{jn}q^{\ell n(n+1)}R(t,t^{2N},q^j,q^n)
\]
over admissible \(R\). Its span consists of finite combinations
with fixed rational coefficients and fixed signed shifts.
The estimates may depend on that fixed element.
The weight \(\mathcal R_\eta/\Delta\) in
\eqref{eq:physicalforcing} is admissible by the pole and support
checks there.

Choose a fixed Laurent-polynomial numerator and denominator for
the weight. After evaluation, their exponents come from a finite
list of affine expressions in \(N,j,n\) and lie in
\([-O_R(N),O_R(N)]\) on the triangle. Cancellation only deletes
terms, and admissibility keeps the denominator nonzero.
The unweighted summand valuation is
\[
 \frac{j(j+1)}2-jN+j+\Big(\ell+\frac12\Big)n(n+1).
\]
If \(\ell=-h<0\), the minimum is at \(n=j=N-1\)
and equals \(-(N-1)(hN-1)\).
If \(\ell>0\), the minimum in \(n\) is at zero. The
remaining expression attains its minimum
\(-(N-1)(N-2)/2\) at \(j=N-2,N-1\) for \(N\ge2\).
The valuation inequality for sums permits cancellation at these
tied terms. Multiplication by a fixed coefficient, fixed shifts,
and finite combinations preserve the resulting lower bound
\[
 \nu(w_N)\ge
 \begin{cases}
  -hN^2-O_w(N),&\ell=-h<0,\\
  -N^2/2-O_w(N),&\ell>0
 \end{cases}
 \qquad(w\in\mathcal W_\ell).
\]
In \eqref{eq:actual-eta-plus}, the leading quadratic coefficient
is \(-a=-(m+h)\) in the first case and \(-a\le-1\) in the
second. In both cases it is strictly below the forcing bound.
This proves \(F\notin\mathcal W_\ell\), including at \(a=1\),
where the actual valuation is \(-(N-1)^2\).

\subsection{Infinity growth for the negative companion at \(a<0\)}\label{sec:select-neg-minus}

For \(\eta=-1\) and \(b=m+\ell=-a\ge1\), use
\(\nu_\infty(f)=\nu(f(t^{-1}))\).
Factor by factor in \eqref{eq:universalbailey},
\[
 \mathcal C_j(N;q^{-1})=q^{-j(j+1)} \mathcal C_j(N;q),\qquad
 K_{jn}(q^{-1})=(-1)^j q^{j(j-1)/2}K_{jn}(q).
\]
The second identity comes from
\((q^{-1};q^{-1})_k=(-1)^kq^{-k(k+1)/2}(q;q)_k\).
The inverted actual input is \(q^{-b n(n+1)}V_n(q)\).
After cancellation of the finite kernel powers, the summand
valuation is \(-jN+(\tfrac12-b)n(n+1)\).
For fixed \(j\) its unique minimum is at \(n=j\).
The successive \(j\)-difference there is \(-N-(2b-1)j<0\),
so \(j=N-1\) is the unique minimum. The resulting value and
leading coefficient are
\[
 \nu_\infty(F_N)=-(b+\tfrac12)N(N-1),\qquad
 \operatorname{lc}_\infty(F_N)=(-1)^{N-1}.
\]

In the weighted module \(\mathcal W_\ell\) the inverted
unweighted exponent is instead \(-jN+(\tfrac12-\ell)n(n+1)\).
For \(\ell\le0\) its minimum is at \(n=0,j=N-1\), for
\(\ell\ge1\) at \(n=j=N-1\).
Inverted fixed weights still contribute only \(O_R(N)\),
by the same affine exponent bound. Thus every fixed
\(w\in\mathcal W_\ell\) has
\[
 \nu_\infty(w_N)\ge-\max(1,\ell+\tfrac12)N^2-O_w(N).
\]
Here \(b+\tfrac12>\max(1,\ell+\tfrac12)\), since \(b\ge1\)
and \(b-\ell=m\). This gives the strict separation for the
negative companion in this range.

\subsection{The finite Gaussian forcing module}\label{sec:select-gaussian}

The remaining ranges are \(\eta=1,\ a\le0\) and
\(\eta=-1,\ a\ge0\). Put \(\gamma=|\ell|\) in these ranges.
Let \(w_n^{(\ell,c)}=q^{\ell n(n+1)}(Q^c-\widehat Q^c)\)
be the half-interval input, and let \(G_c^{(\ell)}\) be the output
of \eqref{eq:universalbailey} on that input. We call these finite
transforms Gaussian periods. They satisfy
\[
 G_0^{(\ell)}=0,\qquad G_{-c}^{(\ell)}=-q^cG_c^{(\ell)}.
\]
The reflection identity follows from \(\widehat Q=q^{-1}/Q\).
For \(z=0,1\) write
\(J^z=\sum_{i=0}^d d_i^{[z]}\Delta R_i\) and let
\(D^{[z]}=\sum_i d_i^{[z]}u^i\).
This operator implements the weight \(J^z/\Delta\) on any
fixed input by the same cleared color-shift relation used
in \eqref{eq:physicalforcing}.
Put \(p_i=\ell+\eta\lambda_i=\eta m+\ell+\eta\beta_i\).
Expanding \(A(Q,J)\) in \eqref{eq:innerforcing}, with its
reflected subtraction, gives
\begin{equation}\label{eq:finite-module-force}
 P_aF=-q^{\eta m+\ell}\sum_{i=1}^4 c_i(q^\eta)
 \left(D^{[0]} G_{p_i}^{(\ell)}-D^{[1]}G_{p_i-1}^{(\ell)}\right).
\end{equation}
Here the \(Q\)-terms are precisely \(Q^{p_i}-JQ^{p_i-1}\).
For \(\eta=1\), \(\ell=\gamma=m-a\ge m\), so
\(1\le p_i\le2\gamma-1\) and \(0\le p_i-1\le2\gamma-2\),
using \(m-r-s\ge1\).
For \(\eta=-1\), \(\ell=-\gamma=-m-a\), so
\(-p_i=\gamma+\lambda_i\), and both \(-p_i\) and \(-p_i+1\)
are in \([1,2\gamma]\).
Thus in each case \(P_aF\) belongs to the left \(\D\)-module
generated by the following \(2\gamma\) periods
\begin{equation}\label{eq:finite-U}
 \mathcal U_\ell=\sum_{c=1}^{2\gamma}\D G_c^{(\ell)} .
\end{equation}
The identity uses only the displayed periods and the two
pole-basis weight operators.

\subsection{Controlled coefficients from the finite clasp sum}\label{sec:select-coefficients}

Coefficient growth separates the satellite sequence from \(\mathcal U_\ell\)
once equality on a tail forces equality of expansion coefficients. At any
fixed cutoff, the remainder estimates give higher order than the retained
nonzero terms for large \(N\).

Fix \(\theta\in\{0,1\}\) and a sequence \(w=(w_n)_{n\ge0}\)
with values in \(\Q((q))\), independent of \(N\).
Treat \(X\) as a separate formal variable. To form the actual
output, first restrict the sum to the finite triangle \(i+n<N\),
then set \(X=q^N\).
\begin{equation}\label{eq:coefficient-kernel}
 Y_{\theta,N}(w)=\left.\frac{1-q}{1-X}
 \sum_{\substack{i,n\ge0\\ i+n<N}}
 X^i(-1)^{\theta n}q^{\theta n(n+1)/2}w_n
       \frac{(q^{i+1};q)_n(Xq^{-i-n};q)_n}{(q;q)_n}\right|_{X=q^N}.
\end{equation}
To form \(A_{\theta,d}(w)\), allow all \(i,n\ge0\) and take
\([X^d]\) with \(1/(1-X)\) expanded at zero and \(X\) still formal.
The lemma proves coefficientwise convergence in \(\Q((q))\) and
compares fixed truncations of \(\sum_{d\ge0}A_{\theta,d}(w)X^d\)
with the finite output.
For the output \(F\) on input \(v\), reverse the \(i\)-sum in
\eqref{eq:claspfinite} to get \(F_N=q^{N-1}Y_{0,N}(v)\).
Inverting that sum instead, set
\((\mathcal I F)_N=q^{N(N-1)}F_N(t^{-1})\).
The prefactor in \eqref{eq:claspfinite} contributes
\(q^{-N(N-1)}(1-q)/(1-X)\) after inversion, and the ratio of
Pochhammers contributes
\((-1)^n q^{-nN+n(n+1)/2}\) times the ratio in
\eqref{eq:coefficient-kernel}. The inverted monomial contributes
\(X^i q^{Nn}\), giving \((\mathcal I F)_N=Y_{1,N}(v(q^{-1}))\)
on the same finite triangle.

\begin{lemma}[Finite-cutoff estimates]\label{lem:coefficient-cutoff}
Let \(w_n\in\Q((q))\) be fixed independently of \(N\).
For \(\theta=0\) assume \(\nu(w_n)\ge n^2-Bn-C\)
with fixed \(B,C\). For \(\theta=1\) assume \(\nu(w_n)\ge-C\)
with fixed \(C\). In either case the sums defining
\(A_{\theta,d}(w)\) converge in \(\Q((q))\) for every fixed \(d\ge0\).
For every fixed \(H\ge1\),
\begin{equation}\label{eq:coefficient-remainder}
 \nu\left(Y_{\theta,N}(w)-\sum_{d=0}^{H-1} A_{\theta,d}(w)q^{dN}\right)
 \ge HN-O_H(1)
\end{equation}
as \(N\) increases.
\end{lemma}
\begin{proof}
Write \(d=i+j+e\), with \(j\) the degree from the second
Pochhammer and \(e\) the distinct geometric degree from
\((1-X)^{-1}\). For \(j\le n\) the least contribution of
the Pochhammer comes from the unique subset of indices
\(0,\ldots,j-1\). It has exponent
\(-j(i+n)+j(j-1)/2\) and sign \((-1)^j\).
The other factor \((1-q)(q^{i+1};q)_n/(q;q)_n\) has
order zero and leading coefficient one. Thus the exact order
of a nonzero term in \(A_{\theta,d}(w)\) is
\begin{equation}\label{eq:coefficient-valuation}
 \theta\,n(n+1)/2-j(i+n)+j(j-1)/2+\nu(w_n).
\end{equation}
This grows quadratically in \(n\) for fixed \(d\) under either
hypothesis. Hence the auxiliary coefficient sums converge.

\emph{High degrees on the finite triangle.} The terms of degree \(d\ge H\)
contribute the exponent
\[
 Nd-j(i+n)=i(N-j)+j(N-n)+eN
\]
in addition to the remaining terms of
\eqref{eq:coefficient-valuation}.
For \(\theta=0\), this is at least \(H(N-n)\) if \(j\ge H\).
If \(j<H\), the bound \(i+e\ge H-j\) gives
\(HN-j(H-j)-jn\ge HN-H^2-Hn\).
Adding the assumed quadratic in \(n\) gives \(HN-O_H(1)\).
For \(\theta=1\) the evaluated order is
\[
 iN+j(N-i-n)+eN+n(n+1)/2+j(j-1)/2+\nu(w_n).
\]
If \(i\ge H\) this is at least \(HN-C\).
Otherwise \(j+e\ge H-i\) gives the lower bound
\(HN-(H-i)(i+n)+n^2/2-C\), again \(HN-O_H(1)\).
The geometric series converges for \(N>0\).

\emph{Omitted terms outside the triangle.} These are the terms with
\(d<H,\ i+n\ge N\).
For \(\theta=0\) their unevaluated order is at least
\(n^2-(B+H)n-C-H^2\) with \(n\ge N-H+1\).
For \(\theta=1\), rewrite \eqref{eq:coefficient-valuation}
as \((n-j)(n-j+1)/2-ij+\nu(w_n)\).
Here \(n-j\ge N-H+1\), so the order is at least
\((N-H)^2/2-H^2-C\).
In each case the bound is \(HN-O_H(1)\), even before adding
\(dN\). The omitted auxiliary tails converge quadratically.
The two estimates prove \eqref{eq:coefficient-remainder}.
\end{proof}

\emph{Controlled expansions on a tail.}
An expansion \(\sum_{j\ge j_0}c_j(t)M^j\), with
\(j_0\in\Z\) and \(c_j(t)\in\Q((t))\), is controlled if,
for every fixed \(J_0>j_0\), subtracting its terms \(j<J_0\)
evaluated at \(M=t^{2N}\) from the actual sequence leaves a
remainder of \(t\)-order at least \(2J_0N-O_{J_0}(1)\).
The constant is for that cutoff and that sequence. At a first
differing coefficient the evaluated term has smaller eventual
slope than all later retained terms and the remainder at a
larger fixed cutoff. Thus equality on a tail implies equality
of the coefficients. We call a tail \emph{flat} for its
controlled expansion if all the coefficients vanish. Its order
then exceeds any fixed linear slope eventually.

The lemma gives such an expansion for the output before
normalization at \(\theta=0\), namely
\(t^{-4}\sum_d A_{0,d}(w)(t^4)M^{2d+2}\).
At \(\theta=1\) the normalized inverted output has expansion
\(\sum_d A_{1,d}(w)(t^4)M^{2d}\).
The intermediate \(M\)-cutoffs follow by moving finitely many
terms into the remainder.

\emph{Rational coefficients and shifts.}
A rational coefficient \(r(t,M)\) has a controlled Laurent
expansion. Its fixed Taylor remainder is a power \(M^{J_0}\)
times a rational function regular at \(M=0\), whose valuation
is bounded below on the tail. Products of a rational coefficient
and a controlled expansion are controlled: for target cutoff
truncate each factor at \(J_0\) minus the other's lowest degree.
The omitted products have the required order, using the
linear lower bounds on the actual tails from coarse truncations.
For a flat tail the same argument uses any needed cutoff. A fixed shift \(u^k\), \(k\in\Z\), multiplies each
coefficient of \(M^j\) by \(t^{-2kj}\).
These observations prove closure under the full rational field
and finite Ore operations, including cancellations.

\emph{Bounds after rational multiplication.} The Taylor coefficients of any fixed rational function have a
linear lower valuation bound. After removing its Laurent power,
normalize its denominator to \(1+\sum_{i=1}^s b_i(t)M^i\)
and choose \(C'\) with \(\ord_{t=0}b_i\ge-C'i\).
The inverse coefficients satisfy the triangular convolution
recursion and therefore have order at least \(-C'j\).
A finite numerator changes only a linear bound.
If a period's coefficients have order at least
\(-\lambda j^2-O(j)\) for fixed \(\lambda>0\), multiplication
by a rational coefficient preserves this quadratic coefficient.
Indeed the convolution is finite and the period index is at
most \(j\) plus a constant. Fixed shifts incur linear loss,
and finite sums can only raise the minimum order.

\emph{Normalized inversion.} The exact intertwining
relations are
\begin{equation}\label{eq:inversion-shift}
 \mathcal I(rF)=r(t^{-1},M^{-1})\,\mathcal IF,\qquad
 \mathcal I(u^kF)=q^{-k(k+1)} M^{4k}u^k\mathcal IF.
\end{equation}
These follow by subtracting \((N-k)(N-k-1)\) from \(N(N-1)\).
Thus the same closure holds for the normalized inversion of the
finite module \eqref{eq:finite-U}.

\subsection{The positive companion at \(a=0\) and at \(a<0\)}\label{sec:select-pos-nonpos}

At \(a=0\), the elementary lowest-term argument of
\eqref{eq:actual-eta-plus} gives tied outer terms with opposite
signs. The reindexed clasp formula instead gives the exact
leading terms for the actual sequence. In \eqref{eq:coefficient-kernel}
with \(\theta=0\), every \(i=0\) summand with \(v_n=V_n\)
has valuation zero and leading coefficient one. Indeed at
\(X=q^N\) all the Pochhammer exponents in the finite support are
positive. All \(i>0\) terms have greater valuation. Hence
\[
 \nu(F_N)=N-1,\qquad \operatorname{lc} F_N=N .
\]
For \(\ell=m\), the input \(w_n^{(\ell,c)}\) has valuation
\(m n^2+(m-c)n-c\), so Lemma~\ref{lem:coefficient-cutoff}
applies to every generator with \(\theta=0\).
A nonflat controlled expansion has a fixed leading coefficient
on a tail, by its first nonzero \(M\)-coefficient.
A flat tail has \(t\)-order greater than \(4N-4\)
for all sufficiently large \(N\). Thus neither can represent the displayed
actual sequence, proving \(F\notin\mathcal U_m\).

For \(a=-b<0\), \(\ell=\gamma=m+b\), apply the same cutoff
lemma with \(\theta=0\) to \(w_n=q^{b n(n+1)}V_n\) and
to every period input. Their valuations are bounded below
by \(n^2-Bn-C\) with fixed constants.
Set \(n=j+l\), \(i=d-j-e\) in \eqref{eq:coefficient-valuation}.
For the actual input the order is
\[
 (b+\tfrac12)j^2+(b-\tfrac12-d)j+ej
                  +(2b-1)jl+b l^2+b l .
\]
For \(d=(2b+1)k\), \(k\ge1\), the \(j\)-part differs from
its value at \(j=k\) by \(h((2b+1)h+2b-1)/2\), \(h=j-k\).
This is positive for every nonzero integer \(h\).
The remaining nonnegative terms vanish at the minimum only for
\(l=e=0\). Hence the unique index is \(n=j=k,\ i=2bk\).
The unit \(V_k\) and the least Pochhammer subset give sign
\((-1)^k\) and exact order
\begin{equation}\label{eq:positive-companion-coeff}
 \nu A_{0,(2b+1)k}(w)
   =-\big((2b+1)k^2-(2b-1)k\big)/2.
\end{equation}

For comparison, the valuation of \(w_n^{(\ell,c)}\) is
\(\gamma n^2+(\gamma-c)n-c\). In the same variables the order
of a contributing term is
\[
\begin{aligned}
 &(\gamma+\tfrac12)j^2+(\gamma-c-\tfrac12-d)j+ej\\
 &\hspace{2em}+(2\gamma-1)jl+\gamma l^2+(\gamma-c)l-c .
\end{aligned}
\]
Discarding the nonnegative \(ej,(2\gamma-1)jl\) and
\((\gamma-\tfrac12)j\), square completion gives
\[
 \nu A_{0,d}\bigl(w^{(\ell,c)}\bigr)\ge
 -\frac{(d+c)^2}{2(2\gamma+1)}
 -c-\frac{\max(c-\gamma,0)^2}{4\gamma}.
\]
Write \(d_M=2d+2\) for the \(M\)-degree. The output
coefficients have \(t\)-order at least
\(-d_M^2/[2(2\gamma+1)]-O(d_M)\), and quantitative closure
gives the bound for every fixed element of
\(\mathcal U_\ell\). At \(d_M=2(2b+1)k+2\) the actual
\(t\)-order from \eqref{eq:positive-companion-coeff}, including
the common \(t^{-4}\) prefactor, is
\(-2(2b+1)k^2+2(2b-1)k-4\).
The quadratic gap in \(k\) is
\(4(2b+1)(\gamma-b)/(2\gamma+1)>0\).
Uniqueness of controlled coefficients now excludes \(F\)
from \(\mathcal U_\ell\).

\subsection{The negative companion at \(a\ge0\)}\label{sec:select-neg-nonneg}

Here \(\ell=-\gamma=-m-a\). In \(\mathcal I F\) the inverted
actual input is \(w_n=q^{a n(n+1)}V_n(q)\).
For the inverted periods it is
\[
 w_n^{(c)}=q^{\gamma n(n+1)}
                 (q^{-cn}-q^{c(n+1)}),\qquad 1\le c\le2\gamma .
\]
The normalized inverted periods are \(Y_{1,N}(w^{(c)})\).
Both types of input satisfy the \(\theta=1\) cutoff hypothesis.
With \(n=j+l\) the actual coefficient order becomes
\[
 (a+1)j^2+(a-d)j+ej+2ajl+(a+\tfrac12)l(l+1).
\]
At \(d=2(a+1)k\), \(k\ge0\), the \(j\)-part differs from
its value at \(k\) by
\((j-k)^2+a(j-k)(j-k+1)\), positive for any nonzero
integer \(j-k\). Equality in the other terms forces
\(n=j=k, e=0\), with \(e=0\) at \(k=0\) from the degree
restriction. Thus \(i=(2a+1)k\), and the signs
\((-1)^n\) and \((-1)^j\) cancel. This proves
\begin{equation}\label{eq:negative-companion-coeff}
 \nu A_{1,2(a+1)k}(w)=-(a+1)k^2+ak,\qquad
 \operatorname{lc}A_{1,2(a+1)k}(w)=1 .
\end{equation}
The period input has order
\(\gamma n(n+1)-cn\ge\gamma n(n-1)\ge\gamma j(j-1)\).
Using \(i\le d-j\), the coefficient order is at least
\((\gamma+1)j^2-(d+\gamma)j\), giving
\[
           \nu A_{1,d}(w^{(c)})\ge-\frac{(d+\gamma)^2}{4(\gamma+1)} .
\]
At \(M\)-degree \(d_M=2d\) this is
\(-d_M^2/[4(\gamma+1)]-O(d_M)\) in \(t\)-order.
The bound holds for every fixed element of
\(\mathcal I\mathcal U_\ell\) by \eqref{eq:inversion-shift}
and the finite closure above. Along \(d_M=4(a+1)k\) the actual
order is \(-d_M^2/[4(a+1)]+ad_M/(a+1)\).
The positive gap
\[
 \frac{1}{4(a+1)}-\frac{1}{4(\gamma+1)}
 =\frac{m}{4(a+1)(\gamma+1)}
\]
overcomes every fixed linear loss, proving \(F\notin\mathcal U_\ell\).
This argument includes the boundary values \(a=0,1\) for the
negative companion. The positive companion at these values uses,
respectively, its varying leading coefficient and the color valuation
\eqref{eq:actual-eta-plus}.

\subsection{The selected right divisor}\label{sec:selected-divisor}

In all cases \eqref{eq:physicalforcing} or
\eqref{eq:finite-module-force} puts \(P_aF\) in the corresponding
module and the estimates give nonmembership of \(F\).
Lemma~\ref{lem:relative} supplies a monic selected right divisor \(R_a\).
For \(a\ne0\), irreducibility makes it the monic form of \(P_a\).
For \(a=0\) the class is nonzero and the row in the satellite
color has no linear right divisor, so the degree of \(R_0\) is two or three.
The monic Gauss argument preserves its degree and endpoints,
and \eqref{eq:Eb} implies \(E_{\rm back}\mid\overline R_0\).
Thus in both situations every actual reduced-tail annihilator
is right-divisible by \(R_a\) with regular nonzero residue
containing \(E_{\rm back}\). The degree computation concerns the
class modulo the chosen forcing module.

\paragraph{The unknot companion.}
For \(r=1\), the companion is the unknot. In the negative-clasp
marking the satellites are the pure closures \(Z_\ell\).
The input of~\eqref{eq:commoninput} and the quantum dimension give
\[
 v_n=q^{\ell n(n+1)}[2n+1],\qquad
 [2n+1]=\frac{q}{q-1}(Q-\widehat Q),\qquad
 F=\frac{q}{q-1}G_1^{(\ell)} .
\]
The last identity uses the same finite transform as
\eqref{eq:universalbailey}, with the half-interval input of
Section~\ref{sec:select-gaussian}. The weight
\(Q-\widehat Q=Q-(qQ)^{-1}\) is admissible on the whole finite
triangle. Thus \(F\in\mathcal W_\ell\). For \(|\ell|\ge1\),
\(G_1^{(\ell)}\) also belongs to the list of generators of
\(\mathcal U_\ell\) in~\eqref{eq:finite-U}. At \(\ell=0\),
the Hopf-belt identity of Section~\ref{sec:family} at dimension
\(p=1\) has \(\mathcal C_j(1)=0\) for \(j\ge1\) and \(b_0=1\), giving
\(F_N=1\) for all \(N>0\). In each relevant quotient the class of
\(F\) is zero. For the nontrivial twist knots, L\^e's AJ theorem
\cite[Cor.~0.2]{Le2006} and the resultants of
Section~\ref{sec:eexclusion} give the absence of \(E\) from
the specialized emph{minimum}; for \(Z_0=U\) this follows
from the first-order recurrence of \([N]\), whose primitive
specialization is associated to \(L-1\). Reflection gives the
corresponding statement for \(E^\vee\) at positive clasp.

\section{Classical exclusion in the preferred coordinates}\label{sec:torusclass}

To exclude \(E\) on the classical side, it suffices to find a
polynomial that vanishes on the generic outer boundary values and
is coprime to \(E\). Restriction to the companion forces two
fillings of the marked pattern, which give such a polynomial.
The next lemma passes from generic vanishing to divisibility by the
full horizontal \(A\)-polynomial.

\subsection{The full boundary polynomial and its horizontal factors}\label{sec:boundary}
The classical family argument gives a containment bound off finitely many
meridian fibers. The following lemma selects the full affine component
equation before passing to its horizontal localization, where that generic
bound implies divisibility.

Fix the exterior $X$ of a knot in $S^3$, a finite group presentation and words for its positive meridian and preferred longitude. The $\SL(\C)$ representation equations and the restriction of characters are defined over $\Q$: matrix multiplication, determinant one, inversion and the traces of the peripheral words use integer coefficients. Write $\mathbb T=(\C^\times)^2$. Let $S_X\subset\mathbb T$ be the set of common-eigenline pairs $(M,L)$ of peripheral meridian and preferred longitude matrices in representations of the knot exterior $X$. Closures are Zariski closures.

\begin{lemma}[Generic boundary containment]\label{lem:aex}
If a nonzero $F\in\Q[M^{\pm1},L^{\pm1}]$ vanishes on $S_X$ off finitely many meridian fibers, the union of the one-dimensional components of the affine closure of $S_X$ has a primitive integral equation $\widetilde A_X(M,L)$ with no repeated irreducible factors. Its horizontal localization $A_X$, the image of $\widetilde A_X$ in $\Q(M)[L^{\pm1}]$, divides $F$. For a knot in $S^3$ this equation includes $L-1$.
\end{lemma}
\begin{proof}
We first identify the curve components of the boundary image.
Each horizontal component has a dense part away from the
exceptional fibers, and vertical factors become units upon
localization in \(M\). This will give the divisibility.

The diagonal eigenvalue torus maps finitely onto the boundary character variety, with fibers the orbits of simultaneous inversion. If a peripheral matrix has distinct eigenvalues the commuting pair has two common eigenlines. If neither matrix does, the single eigenvalue pair is fixed by inversion, and simultaneous triangularization still supplies a common eigenline. Thus over the actual restriction image the full fiber consists of the pairs in $S_X$. The closure of an inversion-stable set is stable, and the finite quotient is closed. Its closure is consequently the full inverse image of the closure of its image.

Choose a nonzero complex polynomial $d(M)$ vanishing on the finitely many exceptional meridians. At this stage only an upper bound on dimension is needed. In the family chart $M^4-1$ itself suffices for the exceptional nonzero meridians. Then $dF$ vanishes on the full set, excluding two-dimensional components. Among the finitely many irreducible character components, those with zero-dimensional restriction image supply at most finitely many points. Those with one-dimensional image give after lifting and closing precisely the curve components. This is the construction of Cooper, Culler, Gillet, Long and Shalen~\cite[Secs.~2.2--2.3]{CCGLS}, with coordinates written $(M,L)$. First take closure in $\mathbb T$, then close each selected torus curve in $\C^2$. Its torus locus is dense, so this last step adds no axis component. For an irreducible affine plane polynomial $g$, its zero set lies in $\overline{S_X}^{\,\mathrm{aff}}$ exactly when it is one of these curve components: an irreducible curve in a finite union of curves and isolated points must equal one of the curves. This identifies which factors occur in $\widetilde A_X$. Field automorphisms over $\Q$ preserve the representation equations, restriction and dimension selection. The reduced equation descends up to scalar to a primitive polynomial over $\Z$, as in the proposition in \cite[Sec.~2.3]{CCGLS}. Abelianization realizes $(M,1)$ for every $M\ne0$; the reduced product includes that line once \cite[Sec.~2.5]{CCGLS}.

The complement of the exceptional fibers is dense in each horizontal component. Vanishing there makes its irreducible factor divide $F$ over $\C(M)$. The distinct factors multiply, and vertical factors become units. Monic polynomial division over $\Q(M)$ gives the same zero remainder after scalar extension, proving divisibility over $\Q(M)$. An exceptional factor introduced in $d$ bounds the image. Equality requires selection of the components that actually occur.
\end{proof}

\subsection{The two companion fillings}

For the right-handed torus knot the group and its preferred pair
are~\cite[Sec.~2.2]{DuboisKashaev}
\[
 \langle A,B\mid c_{\rm fib}=A^r=B^s\rangle,\qquad
 \mu_C=A^vB^w,\quad sv+rw=1,\qquad
 \lambda_C=c_{\rm fib}\mu_C^{-m}.
\]
The regular fiber has linking \(m=rs\), accounting for the
longitude correction. Under a companion mirror the meridian
inverts and the longitude is transported, giving
\(\lambda_C=c_{{\rm fib},\eta}\mu_C^{-\eta m}\), with \(c_{{\rm fib},\eta}\) central.
An abelian companion image has longitude \(I\).
In any nonabelian image, including a reducible nonabelian one,
\(c_{{\rm fib},\eta}\) maps to \(\varepsilon I\) with
\(\varepsilon=\pm1\): the commutant of any nonscalar
\(2\times2\) matrix is commutative. This yields the necessary
filling conditions from \eqref{eq:twistgluing}:
\[
 \mu_A\lambda_A^{-\ell}=I\quad\text{or}\quad
 \mu_A\lambda_A^{\,\eta m-\ell}=\varepsilon I .
\]
In these equations matrices denote their images under the
representation. Multiply the pattern representation by the
sign character with value \(\varepsilon\) on the axis meridian
and \(1\) on the outer meridian. Both preferred longitudes have
zero homology by zero linking, so this preserves both and fixes
the entire outer pair. It makes the second relation a pure
filling. Hence the two signed indices in \(\mu_A\lambda_A^k\)
are \(-\ell\) and \(a=\eta m-\ell\).

\subsection{A containing equation for a signed filling}\label{sec:word}

For a representation \(\rho\) of the pattern group, let
\(X=\rho(x)\) and \(Y=\rho(y)\) be the matrices of Schilder's
positive meridian generators \(x,y\) from Section~\ref{sec:mark}.
Here the capital letters name matrices of those generators.
Put
\[
 Q_0=X^{-1},\quad P=Y^{-1}Q_0Y,\quad
 W=P Q_0^{-1}P^{-1}Q_0 .
\]
In the source words of Section~\ref{sec:mark} the axis longitude
evaluates to \(W\) and the outer longitude to
\(P^{-1}Q_0YQ_0^{-1}Y^{-1}Q_0\). The reflection matrix
\(-D\) on the adapted axis sends the signed filling to
\(YW^{-k}=I\). Thus \(Y=W^k\) and
\(PW^{-k}=W^{-k}Q_0\), with preferred outer meridian \(Q_0\).
At \(k=0\) this gives \(Y=I,\ P=Q_0\), and the longitude word
is \(I\). If \(P,Q_0\) share an invariant line, the original meridians
\(Y=W^k\) and \(X=Q_0^{-1}\) both preserve it. The entire pattern
representation then preserves that line. The longitude has eigenvalue \(1\)
on that line by zero peripheral homology. For generic \(M\),
the other eigenvalue of the longitude is also \(1\).
This includes reducible pattern restrictions.

Otherwise, away from \(M^2=1\), choose the \(M\)-eigenline
of \(P\) and the \(M^{-1}\)-eigenline of \(Q_0\).
They are distinct. Normalize the upper off-diagonal entry
of \(P\), which is nonzero by irreducibility of this pair.
Then
\[
 P=\begin{pmatrix}M&1\\0&M^{-1}\end{pmatrix},\qquad
 Q_0=\begin{pmatrix}M&0\\\zeta&M^{-1}\end{pmatrix}.
\]
Write \(h=M-M^{-1}\) and \(w=h+\zeta/M\). Direct multiplication
gives
\[
 W=\begin{pmatrix}(1-\zeta)^2+M^2\zeta&w\\
                  \zeta w&1+\zeta/M^2\end{pmatrix}.
\]
Cayley--Hamilton writes \(Y=\alpha W+\beta I\) for every signed
\(k\). A scalar \(Y\) would force \(P=Q_0\), so \(\alpha\ne0\).
The upper-right entry of \(YP-Q_0Y\) makes \(Y\) proportional
to \(W-cI\), where \(c=W_{11}-hW_{12}\).
This latter matrix is invertible. Its first column is
\(w(h,\zeta)^t\), so \(w\ne0\).
Use \(x=M^2\) and the coordinates \(J,Z,e\) of
Section~\ref{sec:symbol}. On the \(M\)-eigenline of \(Q_0\), generated by \((h,\zeta)^t\),
multiplication of the longitude word gives
\[
 L=\frac{h-M\zeta}{h+\zeta/M},\qquad
 \zeta=\frac{(x-1)(1-L)}{L+x}=1-J .
\]
These equations use both preferred eigenvalues on that common
eigenline. At \(L=-x\) the cross-multiplied relation forces
\(x^2=1\). The exceptional meridian-only defining factors become units
under horizontal localization.

For \(r_0(U)=hU_{12}+U_{22}\) the same matrix gives
\[
 \tr W=2+\zeta^2+h^2\zeta=Z,\qquad
 r_0(I)=1,\qquad r_0(W^{-1})=Z-1+e .
\]
The upper-right entry of \(PW^{-k}-W^{-k}Q_0\) is
\(r_0(W^{-k})\). Cayley--Hamilton in both directions therefore
gives \(T_k=r_0(W^{-k})=0\).
The initial values and recurrence are exactly those of the
\(T_k\) in Section~\ref{sec:symbol}. The same trace polynomials
thus control the classical fillings and the commuting symbol
\eqref{eq:commonPi}, which also contains the factor \(Z+2\)
corresponding to \(E\).
The positive-index representative is \(B_k\), and the
negative-index representative is \(C_{-k}\), as defined in
\eqref{eq:B}. Their values at \(L=-x\) from
Section~\ref{sec:twistirr} control the clearing chart.
Thus, putting
\[
 \mathcal P_k=
 \begin{cases}B_k,&k>0,\\1,&k=0,\\C_{-k},&k<0,\end{cases}
\]
all generic actual outer values lie on the following
containing equation:
\begin{equation}\label{eq:classicalcontainment}
 A_{Z_\ell(T_\eta)}
       \ \mid\ (L-1)\mathcal P_{-\ell}\mathcal P_{\eta m-\ell}.
\end{equation}
The divisibility follows from Lemma~\ref{lem:aex}. It includes
the abelian component and does not require extending generic
characters back across the companion torus.

\subsection{The resultants with \(E\)}\label{sec:eexclusion}

A nonzero longitude resultant excludes \(E\) as a common
irreducible component. Work in the quadratic longitude algebra
\(\Q(M)[L^{\pm1}]/(E)\), whose trace and norm are the sum and
product of the values on the two longitude roots. The meridian
trace coordinate here is \(K=\tr(Q_0^2)=M^2+M^{-2}\).
On \(E=0\), equation \eqref{eq:pullbacks} gives \(Z=-2\) and
\(J^2-KJ+K+3=0\), hence
\[
 \operatorname{Tr}(e)=2-K,\qquad
 \operatorname{Nm}(e)=4,\qquad
 \operatorname{Nm}(L+x)=x(x-1)^2 .
\]
For the first two identities use \(e=1+J-K\), and for the
last use the sum and product of the two roots of \(E\).
At \(Z=-2\) the trace recurrence gives
\(T_k=(-1)^k(2k+1-ke)\) for all integers \(k\).
Since \(E\) is monic, the resultant is the norm of the
other polynomial. Consequently
\begin{align}
 \Res_L(E,B_k)
 &=x^{4k-1}(x-1)^{4k}
      [k(2k+1)(x+1)^2+x],\label{eq:allbresultant}\\
 \Res_L(E,C_k)
 &=x^{4k-2}(x-1)^{4k-2}
      [k(2k-1)(x+1)^2+x]\qquad(k\ge1).\label{eq:allcresultant}
\end{align}
Indeed the relevant norms of \(T_k,T_{-k}\) are
\(1+k(2k+1)(K+2)\) and \(1+k(2k-1)(K+2)\).
The bracketed polynomials equal \(-1\) at \(x=-1\).
Together with \(E(M,1)=4x\), the nonzero resultants imply
\(\gcd_R(E,A_S)=1\) in \eqref{eq:classicalcontainment},
for both companion signs and all signed filling indices.

\subsection{Distinct untwisted positive-companion knots}

For \(\ell=0\) and \(\eta=1\), the gluing torus separates the
hyperbolic Whitehead link exterior from a Seifert exterior
over a disk with two exceptional fibers of orders \(r,s\)
\cite[(1.4)--(1.6)]{Zheng}. Both pieces are irreducible
with incompressible boundary, and neither is a torus product.
The hyperbolic side is atoroidal. The union cannot be Seifert
fibered, since cutting a Seifert space along an essential torus
would give Seifert pieces, including the torus or interval
bundles of a horizontal cut. Thus the JSJ decomposition
identifies the unique Seifert side under an exterior
homeomorphism~\cite[Secs.~1.2--1.3, Thm.~1.9]{Hatcher}.
Its group modulo its center is \(C_r*C_s\). Indeed the regular
fiber generates an infinite cyclic central subgroup, with the
displayed quotient having trivial center. A finite subgroup of
the free product fixes a vertex in the coset tree, so the
maximal finite subgroups recover the conjugacy classes of the
two factors and their orders. The Seifert side therefore
recovers \(\{r,s\}\), and the order \(r<s\) recovers the pair.

\subsection{The companion-abelian filling factor}\label{sec:pure-filling}

To identify the actual boundary factor of the pure closure
\(Z_{-k}\), one can reverse the matrix calculation of
Section~\ref{sec:word}. At \(k=0\) the filling gives \(L=1\),
so \(A_{Z_0}\doteq L-1\). For \(k\ne0\), take a generic
point of \(\mathcal P_k=0\) and form \(P,Q_0,W\) with
\(\zeta=1-J\) as above. In these coordinates
\(T_k=r_0(W^{-k})=0\). The pair \(P,Q_0\) can share a line
only for \(\zeta=0\) or \(\zeta=-h^2\). At those values the
trace recurrence has \(Z=2\) and gives \(T_k=1-kh^2\) and
\(T_k=1\), respectively. Neither vanishes at generic \(M\).

Set \(Y=W^k\) and \(X=Q_0^{-1}\). The matrix \(U=W^{-k}\)
lies in the span of \(I,W\), so \(U_{21}=\zeta U_{12}\).
For this form of \(U\), the diagonal entries of \(PU-UQ_0\)
vanish and its off-diagonal entries are \(r_0(U)\) and
\(-\zeta r_0(U)\). It follows that \(P=Y^{-1}Q_0Y\).
After this substitution \(W\) is the image of Schilder's
axis longitude word \(\lambda_y\). The pattern relator of
Section~\ref{sec:mark} is a cyclic conjugate of the inverse
of the cyclically reduced commutator \([y,\lambda_y]\),
so it holds because \(Y=W^k\). The filling relation holds
by construction. Evaluating the outer word now gives the
same \(L=(h-M\zeta)/(h+\zeta/M)\) as before.
By Section~\ref{sec:twistirr}, \(\mathcal P_k\) is reduced,
distinct from \(L-1\), and defines a horizontal curve in the
generic chart. The realizations are dense in that curve.
The containing equation for the filling and the abelian line
therefore give
\[
 A_{Z_{-k}}\doteq (L-1)\mathcal P_k
 \qquad(k\in\Z).
\]

For a negative-clasp satellite with disk twist \(\ell\), a
representation \(\rho_P\) factoring through \(Z_\ell\)
kills \(\mu_A\lambda_A^{-\ell}\). Let
\(\operatorname{ab}\colon\pi_1(X_C)\to\Z\) be companion
abelianization, with \(\operatorname{ab}(\mu_C)=1\) and
\(\operatorname{ab}(\lambda_C)=0\). Define
\[
 \rho_C(\gamma)=\rho_P(\lambda_A)^{\operatorname{ab}(\gamma)} .
\]
These matrices agree on the glued pair in~\eqref{eq:twistgluing}.
The resulting satellite representation has the same outer
eigenvalue pair, so \((L-1)\mathcal P_{-\ell}\mid A_S\).
This extension concerns the companion-abelian filling.
Equation~\eqref{eq:classicalcontainment} is the general bound
for the full satellite curve.

\section{The signed minimum and the satellite mirror}\label{sec:conclusion}

\subsection{Return to the signed sequence}\label{sec:signedreturn}

We obtain a signed all-integer annihilator and pass the factor
survival from the reduced tail to the primitive signed operators.
The exact minimum order in the trefoil calculation uses the
stronger tail-to-all-integer result of Appendix~\ref{app:domain}.

By the \(q\)-holonomicity theorem for the unreduced colored Jones
sequence~\cite[Thm.~1]{S02}, the actual \(\Phi_S\) has a
nonzero positive-tail annihilator. Here its quantum variable
and color multiplier are \(q=t^4\) and \(q^N=M^2\).
The color reflection \(\iota(t,M,L)=(t,M^{-1},L^{-1})\) is an
Ore automorphism. On the odd signed extension it gives
\[
       (\iota(H)\Phi_S)(N)=-(H\Phi_S)(-N).
\]
Take the original row and its reflection, form a nonzero
common left multiple, and clear its coefficient denominators.
Both tails vanish. The finitely many remaining colors \(c\)
are killed by the left scalar \(\prod_c(M-t^{2c})\).
This gives a nonzero all-integer row and hence a nonzero
saturated ideal. Conversely any member of that ideal kills
the positive tail where its clearing coefficients are defined
and invertible. Primitivity of a polynomial representative
removes common \(t+1\) content, so the residue is nonzero.

On the positive tail \(\Phi_S=g_{\mathrm{dim}}F\) with
\(g_{\mathrm{dim}}=(M-M^{-1})/(t^2-t^{-2})\).
For any primitive representative \(H\) of a nonzero saturated
row, the conjugate \(g_{\mathrm{dim}}^{-1}Hg_{\mathrm{dim}}\)
annihilates the reduced tail. By Section~\ref{sec:survival}
it is a left multiple \(U R_a\) of the selected monic
right divisor. The coefficient ratios
\[
 \frac{\sigma^i(g_{\mathrm{dim}})}{g_{\mathrm{dim}}}
   =\frac{t^{-2i}M-t^{2i}M^{-1}}{M-M^{-1}}
\]
are \(t+1\) units of residue one with \(M\) independent,
for every signed shift \(i\).
The \(t+1\) coefficient Gauss valuation is additive on
products because reduction of the automorphism is the identity
and the residual Laurent polynomial ring is a domain.
The conjugate and \(R_a\) have Gauss order zero, hence so
does \(U\), and \(\bar H=\bar U\,\bar R_a\).
As \(L^2E_{\rm back}(L^{-1})=E/x^2\), this proves
\(E\mid f_S\ne0\) for each negative-clasp satellite.
The classical coprimality was established at the same
preferred pair in Section~\ref{sec:torusclass}.

\subsection{Transport under the satellite mirror}\label{sec:wholemirror}

Let \(D=\diag(-1,1)\) be the effect of reflection on
the transported outer and companion pairs, and \(-D\)
its effect on the gluing-adapted axis pair as in
Proposition~\ref{prop:mark}. In column notation
\eqref{eq:twistgluing} has matrix
\[
 G_\ell=\begin{pmatrix}0&1\\1&-\ell\end{pmatrix}.
\]
For every integer \(\ell\) the marking commutes with
reflection by
\[
     (-D)G_\ell=G_{-\ell}D .
\]
Reflection changes the clasp and reverses each signed full
twist. The zero-linking preferred outer longitude is preserved
under the transported orientation. The marked identity is
\begin{equation}\label{eq:mirror-marked}
       (W_{-,\ell}(T_\eta))^*=W_{+,-\ell}(T_{-\eta}) .
\end{equation}

In zero framing \(\Phi_{S^*}(N;t)=\Phi_S(N;t^{-1})\)
on all signed colors. The knot mirror acts by the Ore map
\(\Theta(t,M,L)=(t^{-1},M^{-1},L)\), which preserves the defining
relation and satisfies
\[
 (\Theta(H)\Phi_{S^*})(N;t)=(H\Phi_S)(N;t^{-1}).
\]
It identifies the saturated annihilator ideals. Under
\(t^{-1}+1=(t+1)/t\) the primitive parameter order is preserved.
Thus \(f_{S^*}(M,L)\doteq f_S(M^{-1},L)\).
The same meridian inversion on the full boundary eigenvalue
set gives \(A_{S^*}(M,L)\doteq A_S(M^{-1},L)\).
This transports both divisibility and coprimality, with the
obstruction polynomial represented by \(M^4E(M^{-1},L)\).
Since \(-\ell\) ranges over all integers and \(-\eta\) over
both chiralities, this proves Theorem~\ref{thm:four}.

\subsection{The untwisted positive-companion subfamily}

For the untwisted right-handed companions the quadratic survives
in every primitive representative of a nonzero saturated
annihilator. Here we also give the sharper classical bound and
show that distinct ordered companion parameters give distinct
satellites.

\begin{theorem}[Whitehead doubles of positive torus knots]\label{thm:family}
For coprime integers $2\le r<s$, let $K_{r,s}=Z(T(r,s))$. Its
localized all-integer annihilator is nonzero. Every primitive polynomial
representative $H$ of a nonzero row in this ideal satisfies
\[
 H(-1,M,L)\ne0,\qquad E\mid H(-1,M,L)\quad\text{in }R.
\]
In particular this holds for $H=\alpha_{K_{r,s}}$. For $m=rs$ one has
\[
 A_{K_{r,s}}\mid (L-1)B_m,\qquad \gcd_R(E,A_{K_{r,s}})=1 .
\]
Thus every $K_{r,s}$ violates~\eqref{eq:AJ}. Distinct ordered pairs give
distinct knot types.
\end{theorem}
\begin{proof}
The signed return applies to every primitive row.
At \(\ell=0,\eta=1\), \eqref{eq:classicalcontainment}
has indices \(0,m\), giving the bound, and
Section~\ref{sec:torusclass} proves coprimality and distinctness.
\end{proof}

\begin{corollary}[Odd two-strand companions]\label{cor:F}
For every $k\ge1$ the same $E$ occurs in every nonzero primitive saturated
row after specialization for $Z(T(2,2k+1))$ and is coprime to its full
$A$-polynomial. These satellites are pairwise distinct.
\end{corollary}

\subsection{Classical factors in the specialized minimum}\label{sec:classical-factors}

\begin{corollary}[Untwisted classical divisibility]\label{cor:untwisted-reverse}
For either clasp sign and either companion chirality in
Theorem~\ref{thm:four}, if \(\ell=0\), then \(A_S\mid f_S\) in \(R\).
\end{corollary}
\begin{proof}
For the negative clasp, \(a=\eta rs\ne0\) at \(\ell=0\).
Section~\ref{sec:selected-divisor} selects the full monic form
of \(P_a\). The Gauss specialization of Section~\ref{sec:signedreturn}
then gives \((L-1)E\mathcal P_a\mid f_S\), by the symbol
calculation in Section~\ref{sec:symbol}. The classical bound
\eqref{eq:classicalcontainment} at \(\ell=0\) gives
\(A_S\mid(L-1)\mathcal P_a\) since \(\mathcal P_0=1\).
The mirror identity~\eqref{eq:mirror-marked} and meridian
inversion give the divisibility for the positive clasp.
\end{proof}

For a negative-clasp satellite with \(\ell\ne0\), the same
selection and specialization give \((L-1)\mathcal P_a\mid f_S\)
when \(a\ne0\). At \(a=0\) we only need the factor \(L-1\),
since \(\mathcal P_0=1\). L\^e's general
divisibility~\cite[Prop.~2.3]{Le2006} gives it for the primitive
minimum in the signed unreduced convention of
Section~\ref{sec:comparison}. This is a statement about \(f_S\),
even if the selected divisor specializes to the quadratic alone.

\emph{Twisted reverse inclusion.}
For \(\ell\ne0\), does \(A_S\mid f_S\) hold? For the negative
clasp this is equivalent to asking whether \(\mathcal P_{-\ell}\mid f_S\).
Indeed Section~\ref{sec:pure-filling} proves that the actual
companion-abelian factor \((L-1)\mathcal P_{-\ell}\) divides
\(A_S\). Conversely, each factor occurs once in \(A_S\), and
\eqref{eq:classicalcontainment} puts each among the factors of
\((L-1)\mathcal P_a\) or \(\mathcal P_{-\ell}\).
Those from \((L-1)\mathcal P_a\) already divide \(f_S\).
For the positive clasp, \eqref{eq:mirror-marked} replaces the
negative-clasp parameter by \(-\ell\) and inverts \(M\). The
corresponding factor in the question is therefore
\(\mathcal P_{\ell}(M^{-1},L)\) up to a unit in \(R\).

\subsection{The invariant recurrence ideal}\label{sec:sikora}

The same extra curve separates the radical ideals in Sikora's
\(\mathfrak{sl}_2\) formulation
\cite[Conj.~9]{Sikora2008} and the \(t=-1\) formulation in
\cite[Sec.~1.2, Conj.~2]{TranStrong2013}.
To pass to these ideals, write \(q_{\rm S}\) for Sikora's
quantum parameter and \(E_{\rm sh}\) for his shift.
His \(\rho\)-shifted function in zero framing is the signed
unreduced Jones sequence in dimension color. The quantum dimension
in his equation~(4) gives the same \([N]\) under
\[
 q_{\rm S}^{1/2}=t^2,\qquad Q=M,\qquad E_{\rm sh}=L .
\]
Indeed his multiplier on color \(N\) is \(q_{\rm S}^{N/2}\)
and his shift raises \(N\) by one
\cite[Secs.~1.2--1.4]{Sikora2008}. In the coefficient ring
\(\C[q_{\rm S}^{\pm1/2}]\), evaluation at
\(q_{\rm S}^{1/2}=1\) is therefore obtained by setting \(t=-1\).
Tran uses the signed sequence and our \(t,M,L\) directly.

Under this dictionary a raw Laurent-polynomial annihilator with
complex coefficients is a finite complex linear combination of
rational Laurent-polynomial annihilators. To see this, choose a
basis over \(\Q\) for the span of its finitely many scalar
coefficients after substitution. The sequence is \(\Q(t)\)-valued,
so linear independence over \(\Q(t)\) makes each rational
row annihilate the signed sequence. For each nonzero rational row,
clear numerical denominators and divide by its Laurent-polynomial
content. The resulting primitive row lies in the rational saturation
of Section~\ref{sec:comparison}. Section~\ref{sec:signedreturn}
and its mirror give a nonzero specialization divisible by \(G\),
with \(G=E\) or \(E^\vee\) as in Theorem~\ref{thm:four}.
Multiplying the content back and summing shows that every raw
specialized annihilator vanishes at generic points of \(G=0\).
A content factor may itself specialize to zero.

In particular the Weyl-invariant subideal vanishes there.
Its specializations are functions on the boundary-character quotient
of the eigenvalue torus by simultaneous inversion of \(M,L\).
For either quadratic,
\(M^4L^2G(M^{-1},L^{-1})=G(M,L)\), so the extra curve
descends to this quotient. By Lemma~\ref{lem:aex}
and Theorem~\ref{thm:four}, a generic point of the curve lies
outside the closure of the full classical eigenvalue image.
The quotient is finite with fibers the simultaneous-inversion
orbits, and that closure is invariant. The image of the generic
point is therefore outside the classical boundary-character
image closure, including the image closure from the trivial-character
component in Conjecture~9 of~\cite{Sikora2008}. It lies in the
specialized invariant recurrence zero set. This separates the zero
sets and hence the radicals in the two stated formulations.

\section{Quantum fingerprints of satellite patterns}\label{sec:fingerprints}

For a knot \(K\), keep the actual primitive minimal recurrence
\(\alpha_K\) of the signed sequence in the rationally saturated ideal
of Section~\ref{sec:comparison}. For \(f_K=\alpha_K(-1,M,L)\ne0\), put
\begin{equation}\label{eq:extra-support}
 \Delta_K=\frac{\sfL(f_K)}{\gcd_R(\sfL(f_K),A_K)}
 \qquad\text{in } R=\Q(M)[L^{\pm1}].
\end{equation}
Here \(\sfL\) denotes the square-free part in \(R\), the product of the
distinct irreducible factors up to association. The full horizontal
classical polynomial \(A_K\) from Section~\ref{sec:comparison} includes
\(L-1\).
Thus \(\Delta_K\) records the irreducible factors that occur on the
quantum side and are absent from the classical side. If
\(\Delta_K\doteq1\), then \(\sfL(f_K)\mid A_K\). Equality
in~\eqref{eq:AJ} also requires \(A_K\mid\sfL(f_K)\).

\begin{definition}[Torus-companion fingerprint]\label{def:fingerprint}
Let \(P\) be an oriented, peripherally marked pattern in the solid
torus \(V\), with integer disk twists \(P_\ell\). Let \(\mathcal T\)
be the family of all nontrivial torus companions, of both chiralities,
with the preferred zero framing. Use the marked companion embedding
and the preferred boundary coordinates on the resulting knot.
On the domain where the primitive minimum exists and
\(f_{P_\ell(C)}\ne0\) for every \(\ell\in\Z\) and \(C\in\mathcal T\),
define
\begin{equation}\label{eq:pattern-fingerprint}
 F_{\mathrm{tor}}(P)=\gcd_{\ell\in\Z,\ C\in\mathcal T}
                       \Delta_{P_\ell(C)} .
\end{equation}
The gcd retains the quantum-extra factors common to every test in
this specified family. It is taken up to
association in \(R\), whose units are \(c(M)L^j\) with
\(c(M)\in\Q(M)^\times\) and \(j\in\Z\).
\end{definition}

Every common divisor in~\eqref{eq:pattern-fingerprint} divides
\(\Delta_{P_\ell(C)}\) for any one test, in which each factor occurs once.
Of its finitely many irreducible factors, take the product of those
dividing every test.
This is a greatest common divisor, with the empty product equal to \(1\).
An equivalence preserving the orientations and peripheral markings
extends over each marked companion, so it identifies the tested
satellites and their boundary variables. Replacing \(P\) by \(P_j\)
reindexes the twists by \(\ell\mapsto j+\ell\). Both operations
preserve the gcd up to association.

In the transported coordinates of Section~\ref{sec:wholemirror},
reflection of a knot inverts \(M\) on both sides of the AJ comparison.
The involution \(\rho\) of \(R\) given by \(\rho g=g(M^{-1},L)\)
therefore gives \(\Delta_{K^*}\doteq\rho(\Delta_K)\).
Reflecting the marked satellite gluing reflects the companion and
reverses the signed disk twist. Since the testing family is closed under
reflection and all twists are used, \(F_{\mathrm{tor}}(P^*)\doteq
\rho(F_{\mathrm{tor}}(P))\).

\begin{corollary}[Whitehead patterns]\label{cor:whitehead-fingerprint}
For the two marked Whitehead patterns of this paper,
\[
 F_{\mathrm{tor}}(W_-)\doteq E,\qquad
 F_{\mathrm{tor}}(W_+)\doteq E^\vee .
\]
\end{corollary}
\begin{proof}
Theorem~\ref{thm:four} supplies \(E\) as a divisor of every
\(\Delta_{W_{-,\ell}(C)}\) in the testing family.
For the untwisted right-handed trefoil example,
Theorem~\ref{thm:one} gives \(\sfL(f_S)\doteq (L-1)B_6E\)
and \(A_S\doteq (L-1)B_6\), so \(\Delta_S\doteq E\).
This single test bounds the common divisor from above. The
mirror identity~\eqref{eq:mirror-marked} gives the positive-clasp
formula since \(\rho(E)\doteq E^\vee\).
\end{proof}

The corollary determines the common factors. Individual satellites may
have further factors in \(\Delta\).

\begin{conjecture}[Hyperbolic marked patterns]\label{conj:hyperbolic-fingerprint}
Let \(P\) be an oriented, peripherally marked pattern. If the interior
of \(V\setminus\operatorname{int}N(P)\) admits a complete finite-volume
hyperbolic metric, then the primitive specializations for all \(P_\ell(C)\),
\(\ell\in\Z\), \(C\in\mathcal T\), are nonzero and
\(F_{\mathrm{tor}}(P)\) is a nonunit in \(R\).
\end{conjecture}

The conjecture asks for a common extra factor across the entire
torus-companion family. It places no restriction on winding number.
Whitehead gives the case above. The Mazur pattern remains an
unresolved further test.

\paragraph{Cable comparisons.}
The AJ conjecture is known for substantial ranges of cable knots.
Write \(C_{u,v}(K)\) for the cable in the class \(u\mu+v\lambda\)
on the boundary of the zero-framed companion neighborhood, with
\(\gcd(u,v)=1\) and winding number \(v\ge2\).
For the torus companion \(T(p,q)\) with \(\gcd(p,q)=1\) and
\(|p|>q\ge2\),
Ruppe and Zhang prove AJ when \(u\) lies outside the interval
with endpoints \(0\) and \(pqv\)~\cite[Theorem~1.1]{RuppeZhang2015}.
Ruppe proves it for cables over the figure-eight knot with
\(|u|>4v\)~\cite[Theorem~1.1]{Ruppe2015}. Tran gives criteria for
\((u,2)\)-cables that apply to specified classes of two-bridge and
pretzel knots, with restrictions on both the companion and the odd
cable parameter~\cite[Theorems~1.1, 1.2 and~1.4]{TranCable2017}.
For example, \(C_{13,2}(T(3,2))\) satisfies the Ruppe--Zhang theorem
because \(13>3\cdot2\cdot2=12\). Since \(T(3,2)=T(2,3)\), this
winding-two cable and the winding-zero Whitehead example have the same
right-handed trefoil companion and different AJ outcomes. The cable
parameters specify a slope. The Whitehead parameter \(\ell\) counts
signed disk twists of the marked pattern in Section~\ref{sec:mark}.

The cited cable theorems use the zero-framed signed unreduced sequence
and its primitive minimal recurrence after rational localization.
They specialize at \(t=-1\) with \(M\) independent. Both Ruppe and
Ruppe--Zhang include \(L-1\) in \(A_K\). In Tran's general cable
paper~\cite[Sec.~1.4]{TranCable2017}, the classical side is written
\((L-1)A_K\), which is the full polynomial used here.

For comparison, \(C_{9,2}(4_1)\) satisfies AJ by Ruppe's
Theorem~1.1 because \(9>4\cdot2\)~\cite{Ruppe2015}.
Tran independently proves AJ for \((u,2)\)-cables of \(4_1\)
for odd \(|u|\ge9\)~\cite[Theorem~1]{TranFigureEight2014}.
His 2014 comparison uses the same signed unreduced normalization and
includes \(L-1\) in the classical polynomial. The two AJ equalities
are up to meridian-only factors, so
\(\Delta_{C_{9,2}(4_1)}\doteq1\). Its exterior contains the hyperbolic
figure-eight companion exterior and the Seifert-fibered cable space.
The classical cabling formula already records the companion boundary
polynomial through a resultant~\cite[(1.1)]{Ruppe2015}. A hyperbolic
JSJ piece can thus contribute to the classical curve without forcing a
quantum-extra factor. The cable example has a Seifert pattern exterior.
It therefore falls outside the hypothesis of
Conjecture~\ref{conj:hyperbolic-fingerprint}. This distinction gives
no proof of the conjecture.

\paragraph{Peripheral coordinates under disk twists.}
For a general marked pattern let \(w\) be its winding number. Suppose
the completed-axis boundary curve has a defining polynomial \(E_P\) in
\(R\) of positive longitude span. Write \(\mu_P,\lambda_P\) and
\(\mu_A,\lambda_A\) for the oriented preferred pairs of the pattern
component and the axis before twisting. Choose the sign of \(j\) so
that the new axis meridian pulls back to \(\mu_A\lambda_A^j\).
The longitude \(\lambda_A\) stays fixed. In the drilled exterior,
linking gives
\[
 [\lambda_P]=w[\mu_A],\qquad [\lambda_A]=w[\mu_P].
\]
Filling along the new axis meridian kills
\([\mu_A]+jw[\mu_P]\).
Thus the new preferred outer longitude pulls back to
\(\lambda_P\mu_P^{jw^2}\), the longitude with zero linking in the
new filling. For Whitehead the gluing exponent in
\eqref{eq:twistgluing} is \(j=-\ell\), and the longitude correction
vanishes since \(w=0\).

Axis completeness is preserved by the change of peripheral basis.
For even \(j\) the axis meridian eigenvalue in a chosen
\(\SL(\C)\) lift also keeps its sign. For the corresponding
completed-axis component in the new preferred coordinates one obtains
\[
 E_{P_j}(M,L)\doteq E_P(M,M^{-jw^2}L).
\]
Definition~\ref{def:fingerprint} gives
\(F_{\mathrm{tor}}(P_j)\doteq F_{\mathrm{tor}}(P)\).
Normalize \(E_P\) to have lowest longitude exponent zero and
leading coefficient one. If its degree is \(d_E>0\) and its
constant coefficient is \(c(M)\ne0\), then the normalized shear
has constant coefficient \(M^{jw^2d_E}c(M)\).
For \(w\ne0\) these give infinitely many distinct monic polynomials
even when \(j\) is restricted to even integers. A nonzero polynomial
in \(R\) has only finitely many divisors up to association.
Consequently, on the nonzero-specialization domain, a universal
divisibility by the completed-axis polynomial in preferred coordinates
for every twisted marking would be incompatible with the all-twist
gcd at nonzero winding. The argument leaves
Conjecture~\ref{conj:hyperbolic-fingerprint} unchanged. For a
winding-zero pattern, when the marked completed-axis data define
a horizontal factor over the comparison field, does that factor
divide the stated torus-companion fingerprint?

\paragraph{Skein torsion.}
Does the extra factor reflect torsion in the Kauffman bracket skein
module of the satellite exterior? Belletti and Detcherry's criterion
for a separating incompressible torus also assumes an irreducible
representation nonabelian on each side and noncentral on the
torus~\cite[Thm.~4.1]{BellettiDetcherry2024}. The present comparison
establishes no relation between such torsion and factor survival.

Can one construct a fingerprint intrinsically from marked relative
quantum boundary data and give a gluing law that accounts for its
survival in the actual minimal recurrence? Persistence for arbitrary
companions, invariance without peripheral markings, detection of all
extra factors, a product rule over JSJ pieces, and universal
Mahler measure and volume identities are separate questions.

% !TEX root = ../main.tex
\clearpage
\section*{Acknowledgments}

We thank Yidi Qi for running Comparator on FAS Research Computing, which provided an essential independent check of the Lean formalization. Initial exploration, literature synthesis, proof, computational experiments, Lean formalization, and drafting of this manuscript were carried out with the assistance of Anthropic's Claude Opus 5.5, Google's Gemini 4 Hydrogen, OpenAI's GPT-6 Astra, and SpaceX's Grok Bot, under the direction of the authors, with GPT-6 Astra contributing the principal derivations. The skills and workflows used with GPT, developed by Y.L. and H.-Y.S., will be made publicly available. The authors take full responsibility for all mathematical claims and the exposition.
\appendix
\begin{samepage}
\section{The exact untwisted trefoil calculation}\label{sec:whitehead}

\begin{theorem}[The trefoil satellite]\label{thm:one}
For $S=Z(T(2,3))$, the signed unreduced recurrence and the full classical
boundary polynomial satisfy
\[
 \ord_L\alpha_S=16,\qquad f_S\doteq(L-1)^2 B_6 E,\qquad
 A_S\doteq(L-1) B_6 .
\]
The three distinct irreducible factors $L-1,B_6,E$ have respective
longitude degrees $1,12,2$. Thus $A_S\mid f_S$ and
$f_S\nmid A_S^b$ for every positive integer $b$.
\end{theorem}
\end{samepage}

The finite satellite transform gives an order-fifteen
inhomogeneous row for the reduced sequence. Its specialized
factors are $L-1$, $B_6$ and $E$, and the crossing criterion
proves its irreducibility. The actual sequence is nonrational
by its growth. Nonzero scalar forcing then provides a
one-dimensional extension, giving actual homogeneous order
sixteen. Appendix~\ref{app:domain} passes that order and
specialization to the signed unreduced sequence. On the
classical side we restrict to the filled pattern and extend
representations back across the trefoil torus to obtain the
exact boundary polynomial.

In this section $B=B_6$, $(k,z,y)=(x+x^{-1},Z,J)$ in the classical change of variables, and $\Phi_a=T_a$ in the notation of \eqref{eq:B}.
\subsection{The exact satellite transform}\label{sec:wtransform}
\begin{figure}[ht]
\centering\includegraphics[width=.11\textwidth]{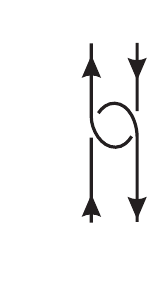}
\caption{The local clasp on printed page 4 of Zheng~\cite{Zheng}.
The untwisted choice and the positive $(3,2)$ companion
input specify the knot in Theorem~\ref{thm:one}.}\label{fig:clasp}
\end{figure}

Put $q=t^4$, $x=M^2$ (so $x=q^N$ at color $N$), and
$(a;q)_j=\prod_{r=0}^{j-1}(1-aq^r)$.
The input $\mathcal T_n(q)=\Phi_{T(2,3)}(2n+1;t)$ is the
unreduced trefoil invariant of dimension $2n+1$.
Write $[N]F_N=\sum_{n=0}^{N-1}\xi_{N,n}\mathcal T_n$ for the
unreduced channel coefficients. Thus the reduced output divides
this entire channel sum by $[N]$.
For the clasp in Figure~\ref{fig:clasp}, Zheng's tangle evaluation
\cite[(2.6)--(2.9), (4.1)]{Zheng} identifies the channel sum. Applying the two-color evaluation in the Hopf basis below gives
\begin{equation}\label{eq:Bailey}
 \begin{split}
 F_N&=\sum_{j=0}^{N-1}\mathcal C_j(N)D_j,\qquad
 \mathcal C_j(N)=(qx;q)_j(q/x;q)_j,\\
 D_j&=q^j(1-q)(q;q)_j
 \sum_{n=0}^j
 \frac{(-1)^nq^{n(n+1)/2}\mathcal T_n}
 {(q;q)_{j-n}(q;q)_{j+n+1}} .
 \end{split}
\end{equation}

For the Hopf inputs $v_s(n)=[s(2n+1)]$, the diagonal two-color
Whitehead evaluation \cite[p.~225]{BBL} is
$[s]\sum_{j<N}b_j\mathcal C_j(N)\mathcal C_j(s)$, where the chosen clasp has
$b_j=q^j(1-q)/(q^{j+1};q)_{j+1}$. Indeed $\mathcal C_j(s;q^{-1})=q^{-j(j+1)}\mathcal C_j(s;q)$ and
$(q^{-j-1};q^{-1})_{j+1}=(-1)^{j+1}q^{-(3j+2)(j+1)/2}(q^{j+1};q)_{j+1}$. Transforming both color factors gives the displayed $b_j$; the top channel below fixes the choice for Zheng's tangle. The other clasp is obtained by
inverting $q$ and has a different coefficient.
Write $v=q^{s/2}$. The finite $q$-binomial theorem gives the coefficient
of $v^{2n+1}$ in $(v-v^{-1})\mathcal C_j(s)$ as
\[
 a_{j,n}=(-1)^nq^{n(n+1)/2}
       \qbinom{2j+1}{j-n},\qquad 0\le n\le j.
\]
Indeed
$(v-v^{-1})\mathcal C_j(s)=(-1)^{j+1}v(v^2)^{-j-1}
q^{j(j+1)/2}(v^2q^{-j};q)_{2j+1}$.
The negative coefficients follow by oddness. Hence
$[s]\mathcal C_j(s)=\sum_{n=0}^j a_{j,n}v_s(n)$.
The $N$ inputs $1\le s\le N$ span the first $N$ channels because
$v_s(n)/[2n+1]$ evaluates a monic Chebyshev polynomial of degree $s-1$
at the distinct points $q^{n+1/2}+q^{-n-1/2}$.
Thus the evaluation identity on those inputs proves the identity
for arbitrary $\mathcal T_0,\ldots,\mathcal T_{N-1}$. Simplifying $b_ja_{j,n}$
gives \eqref{eq:Bailey}. The top channel
$\xi_{N,N-1}=q^{(N-1)/2}(q;q)_{N-1}$ fixes the clasp choice.

\subsection{A finite inner equation with all endpoints}
Extend $\mathcal T_{-n-1}=-\mathcal T_n$. The sum defining $D_j$ is half the sum
over $-j-1\le n\le j$, since the signed kernel and $\mathcal T_n$ both
change sign under $n\mapsto-n-1$. Denote that kernel by $K_{j,n}$,
including its scalar prefactor, and put
$V_j(g)=\frac12\sum_nK_{j,n}g(q^n)\mathcal T_n$, a finite
functional of $g$ at fixed $j$.
For $y=q^j$, $Q=q^n$, direct cancellation of adjacent factorials gives
\[
 \sum_nK_{j,n}Q^a=b_j(1-q^{-a})\mathcal C_j(a),\qquad
 \frac{K_{j,n+1}}{K_{j,n}}
       =-\frac{qQ(1-y/Q)}{1-q^2yQ}.
\]
The first sum is zero for $a=0$ and for $j\ge|a|$.
The trefoil input satisfies
\[
 \mathcal T_{n+1}=q^{-12}Q^{-12}\mathcal T_n+h_C(Q),\qquad
 h_C(q^n)=q^{-3(n+1)}[4n+5]-q^{-9(n+1)}[4n+3].
\]
This follows by taking the two new endpoints in Zheng's finite torus
sum. The Fourier exponents of $h_C$ are $-11,-7,-5,-1$.

Here is the reindexing which turns the trefoil relation into the inner
row. For any Laurent test $g$, the adjacent kernel identity gives
\[
 \sum_nK_{j,n}(1-qyQ)g(Q)\mathcal T_n
       =-q\sum_nK_{j,n}(Q-y)g(qQ)\mathcal T_{n+1}.
\]
Indeed before reindexing its cancellation is
$K_{j,n+1}(1-q^2yQ)=-q(Q-y)K_{j,n}$.
The factor on the left becomes $1-qyQ$ when
the summation variable is replaced by $n+1$.
Take $g(Q)=Q^5$ and substitute the trefoil equation.
After moving its homogeneous part to the left the test polynomial is
\[
 R(Q)=Q^5-qyQ^6+q^{-6}Q^{-6}-q^{-6}yQ^{-7}.
\]
The forcing is $-q^6(Q-y)Q^5h_C(Q)$, whose exponents
are between $-6$ and $5$. Its signed sum vanishes for $j\ge6$.
The functional $V_j$ is invariant under
$Q\mapsto q^{-1}Q^{-1}$, and direct grouping gives
\[
 R(Q)+R(q^{-1}Q^{-1})
                  =S_5+(1-qy)S_6-qyS_7,
\]
where the $S_a$ are evaluated at $w=Q+q^{-1}Q^{-1}$ and
$S_0=2$, $S_1=w$, $S_{a+1}=wS_a-q^{-1}S_{a-1}$.
To turn moments into shifts in $j$, define the operator on
sequences $d=(d_j)$ by
\[
 (\mathcal V d)_j=(q^{-1}y^{-1}+y)d_j-(1-y)d_{j-1}/y .
\]
We apply it to $D$, with $D_j=0$ for $j<0$.
For clarity, the backward kernel ratio, before summing, is
\[
 \frac{K_{j-1,n}}{K_{j,n}}=
 \frac{(1-y/Q)(1-qyQ)}{q(1-y)}
   =\frac{1+qy^2-qy(Q+q^{-1}Q^{-1})}{q(1-y)}.
\]
Its numerator kills both lost endpoints $n=j$ and $n=-j-1$.
Solving this identity for the moment of
$Q+q^{-1}Q^{-1}$ gives the displayed $\mathcal V$.
At $j=0$ its
backward term vanishes, and both remaining channels have
$Q+q^{-1}Q^{-1}=1+q^{-1}$. Consequently
\begin{equation}\label{eq:defect}
 e_j=\big[S_5(\mathcal V)+(1-qy)S_6(\mathcal V)
                     -qyS_7(\mathcal V)\big]D_j
       =0\quad(j\ge6).
\end{equation}
In reindexing the adjacent sum, the missing right endpoint has $Q-y=0$
and the missing left endpoint has $1-qyQ=0$. These are
cross-multiplied finite identities, including the moving endpoints.

\subsection{The pole window and the order-fifteen row}
The inner row has degree seven. The adjoint of $\mathcal V$ against
$\mathcal C_j(N)$ is
\[
 Wg=(q^{-1}y^{-1}+y)g(y)
 -\frac{(1-qy)(1-qxy)(1-qy/x)}{qy}g(qy).
\]
The endpoint identities $\mathcal C_N(N)=0$ at the top and $1-y=0$ at the
bottom eliminate boundary terms in summation by parts. If $g$ is regular at $y=0$, the
possible $y^{-1}$ principal parts cancel. At infinity $W$ raises
the degree by at most two. The simple pole window
$y=q^r/x$, $0\le r\le14$, is stable because substitution $y\mapsto qy$
can add only the pole $r=-1$, canceled by $1-qxy$.

Choose
\[
 \omega(y)=\prod_{r=0}^{14}(1-xy/q^r)^{-1},\qquad
 H=S_5(W)\omega+S_6(W)((1-qy)\omega)-S_7(W)(qy\omega).
\]
The bound $-15+1+2\cdot7=0$ makes $H$ bounded at infinity.
Thus after clearing this fifteen-pole denominator its numerator has
degree at most fifteen. This determines the space spanned by the
external shifts below.

For $0\le i\le15$ put
\begin{equation}\label{eq:Rbasis}
 R_0=1,\qquad
 R_i=\prod_{r=0}^{i-1}
 \frac{(1-x/q^r)(1-q^{r+1}y/x)}
      {(1-q^{r+1}/x)(1-xy/q^r)} .
\end{equation}
Their cleared numerators belong to the same degree-fifteen space.
At $q=1$, $R_i=((x-y)/(xy-1))^i$. Multiplying by the common
denominator gives the degree-fifteen homogeneous monomials in two
independent linear forms, with Möbius determinant $1-x^2$.
They are therefore a basis whenever $x\ne0,\pm1$, with determinant
a unit on the residual disks used below. This proves both existence and
regularity of the unique coefficients in
$H=\sum_{i=0}^{15}a_i(q,x)R_i$.

For $N>15$, the factorial identity
$\mathcal C_j(N-i)=\mathcal C_j(N)R_i(q^j)$ includes the vanishing lost top terms.
Here $U=L^{-1}$ shifts the external color $N$, and its
commuting symbol is $u$. The sixteen basis terms use shifts
$0$ through $15$. With \eqref{eq:defect} they give the exact
finite specification
\begin{equation}\label{eq:PF}
 PF=b,\qquad P=\sum_{i=0}^{15}a_iU^i,\qquad
 b=\sum_{j=0}^{5}\mathcal C_j(N)\omega(q^j)e_j,\quad U=L^{-1}.
\end{equation}
On the right $\mathcal C_j(N)$ is the displayed polynomial in independent $x$.
There is also a direct pole test: at $x=q^{14}$ the $j=0$ denominator of $\omega(q^j)$ has a simple zero, $\mathcal C_0=1$ and $e_0=2(1-q)$ by the endpoint row. For $j=1,\ldots,5$ no denominator factor indexed $0,\ldots,14$ can vanish there. The residue is therefore nonzero. The value $b=31260$ at $q=4$, $q^{1/2}=2$, $N=27$ modulo
$1000003$, with all denominators units, certifies $b\ne0$.
Appendix~\ref{app:finite} specifies the finite evaluation and exact
arithmetic required to check it.

\subsection{The commuting symbol and its factors}
For the fixed trefoil row write
$P=P^{[0]}+hP^{[1]}+O(h^2)$ with $h=q^{-1}-1$,
holding $(x,u)$ fixed in the coefficient expansion.
The bracketed superscripts record deformation degree.
The variable $y=q^j$ for finite evaluation is independent of
$x$ in the auxiliary symbol chart.
At $q=1$, $W$ is multiplication by $z=1+k-ky+y^2$ and $R_i=u^i$
with $u=(x-y)/(xy-1)$. The Chebyshev recursion gives
\begin{equation}\label{eq:Pzero}
 P^{[0]}=(1-xy)^{-15}(1-y)(z+2)\Phi_6(z,e),\quad
 y=\frac{x+u}{1+xu},\quad u=L^{-1}.
\end{equation}
Indeed, let $p_0=1$, $p_1=z$ and $p_{a+1}=zp_a-p_{a-1}$.
At $q=1$, evaluate $S_a$ at $z$. The identities
$S_7-S_5=(z^2-4)p_5$ and $S_6+S_7=(z+2)(p_6-p_5)$
follow from the Chebyshev recursions. Using $z-2=(y-1)e$ gives
$S_5+(1-y)S_6-yS_7=(1-y)(z+2)(p_6+(e-1)p_5)$.
The last factor is $\Phi_6$ by \eqref{eq:B}.
In forward boundary coordinates the exact
normalizations are
\begin{equation}\label{eq:sign}
 \begin{split}
 1-y&=\frac{(1-x)(L-1)}{L+x},\\
 z+2&=\frac{(x+1)^2 E}{x(L+x)^2},\\
 1-xy&=\frac{L(1-x^2)}{L+x},\qquad
 \Phi_6=\frac{B}{x^6(L+x)^{12}}.
 \end{split}
\end{equation}
All fifteen powers of $L+x$ cancel, so
$P^{[0]}\doteq u^{15}(L-1)BE$.

The discriminant of $E$ is $(M^2-1)^2(M^4-6M^2+1)$,
whose final factor is not a square in $\Q(M)$.
For $B$, the degree-preserving specialization $M=2$ modulo $7$ is
\begin{equation}\label{eq:Bmod}
 L^{12}+3L^8+2L^7+3L^6-L^5+3L^3+L^2-L+1.
\end{equation}
Its gcds with $L^{7^4}-L$ and $L^{7^6}-L$ are one, and
$L^{7^{12}}\equiv L$ modulo it. The finite-field irreducibility
criterion proves that the reduction is irreducible.
To check degree preservation, put $c=x(L+x)^2$ and $Z=cz$.
The polynomials $B_a=c^a\Phi_a$ satisfy $B_0=1$,
$B_1=c(z-1+e)$ and $B_{a+1}=ZB_a-c^2B_{a-1}$.
Here $Z,B_1\in\Z[x,L]$, and this recursion gives $B_6=B$
with degree twelve, leading coefficient $x^{12}$ and constant
coefficient $x^{12}$. Thus $B/x^{12}$ is monic in $L$ over
$\Z[M,M^{-1}]$. A monic factorization over $\Q(M)$ would, by
Gauss's lemma, specialize at $M=2$ modulo $7$ with both
degrees intact. This proves irreducibility of $B$. The factors
$L-1,B,E$ are distinct by their degrees.

\subsection{Irreducibility of the order-fifteen operator}\label{sec:crossings}
The two edges needed in Theorem~\ref{thm:graph} join
$E$ to $B$ and $B$ to $L-1$. The following calculation fixes the
quantum sign before computing an index. Differentiate $H$ at fixed
$x,y$ with respect to $q$, and then use $h=q^{-1}-1$.
At a crossing, $P^{[0]}=(P^{[0]})_u=0$. Put
\[
 \rho=\frac{x}{1-x}-\frac{xy}{1-xy},\qquad
 \eta=\frac1{x-1}-\frac{y}{x-y}.
\]
Logarithmic differentiation of \eqref{eq:Rbasis} gives
$\partial_q\log R_i=\rho\,i(i-1)/2+\eta\,i(i+1)/2$.
For local factors $G,D$ and $T=P^{[0]}/(GD)$ it follows that
\begin{equation}\label{eq:physicaljet}
 P^{[1]}/T=-H_q/T+(\rho+\eta)u^2G_uD_u.
\end{equation}
The sign in $-H_q$ is forced by $dq/dh=-1$.

Here is an explicit way to compute $H_q$ without expanding the
sixteen unknown coefficients. To first order
$W=z+(q-1)(d+v\partial_y)$, where
\[
 d=y-k-(1-y)(2y-k),\qquad
 v=-(1-y)(1-xy)(1-y/x).
\]
For a seed $g(q,y)$ and Chebyshev term $S_a(W)g$, its derivative
at $q=1$ is
\begin{equation}\label{eq:jetalgorithm}
 S_ag_q+\frac{zS_a'-aS_a}{2}g+
 S_a'(dg+vg_y)+\frac{vz_y S_a''}{2}g .
\end{equation}
The first two terms differentiate the seed and the Chebyshev
parameter. The last two follow by differentiating each occurrence
of $W$ in its powers. Formula \eqref{eq:jetalgorithm} and the
three seeds $\omega,(1-qy)\omega,-qy\omega$ determine the entire first jet.

For the $E/B$ crossing take $G=z+2$, $D=\Phi_6$.
The elementary evaluations
$\Phi_6(-2,e)=13-6e$ and $z-2=(y-1)e$ give
\[
 z=-2,\quad e=13/6,\quad y=-11/13,\quad
 78x^2+157x+78=0 .
\]
Work in this quadratic field and use $u=(x-y)/(xy-1)$.
The recurrence for $\Phi_6$ gives $\Phi_{6,z}=-91/6$,
$\Phi_{6,e}=-6$, while $z_y=25/78$.
For the Chebyshev numerator, write $\mathcal N=S_5+(1-y)S_6-yS_7$, treating $z,y$ as independent. Its derivatives
at this point are $\mathcal N_z=\mathcal N_y=0$, $\mathcal N_{zy}=-13$, $\mathcal N_{zz}=-44$.
Consequently \eqref{eq:jetalgorithm} reduces to
\[
 H_q/\omega=v(-13-22z_y)-1-y=\frac{2114}{6591}.
\]
Substitution into \eqref{eq:physicaljet}, with $T=(1-y)\omega$,
gives
\[
 P^{[1]}/T=\frac{10643}{146484},\qquad
 G_uD_x-D_uG_x=-\frac{25x}{6}-\frac{169}{36}.
\]
Both $u$ derivatives are nonzero. Reducing
\eqref{eq:crossindex} modulo $78x^2+157x+78$ gives
$I_G=-1/2$, and $I_G+I_D=-1$ gives $I_D=-1/2$.

For the second edge set $G=1-y$, $D=\Phi_6$.
At $y=u=1$, $z=2$ and $\Phi_6(2,e)=1+6e$.
Thus $k=13/6$ and $x=2/3$ or $3/2$.
At $x=2/3$ the first-order multiplier contributes $28$ to $H_q/\omega$,
the Chebyshev parameter contributes $-22$, and the seed contributes
$-4$, giving $H_q/\omega=2$. Here $\rho+\eta=0$,
$T=(z+2)\omega=4\omega$, $G_x=0$, $G_u=-1/5$, $D_u=25/36$,
and the Jacobian is $-3/2$. Hence $P^{[1]}/T=-1/2$ and
\eqref{eq:crossindex} again gives $I_G=I_D=-1/2$.
At $x=3/2$ the Jacobian is $-2/3$ and the same formula gives
the same indices.

We check the disk hypotheses explicitly. At $u=0$ and $\infty$
the values of $y$ are $x$ and $1/x$, with $z=k$ and
$\Phi_6=x^6,x^{-6}$, respectively. Since
$1-xy=(1-x^2)/(1+xu)$ and the numerator of \eqref{eq:Pzero}
is nonzero at $y=1/x$, $P^{[0]}$ has degree fifteen in $u$ with
nonzero leading coefficient. Its constant coefficient is also nonzero. The values $x=0,\pm1$ and $k=-2$ are absent from
the chosen fibers. In detail, if $\Delta=\prod_{a=0}^{14}(1-xy/q^a)$, applying the adjoint to $b(y)/\Delta$ and clearing $\Delta$ gives
\[
 yb(y)+\frac{b(y)-(1-qy)(1-qy/x)(1-xy/q^{14})b(qy)}{qy}.
\]
The bracket is divisible by $y$ over the joint local ring. Thus the cleared right side of the coefficient solve is integral even when the fifteen poles coalesce. The cleared Möbius basis is a unit basis
there, so the full quantum coefficients are regular.
At $y=\infty$, equivalently $u=-1/x$, the numerator and denominator
in \eqref{eq:Pzero} both have degree fifteen in $y$ and nonzero leading
coefficients. Their ratio is nonzero there.

Finally $E$ never meets $L-1$, since $z+2=4$ at $y=1$.
An $E/B$ crossing forces $k=-157/78$ and a unique $y=-11/13$.
At that meridian $\Phi_6(2,2-k)=326/13\ne0$.
An abelian/$B$ crossing forces $k=13/6$, which excludes the
other edge. Each selected meridian fiber therefore has exactly
one inter-factor crossing. The graph is the connected path
$L-1$---$B$---$E$. Theorem~\ref{thm:graph} proves that $P$ is
irreducible over $\Q(t,M)$.

\subsection{The scalar extension and order sixteen}\label{sec:minimum}
The inhomogeneous equation gives a homogeneous equation for the
class of $F$ modulo rational scalar sequences. We compute the order
of the recurrence for $F$ itself from the resulting extension.

\begin{lemma}[The scalar-line extension]\label{lem:scalar}
Let $PF=b\ne0$ on a positive tail, with $P$ irreducible of
order $d$ over the rational sequence coefficient field.
If $F$ is not a rational scalar sequence, its actual homogeneous
minimum has order $d+1$ and is a left associate of
$(U-\sigma(b)/b)P$.
\end{lemma}
\begin{proof}
The cyclic span $V=\mathcal D F$ contains the scalar line
$\mathcal S=\Q(t,M)\cdot 1$ because $b^{-1}PF=1$.
The class of $F$ generates $V/\mathcal S$ and is killed by $P$.
The map $\mathcal D/\mathcal D P\to V/\mathcal S$ sending the class
of $1$ to that of $F$ is onto. Its source is simple by irreducibility
of $P$, and its target is nonzero because $F$ is nonscalar. Thus
it is an isomorphism. The exact sequence
\[
 0\longrightarrow \mathcal S\longrightarrow V
 \longrightarrow \mathcal D/\mathcal D P\longrightarrow 0
\]
gives $\dim_{\Q(t,M)} V=1+d$.
The displayed homogeneous operator annihilates $F$ and has this
order, so it is minimal.
\end{proof}

To show that the class of $F$ is nonzero, fix a generic
real $0<q<1$. Zheng's torus input is
\[
 \mathcal T_n=q^{-6n(n+1)}A_n/(1-q),\qquad
 A_n=\sum_{\ell=-n}^n(q^{6\ell^2+\ell}-q^{6\ell^2+5\ell+1}).
\]
Euler's pentagonal identity gives $A_\infty=(q;q)_\infty>0$.
The increments are
$2q^{6n^2}(\cosh(n\log q)-q\cosh(5n\log q))$.
The ratio $\cosh(5s)/\cosh s$ increases for $s>0$, so these
increments change sign at most once, from positive to negative.
Since $A_0=1-q\ge(q;q)_\infty$, all $A_n\ge(q;q)_\infty$.
Positivity of the clasp coefficients follows from
\[
 \xi_{N,n}=q^{N^2-(N+1)/2}x^{-n}(q;q)_n
 \sum_{a+b=N-1-n}\qbinom{n+a}{n}\qbinom{n+b}{n}x^{-b}.
\]
Taking only its top channel yields
\[
 F_N\ge
 \frac{q^{(N-1)/2-6N(N-1)}(q;q)_\infty^2}{(1-q)[N]} .
\]
The logarithm of this lower bound grows quadratically in $N$.
A rational function of $t,M$, evaluated at a fixed generic $t$
and $M=q^{N/2}$, grows at most exponentially in $N$ because
its Laurent order at $M=0$ is finite and its remaining analytic
factor is bounded and nonzero there. Choosing $q$ away from the
finitely many identically exceptional coefficient values proves
that $F$ is not rational. Lemma~\ref{lem:scalar} now gives order
sixteen.

The dimension gauge $[N]$ has shift ratios of residue one.
Likewise $\sigma b/b$ has residue one after removing any scalar
$h$-valuation of $b$. Conjugate the homogeneous row by $[N]$,
left-shift by sixteen to use nonnegative powers of $L$, clear
denominators on the left, and remove common content.
Appendix~\ref{app:domain} extends its tail annihilation to all integers.
A smaller all-integer relation would restrict to a smaller tail
relation. The primitive symbol is therefore $(L-1)^2BE$.

\subsection{The classical boundary curve}\label{sec:classical}

\begin{proposition}[Classical boundary curve of the marked satellite]
\label{prop:classical}
For the zero-framed satellite $S$ of Theorem~\ref{thm:one},
the horizontal boundary curve in the preferred coordinates has equation
$A_S\doteq (L-1)B$. In particular, $E\nmid A_S$.
\end{proposition}

\begin{proof}
Restriction gives the necessary filling as in
Section~\ref{sec:torusclass}. For equality we will realize a
dense part of its nonabelian boundary curve and extend those
representations over the trefoil piece.

Cut the exterior along the axis torus of the Whitehead pattern
in Figure~\ref{fig:filling}. Denote the preferred axis pair by
$(\mu_{\rm ax},\lambda_{\rm ax})$ and the trefoil pair by
$(\mu_C,\lambda_C)$. The untwisted satellite identifies
$(\mu_{\rm ax},\lambda_{\rm ax})=(\lambda_C,\mu_C)$.
On the outer boundary, $M,L$ are the eigenvalues of the meridian
and preferred longitude of the Whitehead component.

The trefoil group has
$a_C^2=b_C^3$, $\mu_C=a_Cb_C^{-1}$,
$\lambda_C=a_C^2\mu_C^{-6}$.
The centralizer of a nonscalar matrix in $\SL(\C)$ is abelian.
For a nonabelian trefoil image, the central matrix $a_C^2=b_C^3$
is therefore scalar. If it were $I$, then $a_C^2=I$ would make
$a_C$ scalar and the image abelian. Hence $a_C^2=b_C^3=-I$.
The boundary identification gives $\mu_{\rm ax}\lambda_{\rm ax}^6=-I$
on the pattern representation. Multiply this representation by the
sign character that takes the axis meridian to $-1$ and the outer
meridian to $1$. Both preferred longitudes lie in the kernel of this
character because the Whitehead components have linking number zero.
This changes the relation to $\mu_{\rm ax}\lambda_{\rm ax}^6=I$
and preserves the outer eigenvalues.

To determine the twist sign, orient two test strands at the twist
site in parallel. A positive full twist increases their linking
number by one. Surgery with slope $\mu_{\rm ax}\lambda_{\rm ax}^{s}$
has coefficient $1/s$ and changes their linking number by $-s$.
For $s=6$ this gives the twist knot $Z_{-6}$ in the signed convention
of the twist-polynomial calculation below. If the trefoil image is
abelian, then $\lambda_C=I$. Killing $\mu_{\rm ax}$ gives the outer
unknot and hence $L=1$.

Conversely, let a nonabelian representation of the filled pattern
$Z_{-6}$ have axis longitude $D_{\rm ax}$. If $D_{\rm ax}=\pm I$,
the filling forces $\mu_{\rm ax}=I$, so the representation factors
through the unknot. For nonscalar $D_{\rm ax}$, choose
$b_C\in\SL(\C)$ with $\tr b_C=1$ and $\tr(D_{\rm ax}b_C)=0$.
Set $a_C=D_{\rm ax}b_C$. If $D_{\rm ax}=\diag(d,d^{-1})$ has
distinct eigenvalues, take
$(b_C)_{11}=-d^{-1}/(d-d^{-1})$, $(b_C)_{22}=1-(b_C)_{11}$,
and $(b_C)_{12}=1$. The determinant condition fixes $(b_C)_{21}$.
If $D_{\rm ax}=\pm\left(\begin{smallmatrix}1&1\\0&1\end{smallmatrix}\right)$,
take $(b_C)_{21}=-1$, impose the trace, and solve for $(b_C)_{12}$.
Cayley--Hamilton gives $a_C^2=b_C^3=-I$. Thus the trefoil pair
is $(D_{\rm ax},-D_{\rm ax}^{-6})$. Applying the same sign
character to the filled pattern gives this matching pair in the
exchanged order. The two representations glue and realize the
same outer eigenvalues.

For the twist-polynomial calculation, put $\zeta=1-y$ and
$\delta=(M-M^{-1})^2=k-2$. In the marked coordinates the
trace in Mathews's recursion~\cite[Section 2]{Mathews} is
$2+\zeta^2+\delta\zeta=z$. The recursion is
$r_{m+1}=zr_m-r_{m-1}$ with $r_0=1$, $r_1=k-y=1-e$.
The negative index gives $r_{-1}=z-1+e$ and hence
$r_{-6}=\Phi_6(z,e)$. Clearing denominators gives $B$ in
\eqref{eq:B}. Mathews's erratum changes the closed formula for
$n\le0$ in his Theorem~1.1 and retains the recursion in
Lemma~2.1~\cite{MathewsErratum}.

For the signed calculation, write $X=\rho(x)$ and $Y=\rho(y)$
for the positive-generator meridian matrices in Schilder's
presentation, as in Section~\ref{sec:word}. Set $Q=X^{-1}$ for
the reflected meridian, $P=Y^{-1}QY$ and $W=PQ^{-1}P^{-1}Q$. Substitution gives axis longitude $W$ and outer longitude $P^{-1}QYQ^{-1}Y^{-1}Q$. The sign-corrected filling has $Y=W^6$, so $PW^{-6}=W^{-6}Q$. At a generic irreducible point conjugate $P,Q$ to
\[
 P=\begin{pmatrix}M&1\\0&M^{-1}\end{pmatrix},\qquad
 Q=\begin{pmatrix}M&0\\\zeta&M^{-1}\end{pmatrix}.
\]
Their indicated distinct eigenlines supply the basis, and irreducibility permits the off-diagonal rescaling. The upper-right relation is $r(W^{-6})=0$, where $r(U)=(M-M^{-1})U_{12}+U_{22}$. Multiplication gives $\tr W=z$, $r(I)=1$, $r(W^{-1})=z-1+e$, and $(W^{-1})_{21}=\zeta(W^{-1})_{12}$. By Cayley--Hamilton the last relation also holds for $W^{-6}$, so the upper-right equation implies all four entries of the matrix relation. Conversely on its generic locus these matrices with $X=Q^{-1},Y=W^6$ give the filled presentation, equivalently the negative-six twist presentation of \cite[Sec.~2, (1)--(2)]{Mathews}. Evaluating the outer longitude word on $(M-M^{-1},\zeta)^t$ gives $L=(M-M^{-1}-M\zeta)/(M-M^{-1}+\zeta/M)$, equivalently $\zeta=1-y$. This map is invertible for $M^4\ne1$ on the generic curve. Reducible representations have outer longitude eigenvalue one, whereas the irreducible degree-twelve factor is distinct from that line. Thus the generic nonabelian twist curve supplies a dense set on $B=0$ and the extension above realizes it on the satellite. Taking Zariski closure of the boundary
images gives $A_S\doteq(L-1)B$. The factor $E$ is distinct
from both irreducible classical factors.
\end{proof}

The classical proposition and the homogeneous minimum in
Section~\ref{sec:minimum} prove Theorem~\ref{thm:one}.

\section{The quadratic curve of the Whitehead pattern}\label{sec:geometry}
The curve $E=0$ is excluded as an irreducible component of the satellite boundary curve. Its intersections with other curves remain possible. Here we normalize
the curve, identify it in the marked pattern before gluing, and use
the actual root-of-unity sum to distinguish the two lifts at a node.
In this appendix \(d\) is the odd companion label in \(T(2,d)\),
with \(d\ge3\), for the untwisted negative-clasp satellite with
right-handed companion. Thus the effective transfer index
\(\eta rs-\ell\) of Section~\ref{sec:family} is \(2d\) here.
In the diagonal limit \(N\to\infty\), the companion label stays fixed.

For an integral Laurent polynomial $P$ in two variables, its logarithmic Mahler measure uses probability Haar measure:
\begin{equation}\label{eq:haar}
 m(P)=\frac{1}{(2\pi)^2}\int_0^{2\pi}\!\int_0^{2\pi}
 \log|P(e^{\imunit\theta},e^{\imunit\phi})|\,d\theta\,d\phi .
\end{equation}
We use the primitive monic representative $E=L^2-(M^4-4M^2+1)L+M^4$ of \eqref{eq:E}. A unit monomial of coefficient $\pm1$ has log modulus zero. Multiplication by a constant $\gamma\ne0$ adds $\log|\gamma|$, and multiplication by a meridian rational function adds its one-variable logarithmic integral. Thus the representative matters when passing from the localized recurrence ring to \eqref{eq:haar}.

\begin{samepage}
\begin{theorem}\label{thm:geometry}
For the untwisted negative-clasp satellite of \(T(2,d)\), with
\(d\ge3\) odd and right-handed companion, the following hold.
\begin{enumerate}
\item \textit{Stationary curve.}
The homogeneous stationary channel \(Q=-1\) gives \(E=0\)
on \(x=M^2\ne0,\ x^2\ne1\). Its smooth projective normalization
is \(Y^2=X^3-X\) (Sections~\ref{sec:stationary} and~\ref{sec:reg}).
\item \textit{Marked realization.}
In the pattern exterior, this curve is the completed-axis
slice of Tillmann's nondegenerate shape chart, with the preferred
longitude inverted across the source mirror (Section~\ref{sec:shape}).
\item \textit{Mahler measure.}
For the monic integral representative in \eqref{eq:E},
\(m(E)=4\cat/\pi\), where
\(\cat=\sum_{n\ge0}(-1)^n/(2n+1)^2\)
(Section~\ref{sec:mahler}).
\item \textit{The actual diagonal.}
For each fixed \(d\), the dominant window in the actual reduced
Jones sum at \(q_N=e^{2\pi\imunit/N}\) contributes a nonzero leading term.
The auxiliary coordinates \(z_\ell,z_{\ell^\prime}\) in
\eqref{eq:dualratio} converge to \(-\imunit\) on that window. The
normalized limit in \eqref{eq:branchreturn} selects \(w=\imunit\)
above \((M,L)=(-1,-1)\) in the marked slice
(Section~\ref{sec:physical}).
\end{enumerate}
\end{theorem}
\end{samepage}

\subsection{The generating function and the stationary equations}\label{sec:stationary}
The actual input is the unreduced companion invariant
\(\torcolor_n^{(d)}(q)=\Phi_{T(2,d)}(2n+1;t)\), so it equals
\(v_n\) in \eqref{eq:commoninput} at zero twist.
At generic $q$ the unreduced clasp sum can be written
\[
 q^{N^2-(N+1)/2}[z^{N-1}]\sum_{n\ge0}
 \frac{(z/x)^n(q;q)_n}{(z;q)_{n+1}(z/x;q)_{n+1}} \torcolor_n^{(d)},\qquad x=q^N.
\]
Indeed the coefficient in the $n$th term before $\torcolor_n^{(d)}$ and the common prefactor is
\[
 x^{-n}(q;q)_n\!\sum_{\ell+\ell^\prime=N-1-n} x^{-\ell}
 \qbinom{n+\ell}{n}\qbinom{n+\ell^\prime}{n}.
\]
The product of the two binomials times $(q;q)_n$ is $\prod_{r=1}^n(1-q^{\ell+r})(1-q^{N-\ell-r})/(1-q^r)$, which is Zheng's finite clasp coefficient in \cite[(2.3), (2.7)--(2.9)]{Zheng}. The power of $z$ makes $n<N$ automatically. Section~\ref{sec:inner} gives the inhomogeneous recurrence of the actual input. Its homogeneous solution for $m=2d$ is
\(v_n^{\mathrm{hom}}(q)=q^{-2dn(n+1)}\).

We use this homogeneous solution to identify a stationary channel.
Section~\ref{sec:physical} will return to the actual reduced sequence
for the physical comparison.
Use $q=e^h$ and replace the input by \(v_n^{\mathrm{hom}}\).
Here \(z\) is the coefficient-extraction variable, \(Q=q^n\)
the channel multiplier, and \(x=q^N\) the external color
multiplier. In the stationary equations we use \(y=zQ\). The leading formal action, with the coefficient extraction and exterior gauge included, is
\begin{align*}
 S={}&\log(z/x)\log Q+\operatorname{Li}_2(1)-\operatorname{Li}_2(Q)
       -\operatorname{Li}_2(z)+\operatorname{Li}_2(zQ)\\
 &-\operatorname{Li}_2(z/x)+\operatorname{Li}_2(zQ/x)
       -2d(\log Q)^2-\log x\log z+(\log x)^2 .
\end{align*}
Here the formal contribution of $(c;q)_n$ to $h\log$ is $\operatorname{Li}_2(c)-\operatorname{Li}_2(cQ)$, obtained by integrating the logarithm of each successive factor. With coherent local logarithms, differentiation using $d\operatorname{Li}_2(v)/d\log v=-\log(1-v)$ gives
\begin{align*}
 y(1-Q)&=Q^{4d+1}(1-y)(x-y),\\
 xQ(1-y)(x-y)&=(Q-y)(xQ-y),\\
 L=\exp(\partial_{\log x}S)&=\frac{x^2 Q(x-y)}{y(xQ-y)}.
\end{align*}
Remove zero among $x,Q,z,y$, remove one among $Q,z,y,z/x,y/x$, and require $xQ-y\ne0$. The first two are exponentiated critical equations. One can subtract integer multiples of $2\pi\imunit\log z,2\pi\imunit\log Q$ locally so that the internal derivatives vanish.

On $Q=-1$ the first equation becomes $T=y^2+(1-x)y+x=0$, and the second has residual $-(x+1)T$. The external shift is $L=x^2(x-y)/[y(x+y)]$. For $x\ne0$, $x^2\ne1$, substitution of $y=0,1,x,-x,-1$ in $T$ gives $x,2,2x,2x^2,2x$, respectively. The resultant with $-Ly(x+y)+x^2(x-y)$ is $2x^3E$. Conversely put $d_y=x^2+(2x-1)L$, $y=x(L+x^2)/d_y$. Since $\operatorname{Res}_L(E,d_y)=2x^3(x-1)^2$, the denominator is nonzero. Substitution gives $2x^3E/d_y^2$ in $T$ and $-2Lx^3E/d_y^2$ in the shift residual. This proves both inclusions and identifies the excluded fibers.

\subsection{Normalization, regulator and stationary action}\label{sec:reg}
We pass to the normalization to distinguish the points above the nodes.
The coordinate introduced below will also express the regulator primitive.
Put $p=1/y=(x+L)/[x(x-1)]$ and $f=(1+p)/(1-p)$. Then $x=f/p$ and $L=f^2$. Conversely these two expressions, with a choice of $M^2=x$, recover the curve. Put $X=-p$, $Y=Mp(1-p)$. Then $Y^2=X^3-X$ and $M=Y/[p(1-p)]$. The four simple zeros or poles of $x(p)$ are $0,1,-1,\infty$, so the double cover of the $p$-line is branched precisely there. Its smooth projective model is elliptic with $j=1728$. The two nodes of the plane curve in the algebraic torus,
at $(M,L)=(\pm1,-1)$, each lift to $p=\pm\imunit$, with the indicated choice of $M$. Indeed the tangent cones in displacements $s,t$ of $M,L$ are $t^2\pm4st+8s^2$; away from these points the inverse coordinates or the simple discriminant roots give smooth points.

To express the stationary action on this curve, first compute a primitive
of its regulator form.
Let $\eta(a,b)=\log|a|\,d\arg b-\log|b|\,d\arg a$
be the real one-form, and let $D$ denote the Bloch--Wigner
dilogarithm, with $dD(z)=\eta(z,1-z)$. By bilinearity and $-\!L/x=-pf$,
\[
 2\eta(M,L)=\eta(f/p,-pf)=-2\eta(p,f)
 =2\,d\bigl(D(p)-D(-p)\bigr).
\]
The primitive $U(p)=D(p)-D(-p)$ is single valued. Composing the continuous extension of $D$ on the Riemann sphere with the morphism to the $p$-line extends $U$ to the projective normalization. It is zero at the four punctures, and equals $\pm2\cat$ on the two lifts at each node. In rational Milnor $K_2$ of the function field over $\overline{\Q}$, bilinearity and antisymmetry give
\[
 \{M,L\}=2\{M,p\}=\{M^2,p\}=\{f,p\}
 =\{1+p,p\}-\{1-p,p\}=0
\]
by the Steinberg relations at $-p,p$ and torsion of symbols involving $-1$. By Guilloux--March\'e this is the criterion for an $A$-factor
of some 3-manifold with torus boundary~\cite[Def.~5.1 and Prop.~5.2]{GM}.
Section~\ref{sec:shape} identifies the particular realization
in the marked pattern.

Returning to the stationary action, on a critical patch $dS_{\rm crit}=\log L\,d\log x$. Hence $V_S=\Im S_{\rm crit}-\arg L\log|x|$ satisfies $dV_S=-\eta(x,L)=-2\eta(M,L)$. Fix the real basepoint $p=1/2$ to have $V_S=0$. Then $V_S=-2U(p)$, including its branch limits. The sign changes if the peripheral orientation is reversed.

\subsection{The integration chain and the known square-weave period}\label{sec:mahler}
The identity
\[
 \frac{E}{M^2L}=4-X_{\square}-X_{\square}^{-1}
                         -Y_{\square}-Y_{\square}^{-1},
 \qquad X_{\square}=M^2,\quad Y_{\square}=-L/M^2
\]
pulls back the square-weave polynomial under a finite torus cover and translation. On angles the integer matrix is $\left(\begin{smallmatrix}2&0\\-2&1\end{smallmatrix}\right)$, of determinant two. Both the resulting surjective torus map and the translation push probability Haar measure to itself. Thus its period agrees with \cite[Sec.~4.1, Thm.~12]{CKL}. The regulator primitive gives a direct evaluation once we identify the endpoints on the normalization. In angle coordinates for $(X_{\square},Y_{\square})$ the right side is $4-2\cos\alpha-2\cos\beta$, comparable near each zero to $\alpha^2+\beta^2$. Its logarithm is integrable there; the unit monomial and finite nonsingular torus cover preserve this property. Jensen's formula for monic $E$ in $L$ on $M=e^{\imunit\theta}$ chooses
\[
 L_{\rm out}=-M^2e^{h(\theta)},\qquad h(\theta)=\operatorname{arcosh}(2-\cos2\theta).
\]
The other root is inside the circle except at the endpoints. On both arcs $(0,\pi)$ and $(\pi,2\pi)$, oriented by increasing $\theta$, $p=(1-e^{h(\theta)})/(e^{2\imunit\theta}-1)$ starts at $\imunit$ and ends at $-\imunit$, since $h\sim2|\theta-\theta_0|$ near each endpoint. Each integral of $\eta(M,L_{\rm out})=dU$ is $-4\cat$. As $|M|=1$, Jensen integration gives $m(E)=-(2\pi)^{-1}\int\eta=4\cat/\pi$. This uses the displayed integral representative of $E$. In particular \(2\pi m(E)=8\cat\) in the square-weave normalization.

\subsection{The completed axis in the shape chart}\label{sec:shape}
We identify the marked geometric realization by imposing axis completeness
in Tillmann's right-handed link triangulation~\cite[Ch.~5, pp.~137--149, 161--162]{Tillmann}.
Write $(w,b_1,b_2,b_3)$ for the four nondegenerate shapes, $z'=1/(1-z)$, $z''=(z-1)/z$. The gluing and preferred peripheral equations are
\begin{align*}
 wb_1b_2b_3&=1,& wb_1(1-b_2)(1-b_3)&=(1-w)(1-b_1)b_2b_3,\\
 \bar s&=b_1'b_3''/b_2,& t_0&=b_1b_2,\\
 \bar u&=w'b_3''/b_2,& v&=wb_2.
\end{align*}
In the source marking the pairs are
\[
\begin{array}{ccl}
 \text{cusp }0 & (\bar s,t_0) & \text{axis},\\
 \text{cusp }1 & (\bar u,v) & \text{outer pattern component}.
\end{array}
\]
The bars denote squared meridian eigenvalues, and \(t_0,v\)
are preferred longitude eigenvalues. Completing the axis
means imposing complete cusp holonomy in this nondegenerate
chart. For the chosen lift these equations are
\((\bar s,t_0)=(1,-1)\). The axis remains a cusp of the
pattern exterior. Thus $b_2=-1/b_1$, then $b_1=b_3$, then $b_1=b_3=-1/w$, $b_2=w$. Conversely these values solve both gluing equations and both axis conditions for all $w\ne0,1,-1$. They give $u^2=\bar u=(1+w)/[w(1-w)]$, $v=w^2$. The axis meridian is nontrivially parabolic: its unipotent model has off-diagonal entry $-(1+w)/w$. Keeping an $\mathrm{SL}_2$ lift $u$, elimination gives $u^4v^2-(u^4-4u^2+1)v+1=0$, with inverse $w=(u^2v+1)/(u^2-1)$ away from the nodal fiber. At $u^2=1$, $w=\pm\imunit$. The preferred-word map in Section~\ref{sec:mark} identifies
Schilder's outer pair with the source pair of cusp \(1\),
and inverts both elements of the source axis pair. The mirror
in Proposition~\ref{prop:mark} takes the outer eigenvalues
\((u,v)\) to \((u^{-1},v)\) on a corresponding eigenline.
Passing to the other common eigenline inverts both eigenvalues.
With $p=1/y$ this gives the marked return map
\begin{equation}\label{eq:shapereturn}
 (M,L)=(u,v^{-1}),\qquad w=\frac{1-p}{1+p}=\frac{y-1}{y+1},\qquad
 p=\frac{1-w}{1+w}.
\end{equation} The factor selected by this marking is exactly $E$.

This coordinate correspondence also lets us compare the signed shape
volume with the stationary action.
The signed shape volume is $2D(w)+2D(-1/w)=2U(w)$. The relations $x(w)=x(p)$ and $L(p)=L(w)^{-1}$ in the two rational parametrizations give $d[U(w)+U(p)]=0$. Continuity on the connected parameter sphere and the zero at $p=1/2$ yield $2U(w)=-2U(p)=V_S$. At $w=\imunit$ all four ideal tetrahedra have shape $\imunit$,
with angles $\pi/2,\pi/4,\pi/4$. The volume refers to the
pattern piece with cusp \(0\) complete. The constants in the
two calculations are
\[
\begin{array}{rl}
 \text{signed volume in the right-handed chart at }w=\imunit & 4\cat,\\
 \text{signed value under the source mirror} & -4\cat,\\
 \text{unsigned volume of the two-cusped Whitehead link} & 4\cat,\\
 \text{square-weave }2\pi m(E) & 8\cat .
\end{array}
\]

\subsection{The branch selected on the physical diagonal}\label{sec:physical}
The two lifts over a plane node have the same plane coordinates.
The stationary calculation yields the candidate lifts.
The dominant-window and noncancellation estimates determine
which lift the auxiliary coordinates of the finite sum approach.
Fix an odd $d\ge3$ as $N\to\infty$ and use
\[
 q_N=e^{2\pi\imunit/N},\qquad t_N=-e^{\pi\imunit/(2N)},\qquad
 t_N^4=q_N,\qquad t_N^{2N}=-1 .
\]
Take \(q_N^{1/2}=e^{\pi\imunit/N}\).
This diagonal evaluation fixes \(M=-1\), while the AJ symbol
keeps \(M\) independent at \(t=-1\). Zheng's variable $t$ is our $q$, his torus pair is $(d,2)$, and the twist parameter of his clasp is zero. In the finite triangle
below, \(\ell\) is a summation index. Write $\widehat{\torcolor}_n^{(d)}=(q^{1/2}-q^{-1/2})\torcolor_n^{(d)}$. His generic reduced formula \cite[(2.9)]{Zheng} is
\begin{equation}\label{eq:genericdiag}
 F_N(q)=\frac{\sum_{n=0}^{N-1}\xi_{N,n}(q)\,\widehat{\torcolor}_n^{(d)}(q)}
                  {q^{N/2}-q^{-N/2}} ,
\end{equation}
where $\xi_{N,n}$ is his clasp coefficient, equivalently the $n$th coefficient of the finite sum in Section~\ref{sec:stationary}, including its common prefactor. This quotient defines the reduced polynomial before specialization. The denominator vanishes at $q_N$.

\emph{The derivative weight.} In the phase and norm conventions of \cite[(3.2)--(3.3)]{Zheng}, direct differentiation of \eqref{eq:genericdiag} gives
\begin{equation}\label{eq:physicaldiag}
 F_N(q_N)=\frac{q_N^{-1/2}}{N}\!
 \sum_{n=0}^{N-1}\sum_{\ell=0}^{N-1-n} A_{n,N}^{-1}S_{n,\ell}
 [b_{n,\ell}\widehat{\torcolor}_n^{(d)}+q\partial_q\widehat{\torcolor}_n^{(d)}]_{q_N}.
\end{equation}
Here $A_{n,N}=\exp(\pi\imunit n(n+1-N)/(2N))$ and
$S_{n,\ell}=\prod_{h=1}^n 4\sin^2(\pi(\ell+h)/N)/(2\sin(\pi h/N))$.
To check the signs in $b_{n,\ell}$, remove the common power $q^{N^2-(N+1)/2}$ from Zheng's (2.7). It multiplies a numerator that vanishes at $q_N$, so its derivative contributes zero to the sum. The remaining clasp summand and its logarithmic derivative are
\begin{align}
 H_{n,\ell}(q)&=q^{-N(\ell+n)}
 \prod_{h=1}^n\frac{(1-q^{N-\ell-h})(1-q^{\ell+h})}{1-q^h},\nonumber\\
 b_{n,\ell}&=\left.\frac{q H_{n,\ell}'(q)}{H_{n,\ell}(q)}\right|_{q_N}
 =-N(\ell+n)+\sum_{h=1}^n\left(
 \frac{N-\ell-h}{1-q_N^{-N+\ell+h}}+\frac{\ell+h}{1-q_N^{-\ell-h}}
 -\frac{h}{1-q_N^{-h}}\right).\label{eq:blog}
\end{align}
Indeed $q\partial_q\log(1-q^k)=k/(1-q^{-k})$; all exponents $k$ here lie strictly between $0$ and $N$. At $q_N$, $H_{n,\ell}=A_{n,N}^{-1}S_{n,\ell}$ and the common power equals $-q_N^{-1/2}$. The operator $q\partial_q$ applied to the denominator of \eqref{eq:genericdiag} gives $-N$ there. These evaluations give the sign and prefactor in \eqref{eq:physicaldiag}.

The three fraction signs in Zheng's printed (5.4) are opposite to those obtained from (2.7). Writing $b^{\rm Z}_{n,\ell}$ for that printed zero-twist weight gives the exact relation
\begin{equation}\label{eq:weightreturn}
 b_{n,\ell}=-b^{\rm Z}_{n,\ell}-2N(\ell+n).
\end{equation}
For instance, $N=3,n=1,\ell=0$ leaves $H=q^{-3}(1-q^2)$ and hence $b=-3+2/(1-q_3^{-2})$. We use \eqref{eq:blog} for the exact finite identity and the following separate analytic estimates for its limit.

\emph{Analytic input (Zheng, Sections 3--5).} Take $1/2<\delta<2/3$ and put $\mathcal U_N=\{\,|n-N/2|+|\ell-N/4|<N^\delta\,\}$ inside the triangle. For fixed $d$, this window eventually has $|2n+1-N|<N/(4d)$. Thus $\widehat{\torcolor}_n^{(d)}=O_d(1)$ there by \cite[Lemma~5.2, (5.12)]{Zheng}. The printed weight satisfies $b^{\rm Z}_{n,\ell}=O_d(N^2)$ there by (5.18), and $N(\ell+n)=O(N^2)$ throughout the triangle. Equation~\eqref{eq:weightreturn} therefore gives the same local bound for the corrected weight. With Zheng's independent derivative estimate (5.20) this yields
\begin{align}
 b_{n,\ell}\widehat{\torcolor}_n^{(d)}&=O_d(N^2),&
 q\partial_q\widehat{\torcolor}_n^{(d)}&=c_{d,N}v_n^{\mathrm{hom}}(q_N)+O_d(N^2),\label{eq:localder}\\
 c_{d,N}&=\frac{-2N^{5/2}}{(4d)^{3/2}}
 e^{-(d/2+2/d)\pi\imunit/(2N)+\pi\imunit/4} A^+_{d,2}(N,1).\nonumber
\end{align}
Here $A_{n,N}^{-8d}=v_n^{\mathrm{hom}}(q_N)$ exactly. The amplitude is the finite sum
\[
 A^+_{d,2}(N,1)=\sum_{j=1}^{2d-1}
 e^{-N j^2\pi\imunit/(4d)}j^2\sin(j\pi/d)\sin(j\pi/2).
\]
It is $8d$-periodic by \cite[(4.3)]{Zheng} and nonzero for every integer $N$ by Proposition~4.1 there. Hence $|c_{d,N}|\asymp_d N^{5/2}$.

Put $\beta=-(8d+1)$ and $S_0=S_{\lfloor N/2\rfloor,\lfloor N/4\rfloor}$. Lemmas~3.2--3.3 and (3.29) of \cite{Zheng} give
\begin{equation}\label{eq:gaussphysical}
 \sum_{\mathcal U_N} A_{n,N}^{\beta}S_{n,\ell}
 =G_\beta Ne^{-\beta N\pi\imunit/8} S_0(1+O_d(N^{3\delta-2})),
\end{equation}
where $G_\beta=e^{\beta\pi\imunit/4}\int_{\mathbb R^2}
 e^{-\pi(s^2+2sr+2r^2)+\beta\pi\imunit s^2/2}\,ds\,dr$. Integrating $r$ first leaves a Gaussian with positive real quadratic coefficient, so $G_\beta\ne0$. The $\beta=0$ estimate gives $\sum_{\mathcal U_N}S_{n,\ell}=O(NS_0)$. Before the exterior division by $N$, the bracket error in \eqref{eq:localder} therefore sums to $O_d(N^3S_0)$ against a main term of size $\asymp_d N^{7/2}S_0$. Globally (5.16) uses the individual bounds $b^{\rm Z}_{n,\ell}=O_d(N^3)$, $\widehat{\torcolor}_n^{(d)}=O_d(N)$ and $q\partial_q\widehat{\torcolor}_n^{(d)}=O_d(N^3)$. Equation~\eqref{eq:weightreturn} gives an $O_d(N^4)$ bound on the corrected bracket as well. The pointwise complementary estimate (3.25) and at most $N^2$ terms then give the bound of \cite[(5.21)]{Zheng}, namely $O_d(N^6 e^{-\epsilon N^{2\delta-1}}S_0)$ outside $\mathcal U_N$, with $\epsilon>0$. After the division in \eqref{eq:physicaldiag} the main size is $\asymp_d N^{5/2}S_0$, the additive local error is $O_d(N^2 S_0)$, and the complementary error is $O_d(N^5e^{-\epsilon N^{2\delta-1}}S_0)$. The relative Gaussian error in \eqref{eq:gaussphysical} tends to zero. These estimates establish noncancellation on the physical summation domain.

\emph{Identification of the lift.} In the coefficient convolution of Section~\ref{sec:stationary}, the first binomial index is precisely $\ell$ from \eqref{eq:physicaldiag}, and the other is $\ell^\prime=N-1-n-\ell$. Adjacent binomial ratios give the two dual coordinates
\begin{equation}\label{eq:dualratio}
 z_\ell=x\frac{1-q^{\ell+1}}{1-q^{n+\ell+1}},\qquad
 z_{\ell^\prime}=\frac{1-q^{\ell^\prime+1}}{1-q^{n+\ell^\prime+1}} .
\end{equation}
Uniformly on $\mathcal U_N$, $q_N^n\to-1$ and $q_N^\ell,q_N^{\ell^\prime}\to\imunit$. Both denominators tend to $1+\imunit$, so both coordinates tend to $-\imunit$ with error $O(N^{\delta-1})$. With the marked return map \eqref{eq:shapereturn}, the dominant indices therefore have the normalized limit
\begin{equation}\label{eq:branchreturn}
 (x,Q,z,y)\longrightarrow(1,-1,-\imunit,\imunit),\qquad
 p=1/y=-\imunit,\qquad w=(y-1)/(y+1)=\imunit .
\end{equation}
In the stated \(t_N\)-convention, \(M=-1\) and the shift formula
at this limit gives \(L=-1\). The node has two lifts. Equations~\eqref{eq:physicaldiag}--\eqref{eq:dualratio} select the one in \eqref{eq:branchreturn} for the fixed odd two-strand companion.

\section{Signed colors, coefficient descent, and primitive gauges}\label{app:domain}

The dimension gauge relates the order of the reduced tail
recurrence to that of the unreduced tail. For the exact signed
order we prove that the unreduced positive tail and its signed
extension have the same localized annihilator ideal.
Section~\ref{sec:signedreturn} uses a weaker transfer for factor
survival. We write
\(J_K(n;t)=\Phi_K(n;t)\) on all integer colors and
\(j_K(n;q)=j_{K,n}(q)\) on positive colors, with the normalization
of Section~\ref{sec:comparison}.

The tail argument interpolates the color coefficients in a formal
expansion at \(t=1\). The AJ specialization of the resulting
primitive row is at \(t=-1\).

\begin{lemma}[Tail to all integers]\label{lem:tail}
If a Laurent-polynomial coefficient operator annihilates $J_K(n;t)$
for every sufficiently large positive integer $n$, it annihilates
the prescribed signed unreduced sequence at every integer.
\end{lemma}
\begin{proof}
The cyclotomic expansion \cite[Theorems 4.5 and 6.4,
(6.5), (7.1)--(7.2)]{Habiro}, in the product form
\cite[Section 3.2]{GLloop}, gives for positive dimension colors
\[
 j_K(n;q)=\sum_{k\ge0}H_k(q)
       \prod_{a=1}^k(q^n+q^{-n}-q^a-q^{-a}),\qquad
 H_k\in\Z[q^{\pm1}],\quad q=t^4.
\]
The sum terminates at $k=n-1$. For a formal color $X$ put
$t=e^\hbar$, $q=e^{4\hbar}$. Each factor of the product belongs to
$\hbar^2\Q[X^2][[\hbar]]$, so each coefficient of the infinite
expansion is a finite sum in $\Q[X^2]$.
Multiply by
\[
 \frac{\sinh(2X\hbar)}{\sinh(2\hbar)}
       \in X\Q[X^2][[\hbar]] .
\]
This gives one series $\widehat J(X,\hbar)$ with odd polynomial
color coefficients, agreeing with the Jones expansion at
every positive and negative color and giving zero at color
zero.

For $P=\sum_i a_i(t,M)L^i$ the formal residual
\[
 \sum_i a_i(e^\hbar,e^{2X\hbar})\widehat J(X+i,\hbar)
\]
has polynomial coefficients in $X$. They vanish on infinitely
many positive integers and hence identically. At any integer
color the Laurent polynomial $(PJ)(n;t)$ thus has zero Taylor
series under $t=e^\hbar$. This substitution is injective because
the Vandermonde determinant detects the coefficients of a finite
linear combination of distinct exponentials from its derivatives.
Therefore every color is annihilated. The same proof works
for any infinite set of positive test colors.
\end{proof}

For rational coefficients clear a common denominator on the
left. A nonzero Laurent polynomial in $t,M$ can vanish identically
on only finitely many substitutions $M=t^{2n}$:
distinct finite exponent pairs eventually have distinct exponents.
Thus rational tail and all-integer localized ideals coincide.
In particular their minimal Laurent shift spans agree.

\begin{lemma}[Descent through the meridian square]\label{lem:descent}
For the reduced positive-color sequence, the monic minimum over
$\Q(t,M)$ descends to $\Q(q,x)$, where $q=t^4$ and $x=M^2$.
Its order is unchanged.
\end{lemma}
\begin{proof}
A normalized least-span recurrence is unique because subtracting
two monic recurrences of the same least degree would give a
smaller one. Over the constants extended by $i$, the $1,i$
coefficient parts annihilate the \(\Q(t)\)-valued sequence separately,
so the least order is unchanged.
Write its coefficients as $a_j(t,M)$.

Substitute $t\mapsto it$ in the relation evaluated at
each fixed positive color $n$. All the reduced invariant values
in the relation are unchanged because they are Laurent in $t^4$.
The coefficient's meridian evaluation becomes
$(it)^{2n}=(-1)^nt^{2n}$. The coefficients $a_j(it,M)$
therefore give a monic relation on even positive colors, and
$a_j(it,-M)$ on odd positive colors, both with the original
color shift. The interpolation proof of Lemma~\ref{lem:tail},
without its dimension factor, extends each of these relations
to the positive tail. Uniqueness of the monic minimum there
gives $a_j(t,M)=a_j(it,M)=a_j(it,-M)$.
The joint fixed field in $\Q(i,t,M)$ is
$\Q(i,t^4,M^2)$. Intersect with the original rational function field to
obtain $\Q(t^4,M^2)$. On this subfield the actual shift is
$x\mapsto qx$, as required.
\end{proof}

The unreduced dimension gauge is
$g_{\rm dim}=(M-M^{-1})/(t^2-t^{-2})$.
Conjugating an order-$j$ term multiplies its coefficient by
$g_{\rm dim}/\tau^j g_{\rm dim}$, whose residue at $t=-1$ is one with $M$
independent. Also $q-1=t^4-1$ has a simple zero there.
Consequently descent, extension, and this gauge preserve the
primitive residue up to the declared meridian units.

\section{Exact arithmetic checks}\label{app:finite}

For the nonzero forcing use arithmetic modulo $1000003$,
$q=4$, $N=27$, $x=448586$. Compute $\mathcal T_n$ from the centered torus
sum used in the growth argument. Use \eqref{eq:Bailey} to obtain $D_j$,
padding negative
$j$ by zero. Apply the polynomial
$S_5(\mathcal V)+(1-qy)S_6(\mathcal V)-qyS_7(\mathcal V)$.
Its six initial values are
\[
           (999997,885180,708642,458952,782718,585014).
\]
All factorial and weight denominators are checked nonzero by
inversion in the field. Their six-term weighted sum in
\eqref{eq:PF} is $31260$. This suffices to prove that the
rational function $b$ is nonzero. Equation \eqref{eq:PF} follows
from the finite transform and adjoint identity.
For \eqref{eq:Bmod}, repeated seventh powering modulo the
displayed degree-twelve polynomial gives the two gcd-one
conditions at exponents $4,6$ and the identity at exponent
$12$. These are the irreducibility predicates because
the prime divisors of $12$ are $2,3$.
The first Ore jets use rational arithmetic in
$\Q[x]/(78x^2+157x+78)$ and at $x=2/3,3/2$,
exactly as in Section~\ref{sec:whitehead}.

\subsection{Computational attachments}
The scripts are available as arXiv ancillary files. With Python and
SymPy installed, run \verb|python3 reproduce.py --out replay| in the
directory containing the downloaded scripts. Source supplement S1
(Appendix~\ref{app:verification}) contains the prior run records.

\begin{description}
\item[\texttt{bounded\_forcing\_check.py}]
Evaluates the six defects and the weighted forcing in
$\mathbb F_{1000003}$. It checks denominators by modular inversion.
\item[\texttt{check\_symbol.py}]
Checks the symbol identities, degree and endpoint coefficients of $B$,
and all three Frobenius predicates for \eqref{eq:Bmod}.
\item[\texttt{first\_ore\_index.py}]
Computes both allocations at the $B/E$ crossing and the abelian
crossing check at $x=2/3$ and $3/2$.
\end{description}

The forcing check uses the Python standard library. The symbol and
crossing checks use SymPy for exact polynomial and rational arithmetic.
In these scripts, the transform coordinate $y$ is written
\texttt{J} or \texttt{j}. The crossing script writes $\Phi_6$ as
\texttt{F}.

The script \texttt{family\_jet.py} checks the symbolic $E/B_m$ jet in
Section~\ref{sec:jets}. The script \texttt{family\_words.py} checks the
generic word and matrix identities of Section~\ref{sec:word}.
The resultant and inverse map in Section~\ref{sec:stationary} can be checked
with \texttt{stationary\_curve.py}.
The diagram identification, the local lattice theorem and the analytic
estimates have separate proofs or sources. Appendix~\ref{app:verification}
records the scope of the formal results and their interfaces.

\section{What has been formally checked}\label{app:verification}

The written proof of Theorem~\ref{thm:four} covers both clasp signs,
both companion chiralities, and every integer twist.
The current family formalization treats the explicit braids for the
untwisted negative-clasp doubles of right-handed torus knots.
The new classical divisibility and invariant-ideal consequences
of Section~\ref{sec:conclusion} are written deductions.
The historical comparator acceptance concerns one fixed braid.
Jinzheng Li used \href{https://github.com/MerLeanProver/MerLean}{MerLean}
\cite{ren2026merleanagenticframeworkautoformalization,li2026merleanproverrecursiveloopingharness}
in the formalization work.

\subsection{The untwisted family library}

For each coprime \(2\le r<s\), Jinzheng Li's \texttt{theoremT} at commit
\texttt{78e38735} uses the same explicit braid for both invariants of
the untwisted negative-clasp double of the right-handed \(T(r,s)\).
It gives a primitive minimal recurrence and an \(A\)-polynomial in the
library's predicates. Every nonzero primitive saturated row has nonzero
specialization divisible by \(E\). Every \(A\)-polynomial in the
predicate is coprime to \(E\) in the comparison ring, so radical equality
fails. Theorem \texttt{theoremF} gives the corresponding result for
\(r=2,\ s=2k+1\), \(k\ge1\).

The companion declarations have two sharper conclusions. The declaration
\nolinkurl{release_classicalEquality} gives \(A\doteq(L-1)B_{rs}\) in
\(\Q(M)[L^{\pm1}]\) for the untwisted braid family. In the two-strand
subfamily, \nolinkurl{release_minimalRecurrence} gives degree \(8k+8\)
for every primitive minimal recurrence and specialization associated to
\((L-1)^2E B_{2(2k+1)}\). For \(r\ge3\) the recurrence conclusion here is
existence and factor survival. These are library results for the specified
braid records. The written family statement of Theorem~\ref{thm:family}
uses classical containment, and Appendix~\ref{sec:whitehead} gives the
exact written result for the untwisted trefoil.

The records at the family commit give a local Lean build, per-module
kernel replay and axiom checks using only \texttt{propext},
\texttt{Quot.sound} and \texttt{Classical.choice} or subsets of them.
The separate release comparator configuration has no recorded
acceptance. The November~15 comparator result covers the fixed braid
at the historical commit.

\subsection{The historical single braid}

The November~15 comparator matched the exported statement to the challenge.
Lean's kernel and the independent checker nanoda accepted the proof of
\nolinkurl{AJConjecture.not_ajSquareFreeConjecture} at solution commit
\texttt{8dd21cf}~\cite{Comparator}. It negates universal radical equality
in the library's comparison ring by giving recurrence and boundary
polynomial witnesses for a fixed seven-strand braid.

The identification with the named trefoil double has separate topological
evidence. An oriented complement comparison with the zero-framed reference
diagram uses an explicit face gluing isomorphism. Gordon--Luecke
\cite{GordonLuecke} gives the return from the oriented exterior to the
knot type. The meridian slope and the primitive kernel of boundary
homology recover the peripheral pair up to simultaneous inversion.
Proposition~\ref{prop:mark} fixes the pattern for arbitrary companions.

\subsection{Geometric and analytic interfaces}

The library reaches the diagram \texttt{whDouble} from the cut-open braid
by formal moves and checks zero band framing, negative clasp and one
component. It gives a homomorphism from the companion group to the
braid's knot group with the pattern peripheral relations, and lifts
\(\SL(\C)\) representations with the prescribed boundary matrices.
The geometric satellite identification uses the moves as isotopies and
the marking of Section~\ref{sec:mark}. Interpreting the braid Jones function
and group presentation as the conventional invariant and the knot group
with its preferred peripheral data is a separate definition interface.

For its explicit braids the library constructs witnesses to
\texttt{IsAPolynomialOf}. This predicate requires an integral equation
with no repeated factors over \(\C\). Its factors cut out exactly the
irreducible affine curves in the closure of the actual eigenvalue set. The conventional
identification uses affine component selection and the peripheral marking.
Lemma~\ref{lem:aex} gives the written horizontal comparison off the
exceptional meridian fibers. Formal distinctness concerns Markov classes
and color-two Jones values. Its geometric interpretation uses Markov's
theorem or isotopy invariance, and the formal JSJ route takes a separate
topological premise.

The family library also checks parts of the quotient argument and of
the geometry of \(E\). Its coordinate laws cover monomial meridian
substitutions, ramified meridian pullbacks, and linear rescalings of the
deformation parameter. The formal node tests for the untwisted family
give endpoint conditions and indices. The written proof also checks
joint regularity of the companion matrix and its inverse. The family scope
record in S1 gives the declarations and their scope
(\nolinkurl{support/family/V1_SCOPE.md}).
The geometric interpretation of the shape
equations and the uniqueness of the smooth projective model use external
geometric results. For the fixed odd two-strand diagonal, the exact
finite-color identity is checked in the library. The formal asymptotic
statement uses \nolinkurl{ZhengVerbatim}, which consists of Zheng's
Lemmas~3.2--3.3, (5.12) at odd colors, (5.20), and
Proposition~4.1~\cite{Zheng}.
The written proof in Section~\ref{sec:physical} identifies the analytic
inputs where they are used. These estimates are independent inputs to
the diagonal analysis and are unnecessary for the core AJ comparison.

\subsection{Review evidence and source access}\label{sec:source-supplement}

The written proof carries range-specific review qualifications, with
the reflected cases inheriting those of their source ranges.
The scope, qualifications and retained evidence are documented in
\nolinkurl{SOURCE_AND_ASSURANCE.md}.

\emph{Source supplement S1} is
\nolinkurl{Whitehead-AJ-four-families-r6-review-sources.zip}.
It includes the r6 review manuscript, the Appendix~\ref{app:finite} checks,
and the family scope record.

The S1 comparator report is
\nolinkurl{support/single-braid/comparator-2026-11-15.md}.

\emph{Source supplement S0} is
\nolinkurl{Whitehead-AJ-four-families-r4-sources.zip},
the archived \(43\)-page baseline.

A public release of the Lean sources and check records is planned for the
\href{https://github.com/DavetheSlayer/aj-conjecture-formalization}{formalization repository}.
A public locator for S1 is pending.

\bibliographystyle{amsplain}
\bibliography{references}
\end{document}